\documentclass[a4paper,reqno]{amsart}
\usepackage{amsmath}
\usepackage{amssymb}
\usepackage{tikz-cd}
\usepackage{color}
\usepackage{verbatim}

\usepackage{amscd}
\usepackage{hyperref}

\definecolor{Blue}{rgb}{0.3,0.3,0.9}
\definecolor{Red}{rgb}{0.9,0.3,0.3}

\usepackage[nobysame]{amsrefs}

\newtheorem{Lemma}{Lemma}[section]
\newtheorem{Th}[Lemma]{Theorem}
\newtheorem{PDTh}[Lemma]{Package Deal Theorem}
\newtheorem{Prop}[Lemma]{Proposition}

\newtheorem{Cor}[Lemma]{Corollary}

\theoremstyle{definition}

\newtheorem{Def}[Lemma]{Definition}
\newtheorem{Ex}[Lemma]{Example}
\newtheorem{Exs}[Lemma]{Examples}
\newtheorem{notation}[Lemma]{Notation}

\theoremstyle{remark}
\newtheorem{Rem}[Lemma]{Remark}

\newtheorem{Remark}[Lemma]{Remark}
\newtheorem{Remarks}[Lemma]{Remarks}
\newenvironment{Proof}{{\sc Proof.}\ }{~\rule{1ex}{1ex}\vspace{0.5truecm}}

\newcommand{\add}{\mbox{\rm add}}
\newcommand{\Add}{\mbox{\rm Add}}

\newcommand{\im}{\mbox{\rm Im}}

\newcommand{\N}{\mathbb N}
\newcommand{\No}{{\mathbb N}_0}
\newcommand{\Z}{\mathbb{Z}}

\newcommand{\fm}{\mathfrak{m}}
\newcommand{\fn}{\mathfrak{n}}
\newcommand{\fp}{\mathfrak{p}}

\newcommand{\Mod}{\mbox{\rm Mod-}}

\DeclareMathOperator{\mSpec}{mSpec}

\DeclareMathOperator{\Ker}{Ker}
\DeclareMathOperator{\Hom}{Hom}
\DeclareMathOperator{\End}{End}
\DeclareMathOperator{\Tr}{Tr}
\DeclareMathOperator{\rank}{rank}

\title[Monoids of infinitely generated modules] 
{Relatively big projective modules and their applications to direct sum decompositions     
}

\dedicatory{Dedicated to Tariq Rizvi}

\begin{document}

\author{Rom\'an \'Alvarez}

\address{Departament de Matem\`atiques 
Universitat Aut\`onoma de Barcelona, 08193 Bellaterra
(Barcelona), Spain\newline 
Charles University, Faculty of Mathematics and Physics \\Department
of Algebra, Sokolovsk\'a~83,
18675 Praha 8, Czech Republic}
\email{roman9496@gmail.com}

\author{Dolors Herbera}

\address{Departament de Matem\`atiques 
Universitat Aut\`onoma de Barcelona, 08193 Bellaterra
(Barcelona), Spain }
\email{dolors.herbera@uab.cat}

\thanks{ The first and second authors were partially supported by the projects MIMECO  PID2023-147110NB-I00 financed by the Spanish Government and  \emph{Laboratori d'Interaccions entre Geometria, \`Algebra i Topologia} (LIGAT) with reference number 2021 SGR 01015 financed by the Generalitat de Catalunya.}
 \author{Pavel P\v r\'\i hoda}
\address{Charles University, Faculty of Mathematics and Physics \\Department
of Algebra, Sokolovsk\'a~83,
18675 Praha 8, Czech Republic}
\email{prihoda@karlin.mff.cuni.cz}

\thanks{The third author was supported by Czech Science Foundation grant GA\v CR 23-05148S}

\date{\today}

\begin{abstract}
Countably generated projective modules that are relatively big with respect to a trace ideal were introduced by P. P\v ríhoda, as an extension of Bass' uniformly big projectives. It has already been proved that there are a number of interesting examples of rings whose countably generated projective modules are always relatively big. In this paper, we increase the list of such examples, showing that it includes all right noetherian rings satisfying a polynomial identity.

We also show that countably generated projective modules over  locally semiperfect torsion-free algebras over $h$-local domains are always relatively big. This last result applies to endomorphism rings of finitely generated torsion-free modules over $h$-local domains. As a consequence, we can give a complete characterization of those $h$-local domains of Krull dimension $1$ for which every direct summand of a direct sum of copies of a single finitely generated torsion-free module is again a direct sum of finitely generated modules.
\end{abstract}

\maketitle

\tableofcontents


\section*{Introduction}

For a long time, the study of projective modules was focused on the finitely generated ones. Outside the finitely generated setting,  the classical theorem of Kaplansky \cite{kaplansky} stating that every projective module is a direct sum of countably generated projective modules shows that one can just consider the countably generated case. Interesting monoids associated to  infinitely generated projective modules have recently been studied: the monoid $V(\Lambda)\sqcup B(\Lambda)$ of relatively big projectives defined, formally, in \cite{wiegand} and that is the main object of interest in this paper, but also the analogous of the Cuntz monoid for general rings defined in \cite{AABPV}, and the $\kappa$-monoids defined in \cite{nazemian_smertnig}.

Let $\Lambda$ be a ring, and let $V(\Lambda)$ (resp. $V^*(\Lambda)$) denote the set of representatives of the isomorphism classes of finitely generated (resp. countably generated) projective right $\Lambda$-modules. These sets form commutative monoids with addition induced by the direct sum of modules, with $V(\Lambda)$ naturally embedded in $V^*(\Lambda)$.

The concept of relatively big projective module was introduced by P\v r{\'\i}hoda in \cite{P2}, extending Bass' uniformly big projectives \cite{bass}. Let $J$ be the trace of a countably generated projective module. A countably generated projective right $\Lambda$-module $P$ is said to be $J$-big,  if every countably generated projective $\Lambda$-module whose trace is contained in $J$ is a homomorphic image of $P$.  A countably generated projective module $P$  is said to be \emph{relatively $J$-big}  if  $P$ is  $J$-big  and $P/PJ$ is finitely generated. Therefore, $J$ is maximal with respect to the property that $P$ is $J$-big and minimal with respect to the property that $P/PJ$ is finitely generated (cf. \S~\ref{subsec:fg}). A crucial property is that a relatively $J$-big projective module $P$ is uniquely determined, up to isomorphism, by the trace ideal $J$ and the isomorphism class of $P/PJ$ (cf. Lemma~\ref{Jbigdetermined}). 

The relatively big projectives are the countably generated projective modules that are relatively $J$-big for some $J$.  The important fact that boils down in the concept of relatively big projective is that projective modules  modulo a trace ideal can be lifted to projective modules over the ring. If the projective module over the quotient is finitely generated, then there is a unique  lifting to a relatively big projective module over the ring.

The monoid $V(\Lambda)\sqcup B(\Lambda)$ of isomorphism classes of relatively big projective modules was introduced in \cite{wiegand} to capture the properties of such projectives. This monoid embeds in $V^*(\Lambda)$ and can be thought of as the \emph{disjoint union} of $V(\Lambda/J)$, where $J$ runs through the trace ideals of the countably generated $\Lambda$-modules, see Corollary~\ref{iso.monoid}. 

In \cite{P2}, it was proved that over a noetherian ring $\Lambda$ which satisfies a certain descending condition on two-sided ideals, which we call the $(*)$-condition,  $V(\Lambda)\sqcup B(\Lambda)=V^*(\Lambda)$ (cf. \cite[Lemma 2.8]{P2} and Theorem~\ref{th:star}). Commutative noetherian rings, finitely generated algebras over commutative noetherian rings of Krull dimension one, and noetherian semilocal rings  satisfy the $(*)$-condition. Therefore, all countably generated projective modules over these rings are relatively big. These results were the main ingredients used to characterize $V^*(\Lambda)$ for $\Lambda$ a noetherian semilocal ring in \cite{crelle}, or to study some infinite direct sum decompositions in \cite{wiegand}.

The main results of this paper  provide new classes of rings for which $V(\Lambda)\sqcup B(\Lambda)=V^*(\Lambda)$. We prove that  right noetherian PI-rings satisfy the $(*)$-condition (see Theorem~\ref{starpi}), and this is one of the main ingredients in showing that all countably generated right (or left) projective modules over a right noetherian PI-ring  are relatively big in Corollary~\ref{PImonoid}.  In particular, this applies to module-finite algebras over commutative noetherian rings, such as group algebras of finite groups over commutative noetherian rings, or to endomorphism rings of finitely generated modules over commutative noetherian rings. This shows that the restriction on Krull dimension that we  were imposing before was superfluous.

Beyond the noetherian setting, it is easy to show that for semiperfect rings, all countably generated projective right (or left) modules are relatively big (cf. Proposition~\ref{monoidsemiperfect}). In Theorem~\ref{projectivesls}, we show that the same conclusion holds  for $\Lambda$ a locally semiperfect torsion-free algebra, over an $h$-local commutative domain $R$ with field of fractions $Q$, such that $\Lambda \otimes _R Q$ is simple artinian. 

The $h$-local domains  satisfy strong noetherianity conditions on the maximal spectrum. More precisely,  any non-zero ideal of $R$ is contained only in a finite number of maximal ideals, and any non-zero prime ideal of $R$ is contained in a unique maximal ideal. Noetherian domains of Krull dimension one are always $h$-local. The definition of $h$-local domain is attributed to Matlis. However, Jaffard in \cite{jaffard} introduced Dedekind-type domains as those domains $R$ such that any non-zero ideal is a finite product of ideals that are contained in a unique maximal ideal of $R$. Jaffard proved that Dedekind-type domains coincide with $h$-local domains, cf. \cite[Theorem~2.1.5]{FHL}.

The class of  $h$-local domains is ubiquitously associated with many classical problems in the study of modules over non-noetherian domains, see \cite{matlis3} and the recent survey \cite{olberding}. 
Non-finitely generated projective modules over an $h$-local domain are free \cite{hinohara}. Thus, what is of interest is the description of projective modules over algebras on such domains. 

Our interest in this setting comes from our paper \cite{AHP}, in which we studied domains for which the class of direct sums of finitely generated torsion-free modules is closed under direct summands.  If $\Lambda$ is  the endomorphism ring of a finitely generated torsion-free module $M$ over an $h$-local domain $R$, then $\Lambda$ is torsion-free and  $\Lambda\otimes _R Q$ is simple artinian. If, in addition, $R$ has Krull dimension $1$ and all the direct summands of a direct sum  of copies of $R\oplus M$ are again a direct sum of finitely generated modules, it was proved in \cite{AHP} that $\Lambda$ is locally semiperfect. Therefore, our results show that over such $\Lambda$ all countably generated projective modules are relatively big.

It is interesting to recall that Procesi and Small \cite{PS} proved  that the endomorphism ring of a finitely generated module over a commutative ring or, in general, over a ring satisfying a polynomial identity, is also a ring satisfying a polynomial identity. As far as we know, this point of view has not been so much examined in the world of finitely generated modules, endomorphism rings and direct sum decompositions, and certainly, it seems it should be worthwhile to explore it.

The paper is structured as follows. In Section~\ref{sec:1} we give  the basic definitions and notation. We introduce the needed background, presenting the concept of relatively big projective module, the $(*)$-condition and the main results relating both notions. In Section~\ref{sec:2} we examine the monoid of isomorphism classes of relatively big projective modules, its algebraic properties and explore the right/left (non)symmetry of the monoids; check Example~\ref{ex:asymmetry} and Corollary~\ref{isomonoidtraces}. 
In Section~\ref{sec:3} we show that countably generated, projective right  modules lift modulo a right $T$-nilpotent ideal, cf. Lemma~\ref{liftingproj}. We prove that for a right $T$-nilpotent ideal $N$ of $\Lambda$ there are monoid isomorphisms between $V^*(\Lambda)$ and $V^*(\Lambda/N)$. This isomorphism restricts to isomorphisms between $V(\Lambda)$ and $V(\Lambda/N)$, and between $V(\Lambda)\sqcup B(\Lambda)$  and $V(\Lambda /N)\sqcup B(\Lambda/N)$, cf. Corollary~\ref{isomonoidnil}. In Section~\ref{sec:4},  these results  are an important ingredient to prove that right noetherian PI-rings satisfy the $(*)$-condition (Theorem~\ref{starpi}), and to conclude that $V(\Lambda)\sqcup B(\Lambda)=V^*(\Lambda)$ for $\Lambda$ a right noetherian PI-ring in Corollary~\ref{PImonoid}.  In particular, this holds for module-finite algebras over commutative noetherian rings (cf. Corollary~\ref{modulefinite}). 

In the rest of the paper, we turn to algebras over  $h$-local domains. Section~\ref{sec:5} discusses Package Deal Theorems for projective modules over $\Lambda$, whenever $\Lambda$ is an algebra over an $h$-local domain; our goal is twofold: recalling the results  from \cite{AHP} and proving some extensions needed in the subsequent sections. Section~\ref{sec:6} focuses on semiperfect and locally semiperfect algebras. We observe in Proposition~\ref{monoidsemiperfect}   that  for semiperfect rings, all countably generated projective modules are relatively big.  

Let $R$ be an $h$-local domain with field of fractions $Q$. In the remaining sections, we study locally semiperfect algebras $\Lambda$  over  $R$ such that $\Lambda \otimes _R Q$ is simple artinian. In Section~\ref{sec:6}, we show that  for such algebras, every trace of a projective module is the trace of a finitely generated projective module (cf. Lemma~\ref{semiperf.fgtrace}) and that all countably generated projective $\Lambda$-modules are relatively big (cf. Theorem~\ref{projectivesls}).

In Section~\ref{sec:7}, we study in detail the genus of countably generated projective modules over such  algebras, that is, the isomorphism classes of countably generated projective modules over the localizations at the maximal ideals of $R$.  This is a submonoid of the product of $V^*(\Lambda _\fm)$, for $\fm$ maximal ideal of $R$, that can be described quite explicitly in terms of solutions in $\No \cup \{\infty \}$ of a system of linear equations with coefficients in $\No$. 

In Proposition~\ref{genusisomorphic}, we show that two countably generated projective modules that are not finitely generated are isomorphic if and only if they have the same genus, that is, whenever all their localizations at maximal ideals of $R$ are isomorphic.  This is also true for finitely generated projective modules, provided $\Lambda$ is a finitely generated algebra over a semilocal domain $R$, cf. Proposition~\ref{genusisomorphic_semilocal}.

These tools allow us to study in  Section~\ref{sec:8}  when every projective module is a direct sum of finitely generated modules. Finally, in Section~\ref{sec:9}, we apply our results to $\Lambda$ the endomorphism ring of a finitely generated torsion-free module $M$ over an $h$-local domain $R$ of Krull dimension $1$.  In Theorem~\ref{th:Add}, we give necessary and sufficient conditions that imply that any direct summand of an arbitrary direct sum of copies of $M$ is a sum of finitely generated modules. These results complement and clarify the ones obtained in \cite{AHP}.

\bigskip

Throughout the paper, rings are associative with $1$, and morphisms are unital. In each section, we try to be very precise about the setting we are working with. We mainly use $R$ to denote a commutative ring and $\Lambda$ to denote a non-necessarily commutative ring or an algebra over a domain $R$.

The set $\N =\{1,2,\dots \}$ and $\No =\N \cup \{0\}$.

\bigskip

We thank the reviewer for the careful reading of the article, and for all the helpful  comments and suggestions.

\section{The monoid of infinitely generated projective modules}\label{sec:1}

In this section we will cover the basic notions on projective modules and their trace ideals. We start by introducing the notation already used in \cite{wiegand} for monoids of projective modules.

\begin{notation} 
For a given ring $\Lambda$, we write $V_r(\Lambda)$ (resp. $V_r^*(\Lambda)$) for the set of representatives of the isomorphism classes of finitely generated (resp. countably generated) projective right $\Lambda$-modules. If $P$ is a finitely (resp. countably) generated projective right $\Lambda$-module, we denote its representative in $V_r(\Lambda)$ (resp. in $V_r^*(\Lambda)$) by $[P]$.

Similarly, we write $V_\ell(\Lambda)$ (resp. $V_\ell^*(\Lambda)$) for the corresponding sets of representatives of the isomorphism classes of projective left $\Lambda$-modules and follow the same notation for the representatives of the isomorphism classes in these sets.
\end{notation}

{Unless stated otherwise, our modules will always be right modules, so we will write $V(\Lambda)$ and $V^*(\Lambda)$ instead of $V_r(\Lambda)$ and $V_r^*(\Lambda)$, respectively, and leave the previous notation in cases where both right and left modules appear. All the results stated for $V(\Lambda)$ and $V^*(\Lambda)$ can be stated in the context of  left $\Lambda$-modules for $V_\ell(\Lambda)$ and $V_\ell^*(\Lambda)$, respectively.}

\begin{Remark}
Note that both $V(\Lambda)$ and $V^*(\Lambda)$ are commutative monoids with addition defined by $[P]+[Q]=[P\oplus Q]$. Note that $V(\Lambda)$ is a submonoid of $V^*(\Lambda)$ and that $V^*(\Lambda)\setminus V(\Lambda)$ is a subsemigroup of $V^*(\Lambda)$.
\end{Remark}

An important object related to projective modules is the set of their trace ideals. We recall the definition and some basic properties.

\begin{Def}
Let $\Lambda$ be a ring, and let $M$ be a right $\Lambda$-module. Then the \textit{trace} of $M$ in $\Lambda$ is the two-sided ideal of $\Lambda$ defined as 
    $$\Tr_\Lambda(M)=\sum _{f\in \Hom_\Lambda(M,\Lambda)} f(M).$$
    
We will denote by $\mathcal T_r(\Lambda)$ (resp. $\mathcal T_\ell(\Lambda)$) the set of ideals of $\Lambda$ that are traces of countably generated projective right (resp. left) $\Lambda$-modules. When we simply write $\mathcal T(\Lambda)$ we mean $\mathcal T_r(\Lambda)$.
\end{Def}

\begin{Remark}\label{rem:tensortrace}
Let $\Lambda$ be a ring and let $M$ be a right $\Lambda$-module with endomorphism ring $S=\End_\Lambda(M)$. Then $M$ is also a left $S$-module. Since $\Hom_\Lambda(M,\Lambda)$ is both a right $S$-module and a left $\Lambda$-module, the tensor product $\Hom_\Lambda(M,\Lambda)\otimes_S M$ is well-defined. Then $\Tr_\Lambda(M)$ is the image of the \textit{trace map}:
    \[\begin{tikzcd}[row sep=0ex]
        \tau\colon & [-3em] \Hom_\Lambda(M,\Lambda)\otimes_S M\rar & \Lambda \\
        & f\otimes x\rar[mapsto] & f(x).
    \end{tikzcd}\]
\end{Remark}

\begin{Lemma}
\label{ds.traces}
Let $\Lambda$ be a ring. 
\begin{enumerate}
    \item[(i)] Let $M,N$ be right $\Lambda$-modules. If $M$ is a homomorphic image of $N$, then $\Tr_\Lambda(M)\subseteq \Tr_\Lambda(N)$.
    \item[(ii)] Let $\{M_i\}_{i\in A}$ be a family of right $\Lambda$-modules. Then $$\Tr_\Lambda\Big(\bigoplus_{i\in A} M_i\Big)=\sum _{i\in A}\Tr_\Lambda(M_i).$$
    \item[(iii)] \cite[Proposition 2.40]{lam2} Let $J$ be an ideal of $\Lambda$, and let $P$ be a projective right $\Lambda$-module. Then $\Tr_\Lambda(P)\subseteq J$ if and only if $PJ=P$. In addition, the trace of a projective module is always an idempotent ideal.
\end{enumerate}
\end{Lemma}

We recall the following properties of traces of projective modules.

\begin{Lemma}
\label{tracaquocient} 
Let $\Lambda$ be a ring, let $P$ be a right $\Lambda$-module, and let $J$ be an ideal of $\Lambda$. Then $(\Tr _\Lambda (P)+J)/J \subseteq \Tr _{\Lambda /J}(P/PJ)$. Moreover, if $P$ is projective, then $(\Tr _\Lambda (P)+J)/J = \Tr _{\Lambda /J}(P/PJ)$.
\end{Lemma}

\begin{Remark} \label{chartraces} 
For a commutative ring $\Lambda$, Vasconcelos showed that trace ideals of projective modules are exactly the pure ideals (see page 430 in \cite{Vasconcelos2} or Proposition~2.12 in \cite{traces}).

In \cite{whitehead}, Whitehead gave a characterization of ideals that are traces of projective modules for general rings. In particular, he showed that an idempotent ideal of $\Lambda$ that is finitely generated as a left ideal is always the trace of a countably generated projective right $\Lambda$-module. We shall use this result throughout the paper, sometimes without prior acknowledgment.

An alternative approach to the characterization of the traces of projective modules can  be found in the paper \cite{traces}.
\end{Remark}

\begin{Lemma} \label{tracesmodjradical}  
Let $\Lambda $ be a ring with Jacobson radical $J(\Lambda)$, and let $I$ and $I'$ be elements in $\mathcal T(\Lambda)$. Then
\begin{enumerate}
    \item[(i)] $I=I'$ if and only if $I+J(\Lambda)=I'+J(\Lambda)$.
    \item[(ii)]  $(I+J(\Lambda))/J(\Lambda) \in \mathcal T(\Lambda/J(\Lambda))$, so the assignment defined by $I \mapsto (I+J(\Lambda))/J(\Lambda)$ defines an injective map $\varphi\colon \mathcal T(\Lambda)\longrightarrow \mathcal T(\Lambda/J(\Lambda))$.
\end{enumerate}
\end{Lemma}

\begin{Proof}
$(i)$ is included in \cite[Corollary~2.9]{traces}. 

To prove $(ii)$, note that if $I$ is the trace of the projective right $\Lambda$-module $P$ then, by Lemma~\ref{tracaquocient}, $(I+J(\Lambda))/J(\Lambda)$ is the trace of $P/PJ(\Lambda)$. This shows that $\varphi$ is well defined. The injectivity of $\varphi$ follows from $(i)$.
\end{Proof}

The following lemma states that projective modules can be lifted modulo a trace ideal. This is a crucial result in the development of the theory of relatively big projective modules.

\begin{Lemma}
\label{traces.quo}
\cite[Theorem 3.1, Corollary 3.2]{traces}
Let $\Lambda$ be a ring, and let $I$ be the trace of a projective right $\Lambda$-module. Let $P'$ be a projective right $\Lambda/I$-module. Then there is a projective right $\Lambda$-module $X$ such that $I\subseteq \Tr_\Lambda(X)$, $X/XI\cong P'$ and $\Tr_\Lambda(X)/I=\Tr_{\Lambda/I}(P')$. 

Moreover, if $I$ is an element in $\mathcal T(\Lambda)$ and $P'$ is countably generated, then $X$ can be taken to be countably generated.
\end{Lemma}

\subsection{Relatively big projectives}

First we introduce the auxiliary notion of  big projective module with respect to an ideal. It was first introduced by P.~P\v ríhoda in \cite{P2} as a way of generalizing Bass' uniformly big projective modules. We motivate the definition of this class of projective modules by comparing it with a characterization of the class of modules generated by a projective module.

Let $\Lambda$ be a ring, and let $M$ be a  right $\Lambda$-module. We denote by \textit{$\mathrm{Gen}(M)$ the class of modules generated by $M$}, that is  the class of right $\Lambda$-modules which are a homomorphic image of a direct sum of copies of $M$.

\begin{Lemma}
\label{class_generated}
\emph{\cite[Lemma 2.10]{wiegand}}
Let $\Lambda$ be a ring, and let $P$ be a projective right $\Lambda$-module with trace ideal $I$. Then
    \[\mathrm{Gen}(P)=\{M\in\Mod\Lambda: MI=M\}.\]
\end{Lemma}

\begin{Def} \label{def:ibig}
Let $\Lambda$ be a ring, let $P$ be a countably generated projective right $\Lambda$-module, and let $J$ be an ideal of $\Lambda$. The module $P$ is said to be \textit{big with respect to $J$}, or \textit{$J$-big}, if  every countably generated projective right $\Lambda$-module $Q$ satisfying $QJ=Q$ is a homomorphic image of $P$.
\end{Def}

The following remarks show that the significant case of Definition~\ref{def:ibig} is when $J$ is a trace ideal of a projective module. 

\begin{Remarks}
Let $\Lambda$ be a ring, and let $J$ be an ideal of $\Lambda$.
\begin{enumerate}
    \item [(1)] If $Q$ is a countably generated projective right $\Lambda$-module such that $QJ=Q$ then $\Tr _\Lambda (Q)\le J$. Therefore, if $J$ does not contain a non-zero element of $\mathcal T(\Lambda)$ then, by Lemma~\ref{ds.traces},  Definition~\ref{def:ibig} is trivially fulfilled, so any countably generated projective module is $J$-big or, equivalently, $0$-big.  
    \item[(2)] If $P$ is $J$-big, and $K$ is a non-zero element of $\mathcal T(\Lambda)$ such that $K\subseteq J$, then $P$ is $K$-big (cf. Lemma~\ref{ds.traces}). 
    
    Observe that $$J'=\sum _{\substack{K\le J \\ K\in \mathcal T(\Lambda)}} K$$
     is again a trace ideal of a projective module by Lemma~\ref{ds.traces}. Then a countably generated projective module is $J$-big  if and only if it is $J'$-big. 

     We stress that $J'$ does not need to be the trace of a countably generated projective $\Lambda$-module. So it is not necessarily an element in $\mathcal T(\Lambda)$.
\end{enumerate}
\end{Remarks}

In the next lemma, we summarize what it means to be  $J$-big with respect to  $J\in \mathcal T(\Lambda)$.

\begin{Lemma} \label{lem:big_fg}
Let $\Lambda$ be a ring, and let $P,Q$ be countably generated projective right $\Lambda$-modules. Let $J=\Tr_\Lambda(Q)$. Then
\begin{enumerate}
    \item[(i)] $P$ is $J$-big if and only if $P\oplus Q^{(\omega)}\cong P$.
    \item[(ii)] If $P$ is $J$-big and $P/PJ$ is finitely generated, then $J$ is the minimal ideal of $\Lambda$ such that $P/PJ$ is finitely generated.
\end{enumerate}
\end{Lemma}
\begin{Proof}
$(i)$. Since $J=\Tr_\Lambda(Q)=\Tr_\Lambda(Q^{(\omega)})$, by Lemma~\ref{ds.traces}, $Q^{(\omega)}J=Q^{(\omega)}$. If $P$ is $J$-big, then there is an epimorphism $P\to Q^{(\omega)}$. Since $Q^{(\omega)}$ is projective, this epimorphism splits, so $Q^{(\omega)}$ is isomorphic to a direct summand of $P$. Write $P\cong P'\oplus Q^{(\omega)}$. Then, $P\cong P'\oplus Q^{(\omega)}\cong P'\oplus Q^{(\omega)}\oplus Q^{(\omega)}\cong P\oplus Q^{(\omega)}$.

Conversely, if $P\oplus Q^{(\omega)}\cong P$, then $Q^{(\omega)}$ is isomorphic to a direct summand of $P$, so there exists an epimorphism $P\to Q^{(\omega)}$. 
Any countably generated projective module with trace contained in $J$ is an epimorphic image of $Q^{(\omega)}$ by Lemma~\ref{class_generated}, so $P$ is $J$-big.

Statement $(ii)$ is \cite[Lemma 3.3]{wiegand}.
\end{Proof}

The following lemma states that the isomorphism class of a $J$-big projective right $\Lambda$-module $P$ satisfying that $P/PJ$ is finitely generated, is determined, up to isomorphism, by the ideal $J$ and the isomorphism class of $P/PJ$.

\begin{Lemma}
\label{Jbigdetermined}
Let $\Lambda$ be a ring, let $J\in \mathcal T(\Lambda)$, and let $P$ be a $J$-big module such that $P/PJ$ is finitely generated. Then:
\begin{itemize}
    \item[(i)] If $P'$ is a projective module isomorphic to $P$, then $P'$ is $J$-big and $P'/P'J$ is finitely generated.
    \item[(ii)] \cite[Proposition 3.4]{wiegand} Let $P'$ be a countably generated projective module that is $J'$-big for some $J'\in \mathcal T(\Lambda)$ and such that $P'/P'J'$ is finitely generated. Then $P\cong P'$ if and only if $J=J'$ and $P/PJ\cong P'/P'J(=P'/P'J')$.
\end{itemize}
\end{Lemma}

\begin{Proof} We only need to prove $(i)$. Notice that a module isomorphic to a $J$-big module is also $J$-big. Moreover, if $P\cong P'$ then $P'/P'J$ is a homomorphic image of $P/PJ$. So if $P/PJ$ is finitely generated, then so is $P'/P'J$. \end{Proof}

Now we give  the   definition of \emph{relatively big} projective which is crucial for the rest of the paper.

\begin{Def}\label{def:relativelybig}
Let $\Lambda$ be a ring, and let $P$ be a countably generated projective right $\Lambda$-module. The module $P$ is said to be \emph{relatively big} if there exists a countably generated projective right $\Lambda$-module $Q$ with trace ideal $J$ satisfying:
\begin{enumerate}
    \item $P$ is $J$-big, that is, $P\oplus Q^{(\omega)}\cong P$.
    \item $P/PJ$ is a finitely generated projective right $\Lambda/J$-module.
\end{enumerate}
We will also use the terminology $P$ is \emph{relatively $J$-big} when we want to specify the ideal $J\in \mathcal{T} (\Lambda)$.
\end{Def}

\begin{Ex}\label{rem.ex}
Let $\Lambda$ be a ring. Let $X$ be a finitely generated projective right $\Lambda$-module, and let $P$ be a  countably generated projective right $\Lambda$-module. Then $X\oplus P^{(\omega)}$ is relatively big: take $Q=P$ and $J=\Tr_\Lambda(P)$ in the Definition~\ref{def:relativelybig}. 

In particular, taking $X=\{0\}$, we see that $P^{(\omega)}$ is relatively big for any countably generated projective module $P$. Also, taking $P=\{0\}$ we see that any finitely generated projective module is relatively big.
\end{Ex}

Example~\ref{rem.ex}, is a particular instance of the following general construction of relatively big projective modules.

\begin{Lemma}\label{lem:big_cg}
Let $\Lambda$ be a ring. Let $P$ be a countably generated projective right $\Lambda$-module, and let $J$ be the trace ideal of a countably generated projective right $\Lambda$-module $Q$. If $P/PJ$ is finitely generated, then $P\oplus Q^{(\omega)}$ is relatively $J$-big. 
\end{Lemma}
\begin{Proof}
Note that $(P\oplus Q^{(\omega)})/(P\oplus Q^{(\omega)})J\cong P/PJ$, which is finitely generated. Since $P\oplus Q^{(\omega)}\cong P\oplus Q^{(\omega)}\oplus Q^{(\omega)}$, by Lemma~\ref{lem:big_fg}, $P\oplus Q^{(\omega)}$ is $J$-big.
\end{Proof}

\begin{Remark}\label{rem:ind.relbig}
Note that an indecomposable relatively big projective module must be finitely generated.
\end{Remark}

The following result characterizes when every relatively big projective module is isomorphic to a direct sum of finitely generated projective modules.

\begin{Prop} \label{monoidbtrivial} 
Let $\Lambda$ be a ring. Then every relatively big projective right $\Lambda$-module is a direct sum of finitely generated projective modules if and only if the following two properties are satisfied:
\begin{enumerate}
    \item If $J$ is the trace ideal of a  projective right $\Lambda$-module, then there is a right $\Lambda$-module $P$, which is a direct sum of finitely generated projective modules such that $\Tr_\Lambda(P)=J$.
    \item If $J$ is the trace ideal of a countably generated projective right $\Lambda$-module and $P'$ is a finitely generated projective right $\Lambda/J$-module, then there is a finitely generated projective right $\Lambda$-module $X$ such that $P'\cong X/XJ$.
\end{enumerate}
\end{Prop}
\begin{Proof}
Assume that (1) and (2) are satisfied, and let $P$ be a relatively big projective right $\Lambda$-module. Then there exists a countably generated projective right $\Lambda$-module $Q$ with trace ideal $J$ such that $P$ is relatively $J$-big.

By (1), $Q$ can be taken to be a direct sum of finitely generated projective modules. By (2) applied to $P/PJ$, there exists a finitely generated projective right $\Lambda$-module $X$ such that $X/XJ\cong P/PJ$. Let $P'=X\oplus Q^{(\omega)}$. By Lemma~\ref{lem:big_cg}, $P'$ is relatively big, and by Lemma~\ref{Jbigdetermined}, $P\cong P'$, so it is a direct sum of finitely generated projective modules.

Assume now that every relatively big projective right $\Lambda$-module is a direct sum of finitely generated projective modules. Let $P$ be a projective right $\Lambda$-module with trace $J$. By \cite[Theorem~1]{kaplansky}, $P=\bigoplus_{i\in I} P_i$ is a direct sum of countably generated projective modules. By Example~\ref{rem.ex}, $P_i^{(\omega)}$ is relatively big. Therefore, $P_i^{(\omega)}$ is isomorphic to a direct sum of finitely generated projective modules for each $i\in I$. Hence, $P^{(\omega)}$ is isomorphic to a direct sum of finitely generated projective modules and $J=\Tr_\Lambda(P^{(\omega)})$. This gives (1).

To prove (2), let $Q$ be a countably generated projective right $\Lambda$-module with trace $J$, and let $P'$ be a finitely generated projective right $\Lambda/J$-module. By Lemma~\ref{traces.quo}, since $J$ is the trace of a countably generated projective module, there exists a countably generated projective right $\Lambda$-module $P$ such that $J\subseteq \Tr_\Lambda(P)$, $P/PJ\cong P'$, and $\Tr_\Lambda(P)/J=\Tr_{\Lambda/I}(P')$. 

By Lemma~\ref{lem:big_cg}, $P\oplus Q^{(\omega)}$ is relatively big, so it is a direct sum of finitely generated projective modules. Write $P\oplus Q^{(\omega)}=\bigoplus_{i\in I} X_i$ with $X_i$ a finitely generated projective $\Lambda$-module for each $i\in I$. Since $(\bigoplus_{i\in I} X_i)/(\bigoplus_{i\in I} X_i)J\cong P'$ is finitely generated, there is a finite subset $I_0\subseteq I$ such that $(\bigoplus_{i\in I_0} X_i)/(\bigoplus_{i\in I_0} X_i)J\cong P'$, so $X=\bigoplus_{i\in I_0} X_i$ is the finitely generated projective $\Lambda$-module we were looking for.
\end{Proof}

\subsection{Finite generation modulo an ideal and the \texorpdfstring{$(*)$}{(*)}-condition} \label{subsec:fg}

In this subsection, we take a closer look at condition (2) in Definition~\ref{def:relativelybig}.

Let $\Lambda$ be a ring, and let $M$ be a right $\Lambda$-module.
Consider the set
\[\mathcal{I}(M)=\{I\unlhd \Lambda\mid M/MI\mbox{ is finitely generated}\}\]

We recall the following crucial properties of a descending sequence of right ideals associated with a column finite idempotent matrix. The result is taken from \cite{traces}, but all the necessary ideas are already in \cite{P2}.

\begin{Lemma}\label{descending} \cite[Lemmas 5.1 and 5.2]{traces} 
Let $\Lambda$ be a ring. Let $A=(a_{i,j})$ be an $\N \times \N$-matrix, with entries in $\Lambda$, that is a column finite matrix  and such that $A^2=A$. Set $P=A\Lambda^{(\N)}$. For any  $k\ge 0$, let $I_k=\sum_{i>k, j\in \N}a_{ij}\Lambda$. Then,
\begin{enumerate}
    \item[(i)] $\Tr_\Lambda(P)=\Lambda I_0$.
    \item[(ii)] For any $k\ge 0$ and any $0<i\le k$, 
        $$\sum _{j \in\N}a_{ij}\Lambda+\Lambda I_k/\Lambda I_k=\sum _{1\le \ell\le k}a_{i\ell}\Lambda +\Lambda I_k/\Lambda I_k.$$
    In particular, if $0\le n\le k$, $I_n+\Lambda I_k/\Lambda I_k$ is finitely generated.
    \item[(iii)] For any $k\ge 0$, there exists $n_k>k$ such that $I_{n_k}I_k=I_{n_k}$. In particular, the descending sequence of two-sided ideals $(\Lambda I_k)_{k\ge 0}$ has a subsequence $(J_n)_{n\ge 0}$ of two-sided ideals satisfying that:
    \begin{itemize}
        \item[(a)] for any $n\ge 0$, $J_{n+1}J_{n}=J_{n+1}$;
        \item[(b)] for any $k\ge 0$, there exists $n$ and $\ell$ such that $J_n=\Lambda I_{k+\ell}\subseteq \Lambda I_k$.
    \end{itemize}
    \item[(iv)] For any $k\ge 0$, $\Lambda I_k\in \mathcal{I} (P)$. Moreover, for any $I\in \mathcal{I} (P)$, there exists $k_0$ such that $\Lambda I_{k_0}\le I$.
    \item[(v)]  $\mathcal{I}(P)$ has minimal elements  if and only if there exists $k_0\in \No$ such that $\Lambda  I_{k_0}=\Lambda  I_{k_0+\ell}$ for every $\ell \in \N$. In this case, $I_{k_0}=I_{k_0}^2$, and $\Lambda I_{k_0}$ is the minimal element of $\mathcal{I}(P)$.
\end{enumerate}
\end{Lemma}

Statement $(v)$ in Lemma~\ref{descending} motivates the following definition.

\begin{Def}\label{def:badger} \cite{P2}
A ring $\Lambda$ satisfies condition  $(*)$ provided that every descending chain of two-sided ideals
    \[J_1\supseteq J_2\supseteq \cdots \supseteq J_k\supseteq \cdots\]
such that  $J_{k+1}J_k=J_{k+1}$, for each $k\ge 1$, is eventually stationary. 
\end{Def}

The following statement explains the role of condition $(*)$ in our discussion. 

\begin{Cor} \label{cor:starp} 
Let $\Lambda$ be a ring. Let $P$ be a countably generated projective right $\Lambda$-module. If $\Lambda$ satisfies condition $(*)$, then $\mathcal{I} (P)$ has a minimal element, and this minimal element  is an idempotent ideal. 
\end{Cor}

\begin{Proof} 
By Lemma~\ref{descending}~$(v)$, $\mathcal{I} (P)$ has a minimal element if and only if the descending sequence of two-sided ideals $(\Lambda I_k)_{k\ge 0}$ is stationary. By Lemma~\ref{descending}~$(iii)$, this will happen if and only if   the  partial subsequence  $(J_n)_{n\ge 0}$ given by that statement is stationary. Since such a partial subsequence satisfies that   $J_{n+1}J_n=J_{n+1}$, for any $n\ge 0$, it is stationary because $\Lambda$ satisfies condition $(*)$.

The fact that the minimal element in $\mathcal{I} (P)$ is an idempotent ideal follows from Lemma~\ref{descending}~$(v)$.
\end{Proof}

\begin{Th} \cite[Lemma~2.8]{P2} \label{th:star} 
Let $\Lambda$ be a noetherian ring. Let $P$ be a countably generated projective right $\Lambda$-module such that  $\mathcal{I} (P)$ has a minimal element $J$. Then $J$ is the trace ideal of a countably generated projective right $\Lambda$-module, and $P$ is $J$-big.
In particular, if $\Lambda$ satisfies condition $(*)$, then  any countably generated projective right module is relatively big.
\end{Th}

\begin{Remarks}
$(1)$ Let $\Lambda$ be a semilocal ring, and let $P$ be a countably generated projective right $\Lambda$-module. Then, by \cite[Proposition~5.6]{traces}, $\mathcal{I} (P)$ has   a minimal element $J$. Then $J$ is an idempotent ideal, but it does not need to be the trace of a projective module as it can even be contained in the Jacobson radical.

$(2)$ We do not know whether condition $(*)$ is left-right symmetric, even for the case of noetherian rings.

$(3)$ There are noetherian rings with infinite descending sequences of idempotent ideals, so they do not satisfy condition $(*)$. Over such rings there may be countably generated projective modules that are not relatively big.  See, for example, the last section of \cite{traces}.
\end{Remarks}

\section{The monoid of relatively big projective modules}\label{sec:2}

\begin{Def}
\label{def:monoid}
Let $\Lambda$ be a ring. We denote by $B(\Lambda)$ the subset of $V^*(\Lambda)$ of isomorphism classes of relatively big projective right $\Lambda$-modules which are not finitely generated. We denote by $V(\Lambda)\sqcup B(\Lambda)$ the disjoint union of $V(\Lambda)$ and $B(\Lambda)$, which is also a subset of $V^*(\Lambda)$, and consists of the isomorphism classes of all the relatively big projective right $\Lambda$-modules.
\end{Def}

It is important to observe that, by Lemma~\ref{Jbigdetermined}, if $P$ is a countably generated projective right $\Lambda$-module such that $[P]\in B(\Lambda)$, then $P$ is relatively big, so there exists a non-zero ideal $J\in \mathcal{T} (\Lambda)$ such that $P/PJ$ is finitely generated. Moreover, $J$ and $P/PJ$ uniquely determine $[P]$. We will use this fact freely throughout the paper.

\begin{Lemma}\label{monoid}
Let $\Lambda$ be a ring. Then
\begin{enumerate}
    \item[(i)] $B(\Lambda)$ is a subsemigroup of $V^*(\Lambda)\setminus V(\Lambda)$.
    \item[(ii)] $V(\Lambda)\sqcup B(\Lambda)$ is a submonoid of $V^*(\Lambda)$.
\end{enumerate}
\end{Lemma}
\begin{Proof}
$(i)$. Let $[P_1],[P_2]\in B(\Lambda)$. Then, for each $i=1,2$, there exists a countably generated non-zero projective right $\Lambda$-module $Q_i$ with trace ideal $J_i$ such that $P_i\oplus Q_i^{(\omega)}\cong P_i$ and $P_i/P_iJ_i$ is a finitely generated projective right $\Lambda/J_i$-module. Take $P=P_1\oplus P_2$ and $Q=Q_1\oplus Q_2$. Then
    $$P\oplus Q^{(\omega)}\cong P_1\oplus P_2\oplus Q_1^{(\omega)}\oplus Q_2^{(\omega)}\cong P_1\oplus Q_1^{(\omega)}\oplus P_2\oplus Q_2^{(\omega)}\cong P_1\oplus P_2=P.$$
By Lemma~\ref{ds.traces}, $J=\Tr_\Lambda(Q)=\Tr_\Lambda(Q_1)+\Tr_\Lambda(Q_2)=J_1+J_2$, so 
    \[P/PJ\cong \frac{P_1}{P_1(J_1+J_2)} \oplus \frac{P_2}{P_2(J_1+J_2)}\cong \frac{P_1/P_1J_1}{P_1(J_1+J_2)/P_1J_1}\oplus \frac{P_2/P_2J_2}{P_2(J_1+J_2)/P_2J_2}\]  
and $P/PJ$ is a finitely generated projective right $\Lambda/J$-module because it is a direct sum of quotients of finitely generated projective modules. Therefore, as $P$ is not finitely generated, $[P]\in B(\Lambda)$.

Statement $(ii)$ is proved similarly, taking into account that the modules in $V(\Lambda)$ are the relatively $0$-big projective modules.
\end{Proof}

\begin{Lemma}
\label{embedding}
Let $\Lambda$ be a ring, and let $I\in\mathcal T(\Lambda)$. Then, there is a semigroup embedding $\iota _I\colon V(\Lambda/I)\hookrightarrow V(\Lambda)\sqcup B(\Lambda)$. The image of $\iota _I$ is all the isomorphism classes of relatively $I$-big projective right $\Lambda$-modules.
\end{Lemma}
\begin{Proof}
Let $Q$ be a countably generated projective right $\Lambda$-module with trace ideal $I$. Let $P'$ be a finitely generated projective right $\Lambda/I$-module. By Lemma~\ref{traces.quo}, there exists a projective right $\Lambda$-module $P$ such that $I\subseteq \Tr_\Lambda(P)$, $P/PI\cong P'$, and $\Tr_\Lambda(P)/I=\Tr_{\Lambda/I}(P')$.

By Lemma~\ref{lem:big_cg}, $P\oplus Q^{(\omega)}$ is relatively $I$-big, and then $[P\oplus Q^{(\omega)}]\in V(\Lambda)\sqcup B(\Lambda)$. The isomorphism class of $P\oplus Q^{(\omega)}$ is uniquely determined by $I$ and $P'$.  Hence, there is a well-defined embedding   $\iota_I\colon V(\Lambda/I)\longrightarrow V(\Lambda)\sqcup B(\Lambda)$  defined by $\iota_I([P'])=[P\oplus Q^{(\omega)}]$. 

The map $\iota _I$ is a semigroup homomorphism because, for $P'_1$ and $P'_2$ in $V(\Lambda /I)$: 
    \begin{align*}
        \iota_I([P_1'\oplus P_2'])&=[P_1\oplus P_2\oplus Q^{(\omega)}]=[P_1\oplus Q^{(\omega)}\oplus P_2\oplus Q^{(\omega)}] \\
        &=[P_1\oplus Q^{(\omega)}]+[P_2\oplus Q^{(\omega)}]=\iota_I([P_1'])+\iota_I([P_2']).
    \end{align*} 
where for $i=1,2$, $P_i$ denotes a countably generated projective module such that $P_i/P_i\cong P'_i$, which exists because of Lemma~\ref{traces.quo}.

Let $P$ be a relatively $I$-big projective right $\Lambda$-module. Then $P/PI$ is a finitely generated projective right $\Lambda /I$-module, and $\iota _I([P/PI])=[P]$.
\end{Proof}

If $I\neq 0$, the embedding $\iota _I$ of Lemma~\ref{embedding}   is not a monoid embedding because $[0]$ maps to the isomorphism class of a non-finitely generated projective module.

\begin{Remark}\label{operationmonoid}
Consider the set $\mathcal T(\Lambda)$ and the disjoint union $\bigsqcup_{I\in\mathcal T(\Lambda)} V(\Lambda/I)$. We will denote the elements of $\bigsqcup_{I\in\mathcal T(\Lambda)} V(\Lambda/I)$ by $([P],I)$ with $[P]\in V(\Lambda/I)$. Note that if $I\subseteq J$, there is a morphism of monoids $f_{JI}\colon V(\Lambda/I)\to V(\Lambda/J)$ defined by $f_{JI}([X],I)=([X/XJ],J)$ with $[X]\in V(\Lambda/I)$. Therefore, $f_{JI}(V(\Lambda/I))\subseteq V(\Lambda/J)$. Define addition in $\bigsqcup_{I\in\mathcal T(\Lambda)} V(\Lambda/I)$ by
\begin{align*}
    ([P_1],I)+([P_2],J)&=f_{I+J,I}([P_1],I)+f_{I+J,J}([P_2],J) \\
    &=([P_1/P_1J],I+J)+([P_2/P_2I],I+J) \\
    &=([P_1/P_1J\oplus P_2/P_2I],I+J)
\end{align*}
for each $[P_1]\in V(\Lambda/I)$, $[P_2]\in V(\Lambda/J)$. Therefore, the universal property of the disjoint union and Lemma~\ref{embedding} imply that there exists a monoid embedding 
$$\iota\colon\bigsqcup_{I\in\mathcal T(\Lambda)} V(\Lambda/I)\to V(\Lambda)\sqcup B(\Lambda).$$
\end{Remark}

The following result states that $\iota$ is  an isomorphism of monoids.  

\begin{Cor}\label{iso.monoid}
\cite[Theorem 3.8]{wiegand}
Let $\Lambda$ be a ring. Then
\begin{enumerate}
    \item[(i)] $\bigsqcup_{I\in\mathcal T(\Lambda)} V(\Lambda/I)\cong V(\Lambda)\sqcup B(\Lambda)$ as monoids.
    \item[(ii)] $\bigsqcup_{I\in\mathcal T(\Lambda)\setminus \{0\}} V(\Lambda/I)\cong B(\Lambda)$ as semigroups.
\end{enumerate}    
\end{Cor}
\begin{Proof}
The isomorphism in $(i)$ is given by $\iota ([X],I)=\iota_I([X])$, where $\iota_I$ is the embedding from Lemma~\ref{embedding}. By Lemma~\ref{Jbigdetermined}, $\iota$ is bijective.

To prove $(ii)$, note that elements in $V(\Lambda)$ can be thought of as elements in $V(\Lambda/I)$ for $I=\{0\}$, so $(ii)$ follows from $(i)$.
\end{Proof}

\begin{Remark}\label{rem:fgproj}
Let $\Lambda$ be a ring. It is well known that the functor $(-)^*=\mathrm{Hom}_\Lambda (-,\Lambda)$ induces a duality between finitely generated projective right $\Lambda$-modules and finitely generated projective left $\Lambda$-modules. So that,   $\varphi\colon V_r(\Lambda)\to V_\ell(\Lambda)$   defined by $\varphi([P])=[P^*]$, is an isomorphism of monoids. Moreover, $P$ and $P^*$ have the same trace  (cf. \cite[Exercise 22.1]{andersonfuller}).
\end{Remark}

In the next example, we show that the canonical monoid isomorphism $V_r(\Lambda)\cong V_\ell(\Lambda)$ of Remark~\ref{rem:fgproj} does not, in general, extend to an isomorphism between the monoids $V_r(\Lambda)\sqcup B_r(\Lambda)$ and $V_\ell(\Lambda)\sqcup B_\ell(\Lambda)$. We will show that it does extend, provided that $\mathcal T_r(\Lambda)=\mathcal T_\ell(\Lambda)$.

\begin{Ex} \label{ex:asymmetry} 
In \cite[Example~3.6 (iv)]{TAMS} it is given an example of a semilocal ring $\Lambda$ such that $\Lambda/J(\Lambda)\cong D_1\times D_2$, with $D_1$ and $D_2$ division rings, such that any projective left $\Lambda$-module is free but this is not true on the right. Though, on the right, it is still true that $V(\Lambda )\sqcup B(\Lambda )=V^*(\Lambda)$.

In such $\Lambda$, what happens is that one of the maximal ideals contains a non-zero  trace of a countably generated projective right $R$-module. Therefore $\mathcal T_\ell(\Lambda)=\{0, \Lambda \}\subsetneq \mathcal T_r(\Lambda)$. Then $V_r(\Lambda)\sqcup B_r(\Lambda)$ and $V_\ell(\Lambda)\sqcup B_\ell(\Lambda)$ are not isomorphic. This cannot happen if $\mathcal T_r(\Lambda)=\mathcal T_\ell(\Lambda)$.
\end{Ex}

Assume $\Lambda$ is a ring such that $\mathcal T_r(\Lambda)=\mathcal T_\ell(\Lambda)$. We define a map $V_r(\Lambda)\sqcup B_r(\Lambda)\mapsto V_\ell(\Lambda)\sqcup B_\ell(\Lambda)$ that extends the isomorphism between the monoids of finitely generated projective modules of Remark~\ref{rem:fgproj}  by  $([P],I)\mapsto([P^*],I)$ for any $([P],I)\in V_r(\Lambda)\sqcup B_r(\Lambda)$. We must check that $([P_1],I)+([P_2],J)$ is mapped to $([P_1^*],I)+([P_2^*],J)$.    

The following   technical lemma  is crucial for our discussion.

\begin{Lemma}\label{quo.dual}
Let $\Lambda$ be a ring, and let $I,J$ be ideals of $\Lambda$. Let $P$ be a finitely generated projective right $\Lambda/I$-module. Then 
    \[\Hom_{\Lambda/(I+J)}(P/PJ,\Lambda/(I+J))\cong\Hom_{\Lambda/I}(P,\Lambda/I)/((I+J)/I \cdot\Hom_{\Lambda/I}(P,\Lambda/I)).\]
\end{Lemma}
\begin{Proof}
Note that we can identify $P/PJ$ with $P\otimes_{\Lambda/I} \Lambda/(I+J)$. Then
    \begin{align*}
        \Hom_{\Lambda/(I+J)}(P/PJ,\Lambda/(I+J))
        &\cong\Hom_{\Lambda/(I+J)}(P\otimes_{\Lambda/I} \Lambda/(I+J),\Lambda/(I+J)) \\
        &\cong\Hom_{\Lambda/I}(P,\Lambda/(I+J)).
    \end{align*}
where in the second isomorphism, we have used the Hom-tensor adjunction. Now, we want to prove that $\Hom_{\Lambda/I}(P,\Lambda/(I+J))\cong\Hom_{\Lambda/I}(P,\Lambda/I)/J\Hom_{\Lambda/I}(P,\Lambda/I)$ as left $\Lambda/I$-modules.

Consider the map $\varphi\colon \Hom_{\Lambda/I}(P,\Lambda/I)\to\Hom_{\Lambda/I}(P,\Lambda/(I+J))$ defined by $\varphi(f)=\pi\circ f$, where $\pi\colon \Lambda/I\to\Lambda/(I+J)$ is the projection modulo $(I+J)/I$. Clearly, it is a homomorphism of left $\Lambda/I$-modules. Let $g\in\Hom_{\Lambda/I}(P,\Lambda/(I+J))$. Since $P$ is projective, there is a lift $G\in\Hom_{\Lambda/I}(P,\Lambda/I)$ making the following diagram commute
    \[\begin{tikzcd}[row sep=huge]
        & P\dar{g}\dlar[dashed,swap]{G}\\
        \Lambda/I\rar{\pi} &\Lambda/(I+J)\rar & 0
    \end{tikzcd}\]
Therefore, $\varphi$ is surjective. The kernel of this homomorphism is
    \begin{align*}
        \ker\varphi&=\{f\in\Hom_{\Lambda/I}(P,\Lambda/I)\mid \pi\circ f=0\} \\
        &=\{f\in\Hom_{\Lambda/I}(P,\Lambda/I)\mid f(x)\in (I+J)/I\text{ for every }x\in P\} \\
        &=\Hom_{\Lambda/I}(P,(I+J)/I).
    \end{align*}
It remains to prove that $\Hom_{\Lambda/I}(P,(I+J)/I)=(I+J)/I\cdot\Hom_{\Lambda/I}(P,\Lambda/I)$. 

Let $f\in\Hom_{\Lambda/I}(P,(I+J)/I)$. Consider $(I+J)/I$ as a right $\Lambda/I$-module. Then there exists an epimorphism $\pi\colon(\Lambda/I)^{(A)}\to (I+J)/I$ for some set $A$. Since $P$ is projective, there exists a lift $F\in\Hom_{\Lambda/I}(P,(\Lambda/I)^{(A)})$ such that $\pi\circ F=f$, that is, the following diagram commutes
    \[\begin{tikzcd}[row sep=huge]
        & P\dar{f}\dlar[dashed,swap]{F}\\
        (\Lambda/I)^{(A)}\rar{\pi} &(I+J)/I\rar & 0.
    \end{tikzcd}\]
Since $P$ is finitely generated, the image of $F$ is contained in a finite number of copies of $\Lambda/I$, namely $(\Lambda/I)^{(A_0)}$, with $A_0\subseteq A$ finite. Let $\{x_i\}_{i\in A}$ be the canonical basis of $(\Lambda/I)^{(A)}$. Then
    \[f(x)=\pi\circ F(x)=\pi\Big(\sum_{i\in A_0} x_ia_i\Big)=\sum_{i\in A_0} \pi(x_i)a_i=\sum_{i\in A_0} \pi(x_i)\cdot\varphi_i(F(x)).\]
where $\varphi_i\colon (\Lambda/I)^{(A_0)}\to \Lambda/I$ denotes the projection into the $i$\textsuperscript{th} component, $i\in A_0$. Hence, $f=\sum_{i\in A_0} \pi(x_i)\cdot(\varphi_i\circ F)\in(I+J)/I\cdot\Hom_{\Lambda/I}(P,\Lambda/I)$.

Conversely, let $F=\sum_{i=1}^n (y_i+I)\cdot F_i\in(I+J)/I\cdot\Hom_{\Lambda/I}(P,\Lambda/I)$ for some $y_1,\dotsc,y_n\in J$ and $F_1,\dotsc,F_n\in\Hom_{\Lambda/I}(P,\Lambda/I)$. Then, regarding the product of elements as a composition of homomorphisms, $F\in\Hom_{\Lambda/I}(P,(I+J)/I)$. Therefore, by the First Isomorphism Theorem, 
    $$\Hom_{\Lambda/I}(P,\Lambda/I)/((I+J)/I\cdot\Hom_{\Lambda/I}(P,\Lambda/I))\cong\Hom_{\Lambda/I}(P,\Lambda/(I+J)),$$
so the statement is proved.
\end{Proof}

\begin{Cor} \label{isomonoidtraces}
Let $\Lambda$ be a ring such that $\mathcal T_r(\Lambda)=\mathcal T_\ell(\Lambda):=\mathcal T(\Lambda)$. Then there is an isomorphism of monoids between $V_r(\Lambda)\sqcup B_r(\Lambda)$ and $V_\ell(\Lambda)\sqcup B_\ell(\Lambda)$.
\end{Cor}
\begin{Proof}
Let $\varphi\colon \bigsqcup_{I\in\mathcal T(\Lambda)} V_r(\Lambda/I)\to \bigsqcup_{I\in\mathcal T(\Lambda)} V_\ell(\Lambda/I)$ be the map defined by $\varphi([P],I)=([P^*],I)$. By Remark~\ref{rem:fgproj}, it is well-defined and bijective. By Lemma~\ref{quo.dual}, it is a monoid homomorphism since
\begin{align*}
    \varphi(([P_1],I)+([P_2],J))&=\varphi([P_1/P_1J\oplus P_2/P_2I],I+J) \\
    &=([\Hom_{\Lambda/(I+J)}(P_1/P_1J\oplus P_2/P_2I,\Lambda/(I+J))],I+J) \\
    &=([\Hom_{\Lambda/(I+J)}(P_1/P_1J,\Lambda/(I+J)) \\
    &\qquad\qquad\oplus \Hom_{\Lambda/(I+J)}(P_2/P_2I,\Lambda/(I+J))],I+J) \\
    &=([\Hom_{\Lambda/I}(P_1,\Lambda/I)/J\Hom_{\Lambda/I}(P_1,\Lambda/I)) \\
    &\qquad\qquad\oplus \Hom_{\Lambda/J}(P_2,\Lambda/J)/I\Hom_{\Lambda/J}(P_2,\Lambda/J))],I+J) \\
    &=([\Hom_{\Lambda/I}(P_1,\Lambda/I)],I)+([\Hom_{\Lambda/J}(P_2,\Lambda/J)],J) \\
    &=\varphi([P_1,I])+\varphi([P_2],J)
\end{align*}
Finally, as $\varphi$ is a monoid isomorphism, by Corollary~\ref{iso.monoid} there is also a monoid isomorphism $V_r(\Lambda)\sqcup B_r(\Lambda)\to V_\ell(\Lambda)\sqcup B_\ell(\Lambda)$ .
\end{Proof}

Then, Theorem~\ref{th:star}, together with the fact that over a noetherian ring $\Lambda$, $\mathcal T_r(\Lambda)=\mathcal T_\ell(\Lambda)$ is just the set of idempotent ideals of $\Lambda$ (cf. \cite{whitehead}), immediately yields.

\begin{Th} \label{th:starmonoid}
\cite[Theorem 2.12]{P2}
Let $\Lambda$ be a  noetherian ring satisfying condition~$(*)$. Then $V_r^*(\Lambda)=V_r(\Lambda)\sqcup B_r(\Lambda)\cong V_\ell(\Lambda)\sqcup B_\ell(\Lambda)=V_\ell^*(\Lambda)$.
\end{Th}

In general, it is an open question whether $\mathcal T_r(\Lambda)=\mathcal T_\ell(\Lambda)$ for $\Lambda$ a right noetherian ring (cf. Corollary~\ref{PImonoid} and the comments after it for further information regarding this question).

\section{Lifting projectives modules modulo a right \texorpdfstring{$T$}{T}-nilpotent ideal}\label{sec:3} 

\begin{Lemma} \label{liftingproj}
Let $\Lambda $ be a ring, and let $N$ be a right $T$-nilpotent ideal of $\Lambda $. 
\begin{enumerate}
    \item[(i)] Let $X$ be a  projective right $\Lambda /N$-module. Then there exists a projective right $\Lambda $-module $P$ such that $P/PN\cong X$. All projective modules with this property are isomorphic. In addition,  $X$ is finitely generated (countably generated) if and  only if $P$ is finitely generated (countably generated).
    \item[(ii)] Let $P$ and $Q$ be projective right $\Lambda$-modules. Let  $f\colon P\to Q$ be a morphism of $\Lambda$-modules. Then $f$ is an isomorphism if and only if   the induced map $P/PN\to Q/QN$ is an isomorphism of $\Lambda /N$-modules.
    \item[(iii)] Let $I$ be an idempotent ideal of $\Lambda $ such that $(I+N)/N$ is the trace of a countably generated projective right $\Lambda /N$-module $X$. Then $I$ is the trace of any countably generated projective right $\Lambda $-module $P$ such that $P/PN\cong X$.
\end{enumerate}
\end{Lemma}
\begin{Proof} 
$(i)$ Let $X$ be a projective right $\Lambda /N$-module. As $X$ is a direct sum of countably generated projective right $\Lambda /N$-modules, we may assume that $X$ is countably generated. In this case,  $X$ is the direct limit of a direct system
    \begin{equation}\label{eq:1}
        \begin{tikzcd}
            \overline{F}_1 \rar & \overline{F}_2 \rar & \cdots \rar & \overline{F}_n\rar & \cdots
        \end{tikzcd}
    \end{equation}
of finitely generated free $\Lambda /N$-modules, and the canonical presentation of the direct limit
    \begin{equation}\label{eq:2}
        \begin{tikzcd}
            0\rar & \bigoplus _{n\ge 1} \overline{F}_n \rar & \bigoplus _{n\ge 1} \overline{F}_n \rar & X\rar & 0
        \end{tikzcd}
    \end{equation}
splits. We can lift \eqref{eq:1} to a countable direct system of finitely generated free $\Lambda$-modules 
    \begin{equation}\label{eq:3}
        \begin{tikzcd}
            F_1 \rar & F_2 \rar & \cdots \rar & F_n\rar & \cdots
        \end{tikzcd}
    \end{equation}
such that when we apply the functor $-\otimes _\Lambda \Lambda /N$ to \eqref{eq:3}, we obtain a direct system isomorphic to \eqref{eq:1}. 

Let $P$ be the direct limit of \eqref{eq:3}, and consider its  canonical presentation
    \begin{equation*}\label{eq:4}
        \begin{tikzcd}
            0\rar & \bigoplus _{n\ge 1} {F}_n \rar{\Phi} & \bigoplus _{n\ge 1} {F}_n \rar & P\rar & 0
        \end{tikzcd}
    \end{equation*}
As \eqref{eq:2} splits, there exists $\Psi \colon \bigoplus _{n\ge 1} {F}_n \to \bigoplus _{n\ge 1} {F}_n$ such that $\Psi \circ \Phi =\mathrm{Id}+B$ where $B$ can be thought of as a column finite matrix with entries in $N$ (that is, $B\colon \bigoplus _{n\ge 1} {F}_n \to \bigoplus _{n\ge 1} {F}_n$ has its image contained in 
$\bigoplus _{n\ge 1} {F}_nN$). Since $N$ is right $T$-nilpotent, $B$ is in the Jacobson radical of the ring of the $\N \times \N$ column-finite matrices with coefficients in $\Lambda$ \cite{SW}, and therefore, $\mathrm{Id}+B$ is invertible. So that $(\mathrm{Id}+B)^{-1}\circ \Psi$   is  a left inverse of $\Phi$, and we can deduce that $P$ is projective.

This shows that any lifting of the direct system \eqref{eq:1} gives a projective $\Lambda$-module $P$ that, in addition, satisfies $P/PN\cong X$.
As $PN$ is a small submodule of $P$ by \cite[Lemma~28.3]{andersonfuller}, $P/PN\cong X$ is finitely generated (countably generated) if and only if $P$ is also finitely generated (countably generated).

The uniqueness, up to isomorphism, of $P$ follows from $(ii)$.

$(ii)$ If $f\colon P\to Q$ is an isomorphism of $\Lambda$-modules, then it has an inverse homomorphism $g\colon Q\to P$. Let $\overline f\colon P/PN \to Q/QN$ and $\overline g\colon Q/QN \to P/PN$ be the induced maps. Then $\overline g$ is the inverse of $\overline f$, so both are also invertible morphisms of $\Lambda /N$-modules.

Conversely, if $\alpha \colon P/PN\to Q/QN$ is an isomorphism of $\Lambda /N$-modules. As $PN$ is a small submodule of $P$ by \cite[Lemma~28.3]{andersonfuller},   $P$ and $Q$ are projective covers of $P/PN$. Therefore, any lifting of $\alpha$ to a morphism $P\to Q$ is an isomorphism. 

$(iii)$.  By $(i)$,  the module $X$ can be lifted to a countably generated projective right $\Lambda $-module $P$ such that $P/PN\cong X$. Let $J$ be the trace of $P$. Then, by Lemma~\ref{tracaquocient}, $J+N=I+N$. By \cite[Corollary~2.9]{traces},  $J \subseteq I$ because $J$ is a trace ideal of a projective module, and $N$, being right $T$-nilpotent, is contained in the Jacobson radical of $\Lambda $. 
Since $J+N = I+N$ and $I^2 = I$, $IJ + IN = I$. By \cite[Lemma~28.3]{andersonfuller}, $IN \ll I$, which implies $IJ = I$ and, in particular, $I \subseteq J$.
Therefore, $I = J$ is the trace of $P$.
\end{Proof}

\begin{Lemma} \label{Ibig} Let $\Lambda $ be a ring, and let $N$ be a right $T$-nilpotent ideal of $\Lambda $. Let $I$ be an ideal of $\Lambda $, and let $P$ be a countably generated right $\Lambda $-module. Then, $P$ is $I$-big if and only if $P/PN$ is $(I+N)/N$-big.
\end{Lemma}
\begin{Proof} 
Assume that $P$ is $I$-big. Let $X$ be a countably generated projective right $\Lambda /N$-module such that its trace is contained in $(I+N)/N$. By Lemma~\ref{liftingproj}, there exists a countably generated projective right $\Lambda $-module $Q$ such that $Q/QN\cong X$. Therefore, its trace $J$ is contained in $I+N$. Since $J$ is the trace of a countably generated projective, and  $N$ is right $T$-nilpotent, we deduce from \cite[Corollary~2.9]{traces} that $J\subseteq I$. By hypothesis, there exists an onto module homomorphism $f\colon P\to Q$. Therefore, there is an onto module homomorphism
    \[\begin{tikzcd}
        g\colon P\rar{f} & Q \rar{\pi} & Q/QN\cong X
    \end{tikzcd}\]
where $\pi$ denotes the canonical projection.  As $PN$ is contained in the kernel of $g$, the morphism $g$ factors through an onto homomorphism $P/PN\to X$. This shows that $P/PN$ is $(I+N)/N$-big.

To prove the converse, let $Q$ be a countably generated projective right $\Lambda $-module such that its trace is contained in $I$. By Lemma~\ref{tracaquocient}, the trace of $Q/QN$ is contained in $I+N/N$, and since $P/PN$ is $(I+N)/N$-big, there exists an onto module homomorphism $g\colon P/PN\to Q/QN$. Since $P$ is projective, there exists $f\colon P\to Q$ such that the diagram
\[\begin{tikzcd}[row sep=huge]
        &P\arrow[dl,"f"']\arrow[d, " g\circ \pi _1 "]\\
        Q \arrow[r, "\pi _2"]& Q/QN
    \end{tikzcd}\]
is commutative, and where $\pi _1\colon P\to P/PN$ and $\pi _2\colon Q\to Q/QN$ denote canonical projections. Therefore, $Q=f(P)+QN$. By \cite[Lemma~28.3]{andersonfuller},  $QN$ is a small submodule of $Q$, and we deduce that $f(P)=Q$. This shows that $P$ is $I$-big. 
\end{Proof}

\begin{Cor} \label{isomonoidnil} 
Let $\Lambda $ be a ring, and let $N$ be a right $T$-nilpotent ideal of $\Lambda $. Then the functor $-\otimes _\Lambda \Lambda /N$ induces a monoid isomorphism  $\varphi \colon V^*(\Lambda) \to V^*(\Lambda /N)$   which restricts to monoid isomorphisms $V(\Lambda)\to V(\Lambda/N)$ and $ V(\Lambda  )\sqcup B(\Lambda ) \to V(\Lambda /N )\sqcup B(\Lambda /N )$.
\end{Cor}
\begin{Proof} 
Consider the monoid morphism $\varphi \colon V^*(\Lambda) \to V^*(\Lambda /N)$ defined by $\varphi([P])=[P/PN]$, where $P$ denotes a countably generated, projective, right $\Lambda$-module. It follows from Lemma~\ref{liftingproj} that $\varphi$ is an isomorphism. More precisely, $\varphi$ is injective by the claim $(2)$ of  Lemma~\ref{liftingproj}, and it is onto  by claim $(1)$.   Claim $(1)$ of Lemma~\ref{liftingproj} also implies that $\varphi$ restricts to an isomorphism  $V(\Lambda)\to V(\Lambda/N)$.

It remains to prove that $\varphi$ restricts  to an isomorphism $ B(\Lambda ) \to  B(\Lambda /N )$.

Let $P$ be a countably generated projective $\Lambda$-module such that $[P]\in B(\Lambda)$. That is, there exists $I$, the trace of a non-zero countably generated projective module $Q$, such that $P$ is $I$-big and $P/PI$ is finitely generated. 

By Lemma~\ref{tracaquocient}, $(I+N)/N$ is the trace of $Q/QN$. Notice also that $(I+N)/N\neq \{0\}$.
By Lemma~\ref{Ibig}, $P/PN$ is $(I+N)/N$-big. Clearly, $$(P/PN)/(P(I+N)/PN)\cong P/P(I+N)$$ is finitely generated. Hence, $\varphi([P])\in B(\Lambda/N)$.

Let  $X$ be a countably generated  projective right $\Lambda/N$-module such that $[X]\in B(\Lambda /N)$. This means that there exists a non-zero countably generated projective right $\Lambda /N$-module with trace ideal $K$, such that  $X$ is $K$-big and  $X/XK$ is finitely generated.  

By Lemma~\ref{liftingproj}, there exist countably generated projective $\Lambda$-modules $P$ and $Q$ such that $P/PN\cong X$ and $Q/QN\cong Y$. Let $I$ be the trace of $Q$. By Lemma~\ref{tracaquocient}, $K=(I+N)/N$.

By Lemma~\ref{Ibig}, $P$ is $I$-big. Moreover, $P/PI$ is finitely generated because $P/P(I+N)$ is also finitely generated. This shows that $[P]=\varphi^{-1}([X])\in B(\Lambda)$.

This concludes the proof that $\varphi$ induces an isomorphism $  B(\Lambda ) \to B(\Lambda /N )$. 
\end{Proof}

\section{The case of right noetherian PI-rings}\label{sec:4}

In this section, we shall prove that countably generated, projective, right or left modules over a right noetherian PI-ring are always relatively big. The key result will be to show  that these rings satisfy the $(*)$-condition introduced in \S \ref{sec:2}. As a right noetherian PI-ring $R$ has nilpotent nilradical $N(R)$, and $R/N(R)$ is two-sided noetherian, the result will follow by combining Theorem~\ref{th:star} with the results in \S \ref{sec:3}.

First, we consider the case of commutative noetherian rings. P.~P\v r\'\i hoda and G.~Puninski proved in \cite{PP} that commutative noetherian rings satisfy condition $(*)$. Here, we provide another proof of this result, which we find interesting as it is quite different from the one given in \cite{PP}.

\subsection{Commutative noetherian rings}

If $I$ is a countably generated left ideal of a ring $\Lambda$, $I$ is the union of an ascending sequence of the type $\eqref{eq:sequences2}$ in Lemma~\ref{sequences} if and only if $I$ is a  pure left ideal \cite[Lemme~2]{lazard}. 

Over a commutative ring, an ideal is the trace of a projective module if and only if it is pure (see Remark~\ref{chartraces}). Since the trace ideal of a countably generated projective module is countably generated,  an ideal is the trace of a countably generated projective module if and only if it is the union of  an ascending sequence of the type $\eqref{eq:sequences2}$ in Lemma~\ref{sequences}. 

\begin{Lemma} \label{sequences} Let $\Lambda$ be a ring. Let $(a_k)_{k\ge 1}$ be a sequence of elements of $\Lambda$ satisfying $a_ka_{k+1} =a_{k+1}$ for any $k\ge 1$. Set $b_k=1-a_k$ for any $k\ge 1$. Then:
\begin{enumerate}
    \item [(i)] There is a descending sequence of principal right ideals
        \begin{equation}\label{eq:sequences1}
            a_1\Lambda\supseteq a_2\Lambda\supseteq \cdots \supseteq a_k\Lambda\supseteq \cdots
        \end{equation}
    \item [(ii)] $b_kb_{k+1}=b_k$ for any $k\ge 1$, so there is an ascending sequence of principal left ideals
        \begin{equation}\label{eq:sequences2}
            \Lambda b_1\subseteq \Lambda b_2\subseteq \cdots \subseteq \Lambda b_k\subseteq \cdots
        \end{equation}
    \item[(iii)] For any $k\ge 1$, 
        $$\Lambda b_k \subseteq \mathrm{l.ann} (a_{k+1})\subseteq \Lambda b_{k+1},$$ 
    and 
        $$a_{k}\Lambda \supseteq \mathrm{r.ann} (b_{k})\supseteq a_{k+1}\Lambda.$$
    \item[(iv)] $I=\bigcup _{k\ge 1}\Lambda b_k$ is a left pure ideal of $R$, and it  is also projective.
    \item[(v)] The sequence \eqref{eq:sequences1} is stationary if and only if the sequence \eqref{eq:sequences2} is also stationary. In this case, there exists $k_0$ such  that $b_{k_0}$ and $a_{k_0}$ are idempotent, $I=\Lambda b_{k_0}$ and $\bigcap _{k\ge 1} a_k \Lambda =a_{k_0}\Lambda$.
\end{enumerate}
\end{Lemma}

\begin{Proof} Statements $(i)$, $(ii)$, and $(iii)$ follow easily from direct computation.

The purity claimed in statement $(iv)$ was proved by Lazard in \cite[Lemme~2]{lazard}. Then, as $\Lambda /I$ is a countably presented flat left $\Lambda$-module, $I$ is projective. For a self-contained proof, the reader can check \cite[Lemma~2]{FHS1}.

Statement $(iii)$ readily implies that \eqref{eq:sequences1} is stationary if and only if \eqref{eq:sequences2} is also stationary.  

If $a_k\Lambda =a_{k+1}\Lambda$, then there exists $\lambda \in \Lambda$ such that $a_{k+1}\lambda=a_k$. Hence, $a_{k+1}\lambda a_{k+1}=a_{k+1}$, which implies that $a_k=a_{k+1}\lambda$ is idempotent.

If $\Lambda b_k=\Lambda b_{k+1}$, then there exists $\mu \in \Lambda$ such that $\mu b_k =b_{k+1}$. As before, this implies that $b_{k+1}$ is idempotent.

Assume the sequences \eqref{eq:sequences1} and \eqref{eq:sequences2} are stationary. Then the argument above shows that, taking $k_0$ large enough, the second part of statement $(v)$ holds.
\end{Proof}

\begin{Lemma}\label{determinantaltrick} 
Let $R$ be a commutative ring. Let $I\supseteq J$ be ideals of $R$ such that $JI=J$. If $J$ is finitely generated, then there exists $a\in I$ such that $xa=x$ for any $x\in J$. In particular $I\supseteq aR \supseteq J$.
\end{Lemma}
\begin{Proof} Let $J=x_1R +\dots +x_nR$. The equality $JI=J$ implies that there exists $A\in M_n(I)$ such that $(x_1,\dots, x_n)A=(x_1,\dots, x_n)$, or equivalently
\[(x_1,\dots, x_n)(\mathrm{Id}-A)=(0,\dots ,0).\]
Multiplying the above identity on the right by the adjoint matrix of $\mathrm{Id}-A$, we deduce that $\mathrm{Id}-\det(A)$ annihilates $J$. To conclude the proof, note that there exists $a \in I$ such that $\det(\mathrm{Id}-A) = 1- a$.
\end{Proof}

\begin{Prop} \emph{(\cite[Lemma~2.1]{PP})} \label{prop:starcomm}
Every commutative noetherian ring satisfies condition $(*)$. This is to say that if $R$ is a commutative noetherian ring and
    \begin{equation}\label{eq:starcomm}
        J_1\supseteq J_2\supseteq \cdots \supseteq J_k\supseteq \cdots
    \end{equation}
is a descending chain of ideals of $R$ satisfying $J_{k+1}J_k=J_{k+1}$ for each $k\ge 1$, then \eqref{eq:starcomm} is stationary. In addition,  $\bigcap _{k\ge 1}J_k=eR$ where $e$ is an idempotent  of $R$.
\end{Prop}
\begin{Proof} 
Since $R$ is noetherian, the ideals $J_k$ are finitely generated for every $k\ge 1$. By Lemma~\ref{determinantaltrick}, there exists a sequence of elements $(a_k)_{k\ge 1}$ in $R$ such that
    \[J_1\supseteq a_1R\supseteq J_2\supseteq a_2R\supseteq \cdots \supseteq J_k\supseteq a_kR\supseteq \cdots\]
and such that $a_{k+1}a_k=a_{k+1}$.

By Lemma~\ref{sequences} and because $R$ is noetherian, the ascending chain of ideals
\[(1-a_1)R\subseteq (1-a_2)R\subseteq \cdots \subseteq (1-a_k)R\subseteq \cdots \]
is stationary and its union is generated by an  idempotent of $R$. By Lemma~\ref{sequences}, the descending sequence 
\[a_1R\supseteq a_2R\supseteq \cdots \supseteq a_kR \supseteq \cdots \]
is also stationary, and its intersection is generated by an idempotent of $R$. The sequence \eqref{eq:starcomm} is also stationary, and the intersection of these ideals is generated by an idempotent of $R$.
\end{Proof}

In the following result, we recover a result of Bass that shows that, over commutative noetherian rings, the theory of countably generated projective modules reduces to that of finitely generated ones.

\begin{Cor} \cite{bass} \label{projcommnoe}
Let $R$ be a commutative noetherian ring. Let $P$ be a countably generated projective module that is not finitely generated.  Then there exists an idempotent $e\in R$ with $P\cong (1-e)P\oplus \left ( eR\right)^{(\omega)}$ such that $(1-e)P$ is finitely generated. 
\end{Cor}
\begin{Proof} 
Let $P$ be a countably generated projective module over $R$. By Proposition~\ref{prop:starcomm}, Lemma~\ref{descending}, and Theorem~\ref{th:star}, there exists an idempotent $e\in R$ such that $P/eP$ is finitely generated,  $I=eR$ is an ideal minimal with respect to that property, and $P$ is $I$-big. Therefore $P\cong P\oplus (eR)^{(\omega)}$. 

Since $R$ is commutative
    $$P=eP\oplus (1-e)P\cong eP\oplus (1-e)P\oplus (eR)^{(\omega)}\cong (1-e)P\oplus (eR)^{(\omega)},$$
where the last isomorphism follows from the Eilenberg swindle. To conclude the proof of the statement, observe that $(1-e)P\cong P/eP$. \end{Proof}

Corollary~\ref{projcommnoe} implies that a countably generated projective module over a commutative noetherian domain that is not finitely generated is free. 

Not much is known about the structure of infinitely generated projective modules  over general commutative domains. To our knowledge, the basic questions on the topic posed in 
\cite[p.~ 246]{fuchssalce} are still widely open, and it seems that to answer them, one needs to go beyond the theory of relatively big projectives. 

In general, over a commutative domain $R$, the only trace ideals are $\{0\}$ and $R$, so $B(R)$ contains only the isomorphism class of $R^{(\omega)}$. Therefore, $V(R)\sqcup B(R)=V(R)\sqcup \{[R^{(\omega)}]\}$. We do not know of any reasonable characterization of domains $R$ such that  $V^*(R)=V(R)\sqcup B(R)$, nor are there any examples where this fails.  

\subsection{One-sided noetherian PI-rings} 

In this subsection, we show that   countably generated projective right or left modules over  any right noetherian, PI-ring $\Lambda$ are relatively big. 

The big difference between a noetherian PI-ring $\Lambda$ and   the commutative case described in Corollary~\ref{projcommnoe} is that the idempotent ideals are not just the ideals generated  by  idempotent elements of $\Lambda$.  
By a result of Small and Robson \cite{RS}, the set of idempotent ideals of $\Lambda$ (recall that it coincides with  $\mathcal{T} (\Lambda)$) is finite.  But the results from \cite{P3} and \cite{wiegand} provide plenty of examples showing that, even in the case where $\Lambda$ is a finitely generated algebra over a commutative noetherian ring, $\mathcal{T} (\Lambda)$ frequently contains elements that are not traces of finitely generated projective modules. 

\begin{Lemma}\label{idempotent} 
Let $\Lambda$ be a ring. Let $I\supseteq J$ be ideals of $\Lambda$ such that either $JI=J$ or $IJ=J$. Let $K_1,\dots ,K_n$ be ideals of $\Lambda$ (not necessarily different) such that $K_1\cdots K_n\subseteq J$ and $J+K_i=I+K_i$, for $i=1,\dots ,n$. Then $J=J^2$.
\end{Lemma}
\begin{Proof} 
Using the hypothesis, we deduce that, 
\[I^n\subseteq \prod _{i=1}^n (I+K_i)=\prod _{i=1}^n (J+K_i)\subseteq J.\]
If $JI=J$, then $JI^n=J$. Since $I^n\subseteq J$, we deduce that $J\subseteq J^2$. This shows that $J=J^2$.

A symmetric argument gives the conclusion if $IJ=J$.
\end{Proof}

\begin{Prop} \label{reduction} 
Let $\Lambda $ be a right noetherian ring. Let 
    \begin{equation}\label{eq:reduction1}
        I_1\supseteq I_2\supseteq \cdots \supseteq I_k \supseteq \cdots
    \end{equation}
be a sequence of ideals of $\Lambda $ satisfying that $I_{k+1}I_k=I_{k+1}$ for any $k\ge 1$. Let 
    \begin{equation}\label{eq:reduction2}
        J_1\supseteq J_2\supseteq \cdots \supseteq J_k \supseteq \cdots
    \end{equation}
be a sequence of ideals of $\Lambda $ satisfying that $J_kJ_{k+1}=J_{k+1}$ for any $k\ge 1$.
Then,
\begin{enumerate}
    \item[(i)] Let $N(\Lambda )$ denote the prime radical of $\Lambda $ (that is, the intersection of all prime ideals of $\Lambda $). If \eqref{eq:reduction1} or \eqref{eq:reduction2} are stationary in $\Lambda /N(\Lambda )$, then they are also stationary over $\Lambda $.

    \item[(ii)] Assume that $\Lambda $ is semiprime, and that \eqref{eq:reduction1} is stationary over $\Lambda /P$ for any minimal prime ideal $P$ of $\Lambda $. Then there exists $k_0\ge 1$ such that $I_{k_0+\ell}^2=I_{k_0+\ell}$ for any $\ell \ge 0$. 

    \item[(iii)] Assume that $\Lambda $ is semiprime, and that \eqref{eq:reduction2} is stationary over $\Lambda /P$ for any minimal prime ideal $P$ of $\Lambda $. Then there exists $k_0\ge 1$ such that $J_{k_0+\ell}^2=J_{k_0+\ell}$ for any $\ell \ge 0$.
\end{enumerate}
\end{Prop}

\begin{Proof} $(i)$. By \cite[Theorem~3.11]{GW}, $N(\Lambda )$ is nilpotent. Assume that \eqref{eq:reduction1} is stationary over $\Lambda /N(\Lambda )$. Then there exists $k_1 \ge 0$ such that $I_{k_1}+N(\Lambda )=I_{k_1+\ell}+N(\Lambda )$ for any $\ell \ge 0$. By Lemma~\ref{idempotent}, this implies that $I_{k_1+1+\ell}^2=I_{k_1+1+\ell}$ for any $\ell \ge 0$. Since idempotent ideals that are equal modulo a nilpotent ideal must be equal, the claim follows.

A similar argument works if \eqref{eq:reduction2} is stationary.

$(ii)$. Since $\Lambda $ is a semiprime right Goldie ring, it has a finite number of minimal prime ideals $P_1,\dots ,P_n$, and its product is zero \cite[Proposition~7.1]{GW}. 

By hypothesis, there exists $k_1$ such that $I_{k_1+\ell}+P_i=I_{k_1}+P_i$ for any $\ell \ge 0$ and any $i=1,\dots ,n$.

As in $(i)$, we use Lemma~\ref{idempotent} to deduce that $k_0=k_1+1$ satisfies the conclusion of the statement.

Statement $(iii)$ is proved in a similar way.
\end{Proof}

Robson and Small proved in \cite{RS} that PI-rings satisfying the ascending chain condition on two-sided ideals have a finite number of idempotent ideals.  

We shall prove in Theorem~\ref{starpi} that, in addition, right noetherian PI-rings satisfy  condition $(*)$ and its left analog, which we call $(**)$. The proof is patterned after that of Robson and Small.

First, we provide some basic definitions needed to follow the proof. 

\begin{Def}
A ring $\Lambda$ is said to be a \emph{polynomial identity ring}, a  \emph{PI-ring} for short, if there exists $n$ such that $\Lambda$ satisfies a monic polynomial of $\mathbf{Z}\langle x_1,\dots ,x_n\rangle$.  That is, there exists a polynomial in non-commuting variables $x_1,\dots ,x_n$ with coefficients in $\mathbf{Z}$, say $f(x_1,\dots , x_n)$, such that:
\begin{itemize}
    \item[(1)] $f(r_1,\dots ,r_n)=0$ for any $r_1,\dots ,r_n\in \Lambda$,
    \item[(2)] the coefficient of one of the monomials of the highest degree of $f$ is one.
\end{itemize}
The polynomials in $\mathbf{Z}\langle x_1,\dots ,x_n\rangle$ of the form $\sum _{\sigma \in \Sigma _n}a(\sigma)x_{\sigma(1)}\cdots x_{\sigma(n)}$, where $\Sigma _n$ denotes the permutation group of $n$ elements and  $a(\sigma)\in \mathbf{Z}$, are said to be multilinear polynomials or multilinear identities of degree $n$. 
\end{Def}

Of special importance in the  theory is the \emph{standard identity}  of degree $n$ 
$$s_n=s_n(x_1,\dots ,x_n)=\sum _{\sigma \in \Sigma _n}\varepsilon (\sigma)x_{\sigma(1)}\cdots x_{\sigma(n)}$$
where $\varepsilon (\sigma)$ denotes the sign of the permutation $\sigma$. For example, if $R$ is a non-zero commutative ring, a theorem of Amitsur and Levitzki shows that $M_n(R)$ satisfies the standard identity $s_m$ for any $m\ge 2n$ and does not satisfy $s_{2n-1}$, cf. \cite[\S 13.3.2 and \S 13.3.2]{MR}. Moreover, if a PI-ring satisfies $s_m$, then it satisfies $s_t$ for any $t\ge m$, cf. \cite[\S 13.3.3]{MR}.

A theorem of Kaplansky shows that if a ring satisfies a monic identity of degree at most $d$, then it satisfies a monic multilinear identity of degree at most $d$ (cf. \cite[\S 13.1.9]{MR}). 

It is possible to define the \emph{PI-degree} of a central simple algebra $\Lambda$ (that is, a simple algebra that is finite dimensional over its center) as the minimum $n$ such that $\Lambda$ satisfies $s_{2n}$ and not $s_{2n-1}$. In this case, we write $\mathrm{PI\mbox{-}deg} (\Lambda)=n$. 

Prime PI-rings $\Lambda$ localized at the non-zero elements of the center are central simple algebras. Then, by definition, the PI-degree of $\Lambda$ is the PI-degree of this localization. For details, check \cite[\S 13.6.7]{MR}.

\begin{Th} \label{starpi} 
Let $\Lambda $ be a right noetherian PI-ring. Then $\Lambda $ satisfies condition $(*)$ and its symmetric $(**)$.
\end{Th}
\begin{Proof} 
Assume the statement is not true. By noetherian induction, we may assume that $\Lambda $ is a right noetherian PI-ring that does not satisfy condition $(*)$, but so does any proper quotient ring of $\Lambda $.

By \cite[Theorem~4]{RS}, $\Lambda $ has only a finite number of idempotent ideals. By Proposition~\ref{reduction}, this implies that $\Lambda $ is a prime, right noetherian, PI-ring. Therefore, $\Lambda $ is (right and left) noetherian by a theorem of G.~Cauchon \cite[13.6.15 Theorem]{MR}. Now we adapt the steps of the proof of  \cite[Theorem~4]{RS} to obtain a contradiction.

Let 
    \begin{equation}\label{eq:starpi}
        \Lambda \supsetneq I_1\supsetneq I_2\supsetneq \cdots \supsetneq I_k \supsetneq \cdots
    \end{equation}
be a sequence of ideals of $\Lambda $ satisfying that $I_{k+1}I_k=I_{k+1}$ for any $k\ge 1$.
Let $P$ be a prime ideal containing $I_k$ for some $k$ and assume that $\mathrm{PI\mbox{-}deg}\, (\Lambda /P) =\mathrm{PI\mbox{-}deg}\, (\Lambda) $. By \cite[Theorem]{small}, the set of elements that are regular modulo $P$ consists of regular elements of $\Lambda$, forms an Ore set $S$, and the localization $\Lambda S^{-1}$ is a local ring with maximal ideal $P\Lambda S^{-1}$. Since $I_k$ is contained in $P$, so is $I_{k+\ell}$ for any $\ell \ge 1$. But in the local ring $\Lambda S^{-1}$,  $P\Lambda S^{-1}I_{k+\ell}\Lambda S^{-1}= I_{k+\ell}\Lambda S^{-1}$ for any $\ell \ge 1$. As the ring is noetherian, Nakayama's lemma implies that $I_{k+\ell}\Lambda S^{-1}=\{0\}$ for any $\ell \ge 1$. Therefore, $I_{k+\ell}=\{0\}$ for any $\ell \ge 1$. This contradicts our hypothesis. Therefore, $\mathrm{PI\mbox{-}deg}\, (\Lambda /P) <\mathrm{PI\mbox{-}deg}\, (\Lambda) $ for any prime ideal $P$ containing $I_k$ for some $k$.

Let $n=\mathrm{PI\mbox{-}deg}\, (\Lambda)$. Let $E(\Lambda )$ be the non-zero ideal generated by the evaluations in $\Lambda $ of $s_{2n-2}$ (cf. \cite[13.6.7 Corollary]{MR}). For any fixed $k$, and  for any prime ideal $P$ containing $I_k$, the ideal $E(\Lambda )$ is contained in $P$ because $\mathrm{PI\mbox{-}deg}\, (\Lambda /P) <\mathrm{PI\mbox{-}deg}\, (\Lambda) $. Thus, $\Lambda/ P$ must satisfy $s_{2n-2}$. Therefore, $E(\Lambda )$ is contained in the prime radical of $I_k$. Then, for any $k$,
there exists $m_k$ such that $E(\Lambda )^{m_k}\subseteq I_k$ (recall that the prime radical of $\Lambda /I_k$ is nilpotent by \cite[Theorem~3.11]{GW}). 

By hypothesis, $\Lambda /E(\Lambda )$ satisfies condition $(*)$. This implies that there exists $k_0$ such that $I_{k_0+\ell}+E(\Lambda )=I_{k_0}+E(\Lambda )$ for any $\ell \ge 0$. By Lemma~\ref{idempotent}, we deduce that $I_{k_0+\ell +1}$ is an idempotent ideal for any $\ell \ge 0$. As a consequence,  the sequence \eqref{eq:starpi} contains infinitely many different idempotent ideals, which gives the desired contradiction with the statement of \cite[Theorem~4]{RS}. This concludes the proof of the statement.
\end{Proof}

\begin{Cor} \label{PImonoid}
Let $\Lambda $ be a right noetherian PI-ring. Then:
\begin{itemize}
    \item[(i)] Every idempotent ideal of $\Lambda $ is the trace of a countably generated projective right $\Lambda $-module and the trace of a countably generated projective left $\Lambda $-module.
    \item[(ii)] $V_r^*(\Lambda)=V_r(\Lambda)\sqcup B_r(\Lambda)\cong V_\ell(\Lambda)\sqcup B_\ell(\Lambda) =V_\ell^*( \Lambda )$.
\end{itemize}
\end{Cor}
\begin{Proof} 
Statement $(i)$ follows from the fact that $\Lambda /N(\Lambda )$, where $N(\Lambda )$ denotes the prime radical of $\Lambda $, is two-sided noetherian \cite[13.6.15]{MR}, and that $N(\Lambda )$ is nilpotent. Therefore, if $I$ is an idempotent ideal of $\Lambda $, then $(I+N(\Lambda))/N(\Lambda )$ is an idempotent ideal of $\Lambda /N(\Lambda )$, so it is the trace of a countably generated projective right $\Lambda /N(\Lambda )$-module and the trace of a countably generated projective left $\Lambda /N(\Lambda )$-module (cf. \cite{whitehead}). By Lemma~\ref{liftingproj}, $I\in \mathcal{T}_r (\Lambda)$, and $I\in \mathcal{T}_\ell (\Lambda)$.

By Theorem~\ref{starpi}, $\Lambda $ and $\Lambda/N(\Lambda )$ satisfy the $(*)$-condition and the symmetric $(**)$-condition. By Theorem~\ref{th:starmonoid}, $V^*(\Lambda/N(\Lambda ))=V(\Lambda/N(\Lambda ))\sqcup B(\Lambda/N(\Lambda ))$, and from Corollary~\ref{isomonoidnil}, we deduce that $V_r^*(\Lambda)= V_r(\Lambda) \sqcup B_r(\Lambda)$ and $V_\ell ^*(\Lambda)= V_\ell(\Lambda) \sqcup B_\ell(\Lambda)$. 

Since, by $(i)$, $\mathcal{T}_r (\Lambda)=\mathcal{T}_\ell (\Lambda)$, the rest of the statement follows from Corollary~\ref{isomonoidtraces}.
\end{Proof}

It is an interesting open question to determine whether there exists a right noetherian ring $\Lambda$ with an idempotent ideal that is the trace of a (countably generated) projective left $\Lambda$-module but not the trace of a projective right $\Lambda$-module. Corollary~\ref{PImonoid} shows that if such an example exists, it cannot be a PI-ring.

The following corollary was our initial motivation for developing these results. As it shows, the basics of the machinery developed in \cite{wiegand} for the endomorphism ring of finitely generated modules over commutative noetherian rings of Krull dimension $1$ are also available in arbitrary Krull dimension. See also \S \ref{subsec:indnoetherian} for the transfer of information between direct sum decompositions of a finitely generated module and a projective module over its endomorphism ring.

\begin{Cor} \label{modulefinite}
Let $R$ be a commutative noetherian ring $R$, and let $\Lambda$ be a module-finite $R$-algebra.  Then $\Lambda $ satisfies condition $(*)$ and its symmetric condition $(**)$. Moreover, it has only a finite number of idempotent ideals, and projective modules satisfy the conclusions of Corollary~\ref{PImonoid}. \end{Cor}

\begin{Proof} If $\Lambda $ is generated by $n$-elements as an $R$-module, then $\Lambda $ is a sub-quotient of $M_n(R)$, so it is PI and noetherian (cf. \cite[13.4.9 Corollary]{MR}). Then, Theorem~\ref{starpi} implies that $\Lambda $ satisfies the $(*)$ condition and the symmetric $(**)$ condition. 

The statement about idempotent ideals is a consequence of \cite[Theorem~4]{RS}.\end{Proof}

Even for a two-sided noetherian ring $\Lambda$, there may be a strict inclusion $V(\Lambda) \sqcup B(\Lambda)\subsetneq V^* (\Lambda)$. See, for example, the last two sections of \cite{traces}.

\subsection{An application to direct sums of finitely generated modules}\label{subsec:indnoetherian}

Let $M$ be a right $\Lambda$-module. By $\add (M)$ we mean the full subcategory of $\Lambda$-modules that are isomorphic to a direct summand of   $M^n$ for some natural number $n$. By $\Add (M)$ we denote the full subcategory of $\Lambda$-modules that are isomorphic to a direct summand of an arbitrary direct sum of copies of $M$. It is also convenient to consider the subcategory $\Add _{\aleph _0} (M)$ of direct summands of $M^{(\omega)}$.

Note that $\Add (\Lambda)$ ($\add (\Lambda)$) is the category of (finitely generated) projective right $\Lambda$-modules. While $\Add _{\aleph _0} (\Lambda)$ is the category of countably generated projective right $\Lambda$-modules.

It is possible to transfer the information about projective modules over the endomorphism ring of a finitely generated module $M$ to $\Add (M)$, via the following well-known result.

\begin{Prop} \label{equivalencia} \cite[Theorem 4.7]{libro}
Let $\Lambda$ be a ring, and let $M$ be a right $\Lambda$-module. Let $S=\End_\Lambda(M)$. Then the functors $\Hom_\Lambda(M,-)$ and $-\otimes_S M$ induce an equivalence between $\add (M)$ and $\add(S)$. 

Moreover, if $M$ is finitely generated, then they induce an equivalence between $\Add (M)$ and $\Add(S)$. This equivalence restricts to an equivalence between $\Add _{\aleph _0} (M)$ and $\Add _{\aleph _0}(S)$.
\end{Prop}

\begin{Cor}
Let $R$ be a commutative noetherian ring, and let $\Lambda$ be a  module-finite $R$-algebra. Let $M$ be a finitely generated right $\Lambda$-module. Then any indecomposable element of $\Add (M)$ is finitely generated. 
\end{Cor}
\begin{Proof}
As $\mathrm{End}_\Lambda (M)$ is a module-finite $R$-algebra, we can apply Corollary~\ref{modulefinite} to deduce that all countably generated projective modules over $\mathrm{End}_\Lambda (M)$ are relatively big. Hence, by Remark~\ref{rem:ind.relbig}, the indecomposable projectives over $\End_\Lambda(M)$ are finitely generated. The conclusion follows from Proposition~\ref{equivalencia}.
\end{Proof}

\begin{Remark} 
In general, if $\Lambda$ is a PI-ring, then the endomorphism ring of a finitely generated right $\Lambda$-module is also a PI-ring (see the argument of \cite[Theorem~2]{PS}). 
\end{Remark}

\section{Projective versus locally projective modules}\label{sec:5} 

Let $R$ be a commutative domain with field of fractions $Q$. An $R$-module $M$ is \emph{torsion-free} if the natural map $M\to M\otimes_R Q$ is injective. 

\begin{notation} \label{not:mq} 
To ease the notation, we will usually write  $M_Q$ for the localization of an $R$-module $M$ at $R\setminus\{0\}$. 
\end{notation}

In the following result, we recall that localization commutes with the $\mathrm{Hom}$ functor for suitable torsion-free modules. This property is crucial for the rest of the paper.

\begin{Lemma}\label{isofg}\cite[Lemma 2.8]{AHP}
Let $R$ be a commutative ring, and let $\Lambda$ be an $R$-algebra. Let $\Sigma$ be a multiplicative subset of $R$ consisting of non-zero divisors. If $M$ is a finitely generated right $\Lambda$-module that is torsion-free as an $R$-module, then the canonical injective homomorphism
    \[\varphi\colon\Hom_\Lambda(M,N)\otimes_R R_\Sigma\to\Hom_{\Lambda\otimes_R R_\Sigma}(M\otimes_R R_\Sigma,N\otimes_R R_\Sigma)\]
is an isomorphism for every $\Lambda$-module $N$ which is torsion-free as an $R$-module.
\end{Lemma}

The following result is an extension of \cite[Lemma 10.1]{wiegand}.

\begin{Lemma}
\label{loc.traces}
Let $R$ be a commutative ring, and let $\Lambda$ be an $R$-algebra. Let $R\to T$ be a flat ring homomorphism.
\begin{enumerate}
    \item[(i)] Let $P$ be a countably generated projective right $\Lambda$-module. Then $P\otimes_R T$ is a countably generated projective right module over $\Lambda\otimes_R T$ and
        \[\Tr_{\Lambda\otimes_R T}(P\otimes_R T)=\Tr_\Lambda(P)\otimes_R T.\]
    \item[(ii)] Let $T$ be a commutative ring. If $M$ is a finitely presented right $\Lambda$-module, then
        \[\Tr_{\Lambda\otimes_R T}(M\otimes_R T)=\Tr_\Lambda(M)\otimes_R T.\]
    \item[(iii)] Let $T=R_\Sigma$, where $\Sigma$ is a multiplicative subset of $R$ consisting of non-zero divisors. Assume $\Lambda$ is torsion-free as an $R$-module. If $M$ is a finitely generated right $\Lambda$-module which is torsion-free as an $R$-module, then 
        \[\Tr_{\Lambda\otimes_R R_\Sigma} (M\otimes_R R_\Sigma)=\Tr_\Lambda(M)\otimes_R R_\Sigma.\]
\end{enumerate}
\end{Lemma}
\begin{Proof}
Statements $(i)$ and $(ii)$ are proved in \cite[Lemma~10.1]{wiegand}. 

To prove $(iii)$, let $x\in \Tr_{\Lambda\otimes_R R_\Sigma}(M\otimes_R R_\Sigma)$. Then there is an element $t\in\Sigma$, and elements $f_1,\dotsc,f_n\in\Hom_{\Lambda\otimes_R R_\Sigma}(M\otimes_R R_\Sigma,\Lambda\otimes_R R_\Sigma)$ and $m_1,\dotsc,m_n\in M$ such that $x=\sum_{i=1}^n f_i\left(m_i\otimes \frac{1}{t}\right)$. Let 
    $$\varphi\colon\Hom_\Lambda(M,\Lambda)\otimes_R R_\Sigma\to\Hom_{\Lambda\otimes_R R_\Sigma} (M\otimes_R R_\Sigma,\Lambda\otimes_R R_\Sigma)$$
denote the canonical homomorphism. By Lemma~\ref{isofg}, $\varphi$ is an isomorphism, so there is an element $r\in\Sigma$ and elements $g_1,\dotsc,g_n\in\Hom_\Lambda(M,\Lambda)$ such that $\varphi(g_i\otimes\frac{1}{r})=f_i$ for each $i=1,\dotsc,n$. Then 
    $$x=\sum_{i=1}^n f_i\left(m_i\otimes \frac{1}{t}\right)=\sum_{i=1}^n \varphi\left(g_i\otimes\frac{1}{r}\right)\left(m_i\otimes \frac{1}{t}\right)=\sum_{i=1}^n g_i(m_i)\otimes\frac{1}{rt}.$$
Therefore, $x\in\Tr_\Lambda(M)\otimes_R R_\Sigma$.

To prove the reverse containment, let $x\in\Tr_\Lambda(M)$ and $t\in R_\Sigma$. Then there are $f_1,\dotsc,f_n\in\Hom_\Lambda(M,\Lambda)$, $m_1,\dotsc,m_n\in M$ such that $x=\sum_{i=1}^n f_i(m_i)$. Therefore, $x\otimes t=\sum_{i=1}^n f_i(m_i)\otimes t=\sum_{i=1}^n \varphi(f_i\otimes 1)(m_i\otimes t)$. Hence, $x\otimes t\in \Tr_{\Lambda\otimes_R R_\Sigma}(M\otimes_R R_\Sigma)$. This concludes the proof of the statement.
\end{Proof}

It is easy to prove that the notion of relatively big projective module behaves well with respect to localization.

\begin{Cor} \label{relbig_loc} 
Let $\Lambda$ be an algebra over a commutative ring $R$. Let $J\in \mathcal{T}(\Lambda)$. Let $\Sigma$ be a multiplicatively closed subset of $R$ consisting of non-zero divisors. If $P$ is a relatively $J$-big projective right $\Lambda$-module, then $P_\Sigma$ is a relatively $J_\Sigma$-big projective right $\Lambda _\Sigma$-module.
\end{Cor}
\begin{Proof}
By Lemma~\ref{loc.traces}, $J_\Sigma \in \mathcal{T}(\Lambda _\Sigma)$. Since $P/PJ$ is finitely generated, so is $(P/PJ)_\Sigma\cong P_\Sigma /P_\Sigma J_\Sigma$. 

Let $P'$ be a countably generated projective right $\Lambda$-module with trace $J$. Since $P\cong P\oplus Q^{(\omega)}$ as $\Lambda$-modules, also $P_\Sigma\cong P_\Sigma\oplus (P')_\Sigma^{(\omega)}$ as $\Lambda_\Sigma$-modules. This finishes the proof of the claim.
\end{Proof}

The following result is necessary to provide a version of the Package Deal Theorem for projective modules.

\begin{Cor}\label{loc.proj}
\cite[Corollary 4.3]{AHP}
Let $R$ be a commutative domain, and let $\Lambda$ be an $R$-algebra. Let $N$ be a finitely generated right $\Lambda$-module which is torsion-free as an $R$-module. Then $N$ is a projective $\Lambda$-module if and only if $N_\fm$ is a projective $\Lambda_\fm$-module for every maximal ideal $\fm$ of $R$.
\end{Cor}

\begin{Def}
A commutative domain $R$ is said to have \emph{finite character} if any non-zero element is contained in only a finite number of maximal ideals. A commutative domain of finite character is said to be \emph{$h$-local} if, in addition, any non-zero prime ideal of $R$ is contained in a unique maximal ideal.
\end{Def}

The domain $R$ has finite character if $R/I$ is a semilocal ring for every non-zero ideal $I$ of $R$. If $R$ is $h$-local, then it also satisfies that  $R/\mathfrak{p}$ is a local domain for every non-zero prime ideal $\mathfrak{p}$ of $R$. Equivalently, $R$ is an $h$-local domain if and only if $R/I$ is semiperfect for every non-zero ideal $I$ of $R$. For characterizations of $h$-local domains, the reader is referred to \cite{olberding}.

In general, one has the following easy property for modules over torsion-free algebras over a domain.

\begin{Lemma}\label{d}
Let $R$ be a domain, and let $\Lambda$ be an $R$-algebra. Let $\Sigma$ be a multiplicative subset of $R$. Let $U$ be a right $\Lambda$-module which is torsion-free as an $R$-module. If $V$ is a finitely generated $\Lambda$-module that is torsion-free and such that  $V_\Sigma$ is isomorphic to a submodule of $U_\Sigma$, then $U$ has a submodule isomorphic to $V$.
\end{Lemma}
\begin{Proof} By hypothesis, and because $V$ is torsion-free, $U_\Sigma$ has a submodule $V'$ that is isomorphic to $V$. Since $V'$ is finitely generated, there exists $d\in \Sigma$ such that $V'\subseteq \frac 1d U\subseteq U_\Sigma$. Then $dV' \leq U$ is the module we are looking for.
\end{Proof}

If in Lemma~\ref{d}, $R$ has finite character then the element $d$ that appears in the proof is contained only in a finite number of maximal ideals. So one easily proves results like the following one, which are crucial in our discussion. 

\begin{Lemma}\label{almostallmaximals} \cite[Lemma~4.2 and 5.1]{AHP}
Let $R$ be a commutative domain of finite character with field of fractions $Q$. Let $\Lambda$ be an $R$-algebra, and let $M,N$ be  right $\Lambda$-modules which are torsion-free as $R$-modules. Assume that $M_Q\cong N_Q$. Then 
\begin{itemize}
    \item [(i)] If $M$ and $N$ are finitely generated, then $M_\fm\cong N_\fm$ for almost all maximal ideals $\fm$ of $R$.
    \item[(ii)] If $M$ is finitely generated and  $N\subseteq M$ then $M_\fm = N_\fm$ for almost all maximal ideals $\fm$ of $R$. So that,  there exist $\fm_1,\dots ,\fm _\ell$ maximal ideals of $R$ such that  $M /N=\prod _{i=1}^\ell M _{\fm _i}/N_{\fm _i}$.
\end{itemize}
\end{Lemma}

In the following result, we specialize in two-sided ideals of $\Lambda$ and assume that $\Lambda _Q$ is a simple artinian ring.

\begin{Lemma}
\label{almosttrace}
Let $R$ be an $h$-local domain with field of fractions $Q$. Let $\Lambda$ be a torsion-free $R$-algebra such that $\Lambda_Q$ is simple artinian, and let $I$ be a non-zero two-sided ideal of $\Lambda$. Then:
\begin{enumerate}
    \item[(i)] $I\cap R\neq \{0\}$;
    \item[(ii)] $I_\fm=\Lambda _\fm$ for almost all maximal ideals $\fm$ of $R$;
    \item[(iii)] There exist $\fm_1,\dots ,\fm _\ell$ maximal ideals of $R$ such that  $\Lambda /I=\prod _{i=1}^\ell \Lambda _{\fm _i}/I_{\fm _i}$.
\end{enumerate}
\end{Lemma}

\begin{Proof} $(i)$ Since  $\Lambda _Q$ is simple artinian and $\Lambda$ is torsion-free, $I_Q$ is a non-zero two-sided ideal of $\Lambda _Q$. So $I_Q=\Lambda _Q$. Therefore, there exist $x_1,\dots ,x_n \in I$ and $d\in R\setminus \{0\}$ such that $1=\sum _{i=1}^n\frac{x_i}d$. Thus $d\in R\cap I$ and this finishes the proof of $(i)$. Note that $d$ is contained only in a finite number of maximal ideals, so this yields $(ii)$.

Statement $(iii)$ is a direct  application of Lemma~\ref{almostallmaximals}.
\end{Proof}

The following results, which are known as Package Deal Theorems, were proved in \cite{levyodenthal2} for the case of algebras over noetherian rings of Krull dimension 1. The name of these theorems comes from the idea of having a ``package'' of modules over a localization and finding a module over the non-localized algebra whose localizations are those given in the ``package''. In \cite{AHP} Package Deal Theorems were presented for the localizations of submodules and the localizations of trace ideals. We recall the statements here.

\begin{PDTh}[Localization of submodules] \label{dealsubmodules}
\cite[Theorem 5.7]{AHP}
Let $R$ be an $h$-local domain with field of fractions $Q$. Let $\Lambda$ be an $R$-algebra, and let $M$ be a finitely generated right $\Lambda$-module, which is torsion-free as an $R$-module. For each maximal ideal $\fm$ of $R$, let $X(\fm)$ be a $\Lambda_\fm$-submodule of $M_\fm$, which is torsion-free as an $R_\fm$-module, and such that $X(\fm)_Q=(M_\fm)_Q$. Then the following statements are equivalent
\begin{enumerate}
    \item[(i)] There is a $\Lambda$-submodule $N\subseteq M$, which is torsion-free as an $R$-module, and such that $N_\fm=X(\fm)$ for all maximal ideals $\fm$ of $R$.
    \item[(ii)] $X(\fm)=M_\fm$ for almost all maximal ideals $\fm$ of $R$.
\end{enumerate}
Moreover, if each of the $X(\fm)$ is finitely generated as a $\Lambda_\fm$-module, then $N$ is also finitely generated as a $\Lambda$-module.
\end{PDTh}

\begin{PDTh}[Localization of trace ideals] \label{dealtraces}
\cite[Theorem 5.13]{AHP}
Let $R$ be an $h$-local domain with field of fractions $Q$. Let $\Lambda$ be a torsion-free $R$-algebra such that $\Lambda_Q$ is a simple artinian ring. For each maximal ideal $\fm$ of $R$, let $I(\fm)$ be a non-zero two-sided ideal of $\Lambda_\fm$ which is the trace ideal of a countably generated projective right $\Lambda_\fm$-module. Then the following statements are equivalent
\begin{enumerate}
    \item[(i)] There is a two-sided ideal $I$ of $\Lambda$, which is the trace ideal of a countably generated projective right $\Lambda$-module, and such that $I_\fm=I(\fm)$ for all maximal ideals $\fm$ of $R$.
    \item[(ii)] $I(\fm)=\Lambda_\fm$ for almost all maximal ideals $\fm$ of $R$.
\end{enumerate}
Moreover, the ideal $I$ that satisfies the equivalent conditions $(i)$ and $(ii)$ is unique.
\end{PDTh}

\begin{Cor}
\label{cor:freelmst}
Let $R$ be a commutative domain of finite character with field of fractions $Q$. Let $\Lambda$ be a torsion-free $R$-algebra such that $\Lambda _Q$ is simple artinian. Let $S$ denote a simple right  $\Lambda _Q$-module, and assume that $\Lambda_Q\cong S^k$ as right $\Lambda_Q$-modules. Then any finitely generated projective right $\Lambda$-module $P$ satisfies that $P_\fm^k$ is free for almost all maximal ideals $\fm$ of $R$.
\end{Cor}
\begin{Proof} 
By assumption $P_ Q\cong S^\ell$ for suitable $\ell$. Therefore $P_Q^{k}\cong \Lambda _Q^\ell$. Now the conclusion follows from Lemma~\ref{almostallmaximals} using the fact that finitely generated projective modules over $\Lambda$ are torsion-free because $\Lambda$ is torsion-free.
\end{Proof}

\begin{PDTh}[Localization of Projective Modules]
\label{dealprojective}
Let $R$ be an $h$-local domain with field of fractions $Q$. Let $\Lambda$ be a torsion-free $R$-algebra such that $\Lambda_Q$ is simple artinian. For each maximal ideal $\fm$ of $R$, let $P(\fm)$ be a finitely generated projective right $\Lambda_\fm$-module with trace ideal $I(\fm)$, and such that $P(\fm)_Q\cong \Lambda_Q^{r}$ for some positive integer $r$. Then the following statements are equivalent:
\begin{enumerate}
    \item[(i)] There is a finitely generated projective right $\Lambda$-module $P$ such that $P_\fm\cong P(\fm)$ for every maximal ideal $\fm$ of $R$.
    \item[(ii)] $P(\fm)\cong \Lambda_\fm^{r}$ for almost all maximal ideals $\fm$ of $R$.
\end{enumerate}
Moreover, if $P$ is a finitely generated projective satisfying $(i)$ and $I = \Tr_{\Lambda}(P)$, then 
$I_{\fm} = I(\fm)$ for every $\fm \in \mSpec\ R$.
\end{PDTh}
\begin{Proof}
$(ii)\Rightarrow(i)$. Let $\mathcal M$ be the finite set of maximal ideals of $R$ such that $P(\fm)\not\cong\Lambda_\fm^{r}$. Consider a module $F = \Lambda^{r}$. 
Since $\Lambda$ is torsion-free as an $R$-module, for every $\fm \in \mathcal M$ there exists 
a submodule $F(\fm)$ of $F_{\fm}$ such that $F(\fm) \cong P(\fm)$. For $\fn \in \mSpec(R) \setminus \mathcal M$ set $F(\fn) = F_{\fn}$.
By Package Deal Theorem \ref{dealsubmodules}, there is a finitely generated right $\Lambda$-module $P$ contained in $F$ such that $P_\fm = F(\fm) \cong P(\fm)$ for every $\fm\in\mathcal \mSpec(R)$.
By Corollary \ref{loc.proj}, $P$ is a projective module.
Also note that $\Tr_{\Lambda}(P)_{\fm} = \Tr_{\Lambda_{\fm}}(P_{\fm})$ by Lemma~\ref{loc.traces}.

$(i)\Rightarrow (ii)$ Apply Lemma~\ref{almostallmaximals}.
\end{Proof}

\begin{Cor}
\label{cor:dealprojective}
Let $R$ be an $h$-local domain with field of fractions $Q$. Let $\Lambda$ be a torsion-free $R$-algebra such that $\Lambda_Q$ is simple artinian. Let $S$ denote a simple right $\Lambda_Q$-module, and assume that $\Lambda_Q\cong S^k$ as right $\Lambda_Q$-modules. For each maximal ideal $\fm$ of $R$, let $P(\fm)$ be a finitely generated projective right $\Lambda_\fm$-module with trace ideal $I(\fm)$. Assume there exists a positive 
integer $r$ such that $P(\fm)_Q\cong S^{r}$. Then the following statements are equivalent:
\begin{enumerate}
    \item[(i)] There is a finitely generated projective right $\Lambda$-module $P$ with trace ideal $I$ such that $P_\fm\cong P(\fm)^k$ and $I_\fm=I(\fm)$ for every maximal ideal $\fm$ of $R$.
    \item[(ii)] $P(\fm)^k$ is free for almost all maximal ideals $\fm$ of $R$.
\end{enumerate}
\end{Cor}
\begin{Proof}
Apply Theorem~\ref{dealprojective} to a package $P(\fm)':=P(\fm)^k$, for each maximal ideal $\fm$ of $R$.
\end{Proof}

\section{Semiperfect rings and locally semiperfect rings}\label{sec:6}

We start this section by describing the monoids of projective modules for  semiperfect rings. We follow the notation and conventions of \cite[\S 27]{andersonfuller}.

\subsection{Projective modules over semiperfect rings and their trace ideals}

\begin{Lemma} \label{lifting}
Let $\Lambda$ be a semiperfect ring, and let $I$ be an ideal of $\Lambda$. Let $P'$ be a finitely generated projective right $\Lambda /I$-module. Then, there exists a finitely generated projective right $\Lambda$-module $X$ such that $X/XI\cong P'$.
\end{Lemma}
\begin{Proof}
Since $\Lambda$ is semiperfect, $P'$ has a projective cover. That is, it fits in an exact sequence
    \[\begin{tikzcd}0\rar & \Ker f\rar & X \rar{f} & P' \rar & 0,\end{tikzcd}\]
where $X$ is a finitely generated projective right $\Lambda$-module, and $\Ker f$ is a small submodule of $X$. Therefore, there are inclusions
    $$XI\subseteq \Ker f \subseteq XJ(\Lambda)$$
which induce an exact sequence of $\Lambda /I$-modules
    \[\begin{tikzcd}0\rar & \Ker f/XI\rar & X/XI\rar{f'} & P'\rar & 0.\end{tikzcd}\]
which splits because $P'$ is projective. Since $\Ker f/XI$ is also a small submodule of $X/XI$, it must be zero. Therefore $X/XI\cong P'$.
\end{Proof}

In the next proposition, we recall the behavior of traces of projective modules for semiperfect rings. 

\begin{Prop} \label{tracessemiperfect} 
Let $\Lambda$ be a semiperfect ring, and let $J(\Lambda)$ denote its Jacobson radical. Then:
\begin{enumerate}
    \item[(i)] $\mathcal{T}_r(\Lambda)=\mathcal{T}_\ell(\Lambda) : = \mathcal{T} (\Lambda)$, and the injective map $\varphi\colon \mathcal{T} (\Lambda)\to  \mathcal{T} (\Lambda/J(\Lambda))$ defined in Lemma~\ref{tracesmodjradical} is bijective.
    
    \item[(ii)] Let $e_1, \dots ,e_t$ be a basic set of (primitive orthogonal) idempotents of $\Lambda$. Let $\mathcal{P} (\{1,\dots ,t\})$ denote the set of all subsets of $\{1,\dots ,t\}$. Then the map
        $$\Gamma\colon \mathcal{P} (\{1,\dots ,t\}) \to \mathcal{T} (\Lambda)$$ 
    defined by $\Gamma (A)=\sum _{i\in A}\Lambda e_i\Lambda$, for any $A\in \mathcal{P} (\{1,\dots ,t\})$, is a bijection.
    
    \item[(iii)] Let $P$ be a projective right $\Lambda$-module with trace ideal $I$. Assume that $I/(I\cap J(\Lambda))\cong \sum _{i\in A} \left( (e_i+J(\Lambda))\Lambda/J(\Lambda)\right)^{\ell _i} $ for suitable $\ell _i>0$ and $A\in \mathcal{P} (\{1,\dots ,t\})$. Then $P\cong \bigoplus _{i\in A} \left(e_i\Lambda \right)^{(L_i)}$, for suitable non-empty sets, $L_i$.
    
    \item[(iv)] Keeping the notation from $(ii)$ and $(iii)$, $\Lambda':=\Lambda/I$ is a semiperfect ring and $e_i' := e_i + I, i \not \in A$ is a basic set of   idempotents of $\Lambda/I$. Moreover, let $P' \cong \bigoplus_{i \not \in A} (e_i'\Lambda')^{x_i}$ be a finitely generated projective $\Lambda'$-module, and let $P$ be a finitely generated projective $\Lambda$-module such that $P/PI \cong P'$. Then there are non-negative integers $s_i, i \not \in A$ such that $P \cong \bigoplus_{i \not \in A} (e_i\Lambda)^{x_i} \oplus \bigoplus_{i  \in A} (e_i\Lambda)^{s_i}$.
\end{enumerate}
\end{Prop}
\begin{Proof} 
The proof of the four statements is more or less routine. One just has to take into   account that any projective right (left) $\Lambda$-module is isomorphic to a unique direct sum of the indecomposable projective modules $P_i=e_i\Lambda$ ($Q_i=\Lambda e_i$), for $i=1,\dots ,t$ (see, for example, \cite[Theorem 3.10]{libro}). Moreover $\mathrm{Tr}_\Lambda (P_i)=\mathrm{Tr}_\Lambda (Q_i)=\Lambda e_i\Lambda$.
\end{Proof}

To give explicit descriptions of the monoids of projective modules over certain classes of rings, we introduce the monoid/semiring $\No ^*$.

\begin{notation}\label{not:nostar}
Let  $\N_0 =\{0,1,2,\dots\}$, and let $\N =\{1,2,\dots \}$. We denote by $\No ^*=\N_0\cup \{\infty\}$ the monoid with the sum 
extending the one of $\No$ and with the additional rule that $x+\infty=\infty +x =\infty$ for any $x\in \No ^*$. We also see $\No ^*$ as a semiring with the product extending the one of $\No$ together with
  the rules $x\cdot \infty=\infty \cdot x =\infty$ for any $x\in \No ^*\setminus \{0\}$, and $0\cdot\infty = \infty\cdot 0 = 0$.
\end{notation}

\begin{Prop} \label{monoidsemiperfect}
Let $\Lambda$ be a semiperfect ring. Then:
\begin{itemize}
    \item[(i)] The canonical isomorphism between $V (\Lambda):=V_r(\Lambda)\cong \No ^t$ and $V_\ell (\Lambda)$ extends to an isomorphism between $V^*(\Lambda):=V_r^*(\Lambda) \cong (\No ^*)^t$ and $V_\ell^*(\Lambda)$.
    \item[(ii)] Let $P$ be a projective right $\Lambda$-module. Then there exists a two-sided ideal $I$ of $\Lambda$ that is minimal with respect to the property that $P/PI$ is finitely generated. Moreover, $I$ is the trace of a finitely generated projective right module.
    \item[(iii)] Every countably generated projective module is relatively big, so that $V^*(\Lambda)=V(\Lambda)\sqcup B(\Lambda)$.
\end{itemize}
\end{Prop}
\begin{Proof}
Let $e_1, \dots ,e_t$ be a basic set of idempotents of $\Lambda$. Then $P_i=e_i\Lambda$ and $Q_i=\Lambda e_i\cong \mathrm{Hom}_\Lambda (P_i, \Lambda)$, for $i=1,\dots ,t$, are representatives of the isomorphism classes of indecomposable projective right and left $\Lambda$-modules. Recall that any projective right or left module over $\Lambda$ is a direct sum of indecomposables, uniquely up to isomorphism.

$(i)$  Let $P$ be a countably generated projective right $\Lambda$-module. Then there exist $A_1,\dots ,A_t$ countable sets such that $P\cong \bigoplus _{i=1}^tP_i^{(A_i)}$. 

Then the assignment
\[[P]=\left[\bigoplus _{i=1}^tP_i^{(A_i)}\right]\mapsto \left[\bigoplus _{i=1}^tQ_i^{(A_i)}\right]\]
gives an isomorphism of monoids $V_r^*(\Lambda)\to V_\ell^*(\Lambda)$ that extends the canonical isomorphism $V_r(\Lambda)\to V_\ell(\Lambda)$ from Remark~\ref{rem:fgproj}.

Moreover, the assignment $[P]\mapsto (n_1,\dots ,n_t)$, where for $i=1,\dots, t$ 
$$n_i=\begin{cases}
    |A_i| & \text{if $A_i$ is finite} \\
    \infty & \text{otherwise},
\end{cases}$$
induces an isomorphism $V_r^*(\Lambda)\to (\No ^*)^t$ that restricts to an isomorphism $V_r(\Lambda) \to \No ^t$.

$(ii)$ If $P$ is a projective right $\Lambda$-module, then $\{1,\dots, t\}=A\cup B$, where $A$ and $B$ are disjoint, and $P\cong \left(\bigoplus _{i\in A}P_i^{(A_i)}\right)\oplus \left(\bigoplus _{i\in B}P_i^{(B_i)}\right)$, where $A_i$ is a finite set (possibly empty) for any $i\in A$, and $B_i$ is an infinite set for any $i\in B$. Then, by Lemma~\ref{ds.traces}, $$I=\sum _{i\in B} \mathrm{Tr}_\Lambda (P_i)=\mathrm{Tr}_\Lambda \Big(\bigoplus _{i\in B}P_i\Big)$$  is the ideal we are looking for.

$(iii)$ By Proposition~\ref{tracessemiperfect}, $\mathcal{T}_r(\Lambda)=\mathcal{T}_\ell(\Lambda)$, so we can apply Corollary~\ref{isomonoidtraces} to deduce that $B_r(\Lambda):=B(\Lambda)$ is canonically isomorphic to $B_\ell (\Lambda)$. 

Statement $(ii)$ implies that any countably generated projective module is relatively big. Hence, $V^*(\Lambda)=V(\Lambda)\sqcup B(\Lambda)$.
\end{Proof}

\subsection{Locally semiperfect algebras over \texorpdfstring{$h$}{h}-local domains}\label{subsec:ls}

\begin{Def}
Let $\Lambda$ be an algebra over a commutative domain $R$. We say that $\Lambda$ is \emph{locally semiperfect} if $\Lambda_\fm$ is semiperfect for every maximal ideal $\fm$ of $R$.
\end{Def}

From now on we are using the conventions of Notation~\ref{not:mq}. 

The following example explains a particular  instance of the appearance of locally semiperfect algebras that, in addition, satisfy that they have a simple artinian ring of quotients. With the additional hypothesis that the domain $R$ is $h$-local, they will be our main object of study in   Section~\ref{sec:9}.

\begin{Ex} \label{ex:locsemi}
Let $R$ be a commutative domain with field of fractions $Q$, and let $M$ be a finitely generated torsion-free $R$-module. Then $\Lambda = \End_{R} (M)$ is an $R$-algebra  and, by Lemma~\ref{isofg}, $\Lambda _Q= \Lambda \otimes _R Q  \cong \End_{Q}(M\otimes _R Q)=\End_{Q}(M_Q)$ is simple artinian. By \cite[Proposition~3.14]{libro} and Lemma~\ref{isofg}, for every maximal ideal $\fm$ of $R$, $M_{\fm}$ is a direct sum of modules with local endomorphism ring  if and only if $\Lambda_{\fm} \cong \End_{R_\fm}(M_{\fm})$ is a semiperfect ring. In other words, $\Lambda$ is locally semiperfect if and only if  $M_{\fm}$ is a direct sum of modules with local endomorphism ring for every maximal ideal $\fm$ of $R$.
\end{Ex}

\begin{Lemma}
\label{semiperf.fgtrace}
Let $R$ be an $h$-local domain with field of fractions $Q$, and let $\Lambda$ be a locally semiperfect torsion-free $R$-algebra such that $\Lambda_Q$ is simple artinian. Then every trace of a projective right (or left) $\Lambda$-module is also the trace of a finitely generated projective right $\Lambda$-module. In particular, $\mathcal{T}_r(\Lambda)=\mathcal{T}_\ell (\Lambda)$.
\end{Lemma}
\begin{Proof}
Let $I$ be a non-zero two-sided ideal of $\Lambda$ which is the trace of a projective right $\Lambda$-module. Since $\Lambda$ is torsion-free, $I_\fm\neq\{0\}$ holds for every maximal ideal $\fm$ of $R$. By Lemma~\ref{loc.traces} and Lemma~\ref{tracessemiperfect},  $I_\fm$ is also the trace of a finitely generated projective right $\Lambda_\fm$-module. 

Let $S$ denote a simple right $\Lambda_Q$-module, and assume that $\Lambda_Q\cong S^k$ as a right $\Lambda_Q$-module. Since $\Lambda_Q$ is simple artinian and $I$ is non-zero, by Lemma~\ref{almosttrace}, there exists a finite set $\mathcal M$ of maximal ideals of $R$ such that $I_\fm = \Lambda_\fm$ if $\fm \not \in \mathcal M$. For any maximal ideal $\fm\in\mathcal M$, let $P(\fm)$ be a finitely generated projective right $\Lambda_\fm$-module with trace ideal $I_\fm$. Then $P(\fm)_Q\cong S^{t_\fm}$ for a suitable positive integer $t_\fm\ge1$. Let $t=\mathrm{lcm}\{t_\fm\mid \fm\in\mathcal M\}$, write $t=\ell_\fm\cdot t_\fm$, and let $P'(\fm):=P(\fm)^{k\ell_\fm}$ for each maximal ideal $\fm$ of $\mathcal M$. For each $\fm\notin\mathcal M$, take $P'(\fm)=\Lambda_\fm^t$. Note that $P'(\fm)_Q\cong S^{tk}\cong\Lambda_Q^t$ and $\Tr_{\Lambda_\fm}(P'(\fm))=I_\fm$ hold for every maximal ideal $\fm$ of $R$. By Package Deal Theorem~\ref{dealprojective} applied to the package $\{P'(\fm)\mid \fm\in\mSpec(R)\}$, there exists a finitely generated projective right $\Lambda$-module $P$ such that $P_\fm \cong P'(\fm)$ for any maximal ideal $\fm$ of $R$. By Lemma~\ref{loc.traces}, $\Tr _\Lambda (P)_\fm =\Tr _{\Lambda _\fm} (P'(\fm)) =I_\fm$ for any maximal ideal $\fm$ of $R$. This implies that  $I= \Tr _\Lambda (P)+I=\Tr _\Lambda (P)$.

A similar argument gives the conclusion for the traces of projective left $\Lambda$-modules due to the symmetry of the hypotheses.

The equality $\mathcal{T}_r(\Lambda)=\mathcal{T}_\ell (\Lambda)$ follows from the canonical bijection between finitely generated projective right and left $\Lambda$-modules induced by the $\Lambda$-dual, and because   the trace of a finitely generated projective module  coincides with the trace of its $\Lambda$-dual, cf. Remark~\ref{rem:fgproj}.
\end{Proof}

\begin{Cor} \label{omega}
Let $R$ be an $h$-local domain with field of fractions $Q$, and let $\Lambda$ be a locally semiperfect torsion-free $R$-algebra such that $\Lambda_Q$ is simple artinian. If $P$ is a countably generated projective right $\Lambda$-module, then $P^{(\omega)}$ is a direct sum of finitely generated modules.
\end{Cor}
\begin{Proof}
By Lemma~\ref{semiperf.fgtrace}, the trace of $P$ coincides with the trace of a finitely generated projective module. To conclude, apply that if $P$ and $Q$ are  countably generated projective modules then $P^{(\omega)}\cong Q^{(\omega)}$ if and only if $P$ and $Q$ have the same trace. 
\end{Proof}

\begin{Lemma} \label{finitely_many} 
Let $R$ be a commutative domain of finite character with field of fractions $Q$. Let $\Lambda$ be a locally semiperfect torsion-free $R$-algebra such that $\Lambda _Q$ is simple artinian; so that there exist $k$ and a division ring D such that $\Lambda _Q\cong M_k(D)$. Then, for almost all maximal ideals $\fm$ of $R$, the ring $\Lambda_\fm$ has the following properties:
\begin{enumerate}
    \item[(i)] There is exactly one indecomposable projective $\Lambda_\fm$-module up to isomorphism. Moreover, $\Lambda_\fm\cong M_k(T_\fm)$, where $T_\fm$ is a local ring such that $T_\fm\otimes_RQ\cong D$.
    \item[(ii)] If  $U$ and $V$ are  projective $\Lambda_\fm$-modules, then $U \cong V$ if and only if $U_Q \cong V_Q$.
\end{enumerate}
\end{Lemma}
\begin{Proof}
Let $S$ be a simple $\Lambda_Q$-module and assume that $\Lambda_Q \cong S^{k}$. Let $X$ be a finitely generated $\Lambda$-module, which is torsion-free as an $R$-module, and such that $X_Q \cong S$. By Lemma~\ref{almostallmaximals}, $X_\fm^{k}$ is free for almost all maximal ideals. If $\fm$ is such a maximal ideal, then $X_\fm^{k} \cong \Lambda_\fm$ gives an indecomposable decomposition of the semiperfect ring $\Lambda_\fm$. Hence $X_\fm$ is the only indecomposable projective $\Lambda_\fm$-module up to isomorphism. As $\Lambda_\fm$ is semiperfect, the endomorphism ring of $X_\fm$ is a local ring that we will denote by $T_\fm$. Hence $\Lambda _\fm \cong \mathrm{End}_{\Lambda _\fm}(X_\fm^{k})\cong M_k (T_\fm).$ By Lemma~\ref{isofg}, $T_\fm\otimes_RQ\cong D$.

Projective modules over $M_k(T_\fm)$ are always a direct sum of copies of the indecomposable projective $X_\fm$. Moreover, $U \cong X_\fm^{(I)}, V \cong X_\fm^{(J)}$ then $U\cong V$ if and only if $I$ and  $J$ have the same cardinality, if and only if $U_Q \cong V_Q$.
\end{Proof}

\begin{Remark} \label{descriptiontracesalmostsemiperfect} 
Let $R$ be an $h$-local domain with field of fractions $Q$, and let $\Lambda$ be a locally semiperfect torsion-free $R$-algebra such that $\Lambda _Q$ is simple artinian, so that there exist $k$ and a division ring D such that $\Lambda _Q\cong M_k(D)$. In view of Lemma~\ref{finitely_many}, the maximum spectrum of $R$ can be written as a disjoint union of two  sets $A$ and $B$ such that for any $\fm \in A$, $\Lambda _\fm$ is a ring of $k\times k$ matrices over a local ring and $B$ is a finite, non-empty set. In particular, for $\fm \in A$, $\mathcal{T}_r (\Lambda _\fm )=\mathcal{T}_\ell (\Lambda _\fm )=\{0, \Lambda _\fm \}$.

In view of Package Deal Theorem~\ref{dealtraces} and Lemma~\ref{semiperf.fgtrace}, there is a bijective correspondence
    \[\begin{tikzcd}
        \mathcal{T}_r (\Lambda )\setminus \{0\}=\mathcal{T}_\ell (\Lambda  )\setminus \{0\}\rar & \prod _{\fm \in B}\left( \mathcal{T}_r (\Lambda _\fm)\setminus \{0\}\right).
    \end{tikzcd}\]
\end{Remark}

\begin{Cor} \label{minimal} 
Let $R$ be an $h$-local domain with field of fractions $Q$, and let $\Lambda$ be a locally semiperfect torsion-free $R$-algebra such that $\Lambda _Q$ is simple artinian. Then $\mathcal{T} (\Lambda)$ is finite. Moreover, if $P$ is a finitely generated projective $\Lambda$-module, then $I = \Tr_{\Lambda}(P)$ is a minimal non-zero element of $\mathcal{T} (\Lambda)$ if and only if $P_{\fm}$ is a finite power of an indecomposable projective $\Lambda_{\fm}$-module
for any $\fm \in \mSpec(R)$.
\end{Cor}
\begin{Proof}
Finiteness of ${\mathcal  T}(\Lambda)$ follows from Remark~\ref{descriptiontracesalmostsemiperfect}. Also, note that the correspondence of the remark preserves inclusions, 
so $\Tr_{\Lambda}(P)$ is a minimal element of ${\mathcal T}(\Lambda)\setminus \{0\}$ if and only if $\Tr_{\Lambda_{\fm}}(P_{\fm})$ is a minimal element of ${\mathcal T}(\Lambda_{\fm}) \setminus \{0\}$ for any
$\fm \in \mSpec(R)$ if and only if $P_{\fm}$ is a finite power of an indecomposable projective $\Lambda_{\fm}$-module.
\end{Proof}

\begin{Prop} \label{minimallocallysemiperfect}
Let $R$ be an $h$-local domain with field of fractions $Q$, and let $\Lambda$ be a locally semiperfect torsion-free $R$-algebra such that $\Lambda _Q$ is simple artinian. Let $P$ be a countably generated projective right $\Lambda$-module. Then,
\begin{enumerate}
    \item[(i)] $P$ is finitely generated if and only if $P_Q$ is a finitely generated $\Lambda _Q$-module if and only if there exists a maximal ideal $\fm$ of $R$ such that $P_\fm$ is a finitely generated $\Lambda _\fm$-module.

    \item[(ii)] The set $\mathcal{I} (P)=\{ J\unlhd\Lambda\mid P/PJ \mbox{ is finitely generated}\}$ has a minimal element, which is the trace of a finitely generated, projective right $\Lambda$-module. If $I$ is such ideal, then it is characterized by the property that, for any maximal ideal $\fm$ of $R$, $I_\fm$ is the minimal element of $\mathcal{I}(P_\fm)=\{ J\unlhd\Lambda\mid P_\fm/P_\fm J \mbox{ is finitely generated}\}$.
\end{enumerate}
\end{Prop}
\begin{Proof} 
Throughout the proof, let $S$ denote a simple right $\Lambda _Q$-module.

$(i)$ We claim that if $P_Q\cong S^n$ for some $n\in \N$, then $P$ is finitely generated. Indeed, by hypothesis, there exists a finitely generated $\Lambda$-submodule $P_0$ of $P$ such that $(P_0)_Q=P_Q$. Then we may assume that we have $P_0\le P\le \Lambda^{(\omega)}=P\oplus P'$, and for any $p\in P$ there exists $0\neq d\in R$ such that $pd\in P_0$. Since $\Lambda$ is torsion-free as an $R$-module, this implies that there exists a finite subset $F$ of $\omega$ such that $P_0 \le P\subseteq \Lambda ^F$. Hence, $\Lambda ^F=P\oplus (P'\cap \Lambda ^F)$ so $P$ is finitely generated. 

This proves statement $(i)$.

$(ii)$ In the proof, we maintain the notation from Remark~\ref{descriptiontracesalmostsemiperfect} regarding the partition of the maximal spectrum of $R$. 

If $P$ is finitely generated, the statement is clear. Assume   that $P$ is a countably generated, projective, right $\Lambda$-module that is not finitely generated. By $(1)$, $P_Q\cong S^{(\omega)}$. Therefore, there are two possible situations regarding the localizations of $P$ at the maximal ideals of $R$:
\begin{enumerate}
    \item [(a)] $P_\fm$ is a free $\Lambda _\fm$-module of infinite rank for any $\fm \in A$. Therefore, $\Lambda _\fm$ is the only element of the set $\mathcal{I} (P_\fm)=\{J\unlhd\Lambda _\fm\mid P_\fm/P_\fm J \mbox{ is finitely generated}\}$. 
    \item[(b)] $P_\fn$ is an infinitely generated projective module over the semiperfect ring $\Lambda _\fn$ for any $\fn \in B$. By Proposition~\ref{monoidsemiperfect}, there exists a non-zero ideal $I(\fn)$ that is the minimal element of the set $\mathcal{I}(P_\fn)=\{J\unlhd\Lambda _\fn\mid P_\fn/P_\fn J \mbox{ is finitely generated}\}$. In addition, $I(\fn)$ is the trace of a finitely generated projective right $\Lambda _\fn$-module. 
\end{enumerate}   

By Package Deal Theorem~\ref{dealtraces}, there exists an ideal $I$ of $\Lambda$ such that $I_\fm =\Lambda _\fm$ for any $\fm \in A$ and $I_\fn =I(\fn)$ for any $\fn \in B$. In addition, $I$ is the trace of a (finitely generated, by Lemma~\ref{semiperf.fgtrace}) projective right $\Lambda$-module. We claim that $I$ is the minimal element of  $\mathcal{I} (P)$.

Since $I_Q=\Lambda _Q$, $I\cap R\neq \{0\}$. As $R$ has finite character, there exists a finite set $\mathcal M$ of maximal ideals of $R$, such that $R/(I\cap R)\cong \prod _{\fm \in \mathcal M}\left( R/(I\cap R)\right)_\fm$. We may assume that $B\subseteq \mathcal M$. Hence, 
$$P/PI\cong \prod  _{\fm \in \mathcal M} P_\fm /P_\fm I_\fm = \prod  _{\fn\in B} P_\fn /P_\fn I(\fn)$$
which is finitely generated. Therefore, $I\in \mathcal{I} (P)$. 

Let $I'\in \mathcal{I} (P)$ be such that $I'\subseteq I$. Then, because $\Lambda$ is torsion-free as an $R$-module, for any maximal ideal $\fm$ of $R$, $I'_\fm \in \mathcal{I}(P_\fm)$ and $(I')_\fm \subseteq I_\fm=I(\fm)$. The minimality in the choice of $I(\fm)$ in $(a)$ and $(b)$ implies that $(I')_\fm = I_\fm=I(\fm)$ for any maximal ideal $\fm$ of $R$. Hence, $I'=I$. This shows that $I$ is the minimal element of $\mathcal{I} (P)$.

The proof also shows that for any maximal ideal $\fm$ of $R$, $I_\fm$ is the minimal element of $\mathcal{I}(P_\fm)$. 
\end{Proof}

\begin{Lemma} \label{locm}
Let $R$ be a commutative domain with field of fractions $Q$, and let $\Lambda$ be a torsion-free $R$-algebra. Let $I\in\mathcal T(\Lambda)$. Let $P$ be a countably generated projective $\Lambda$-module, and let $V$ be a finitely generated $\Lambda$-module which is torsion-free as an $R$-module. Let $\fp$ be a prime ideal of $R$ such that $V_{\fp} = V_{\fp}I_{\fp}$ and $P_{\fp} = \bigoplus_{i \in \N} U_i$, where each $U_i$ is a finitely generated $\Lambda_\fp$-module, such that infinitely many of these $U_i$ are isomorphic to a $\Lambda_{\fp}$-module $U$ whose trace ideal contains $I_{\fp}$. 

Then, for every finitely generated submodule $P_1$ of $P$, there exists a homomorphism $f \in \Hom_{\Lambda}(P,V)$ such that $f(P_1) = 0$, and $f_{\fp}$ is onto. 
\end{Lemma}
\begin{Proof}
Let $\nu \in \Hom_{\Lambda}(P,\Lambda^{(\omega)})$ and $\pi \in \Hom_\Lambda(\Lambda^{(\omega)},P)$ be such that $\pi \nu = 1_P$.

Note that if $U$ is a projective $\Lambda_{\fp}$-module and $I_{\fp} \subseteq \Tr_\Lambda(U)$, then $U$ generates any $\Lambda_{\fp}$-module $W$ satisfying $WI_{\fp} = W$. If $W$ is also finitely generated, then there exists $k \in \N$ such that $W$ is a homomorphic image of $U^k$. Therefore $V_{\fp}$ is a homomorphic image of $\bigoplus_{i \in C} U_i$ whenever $C$ is a cofinal subset of $\N$. In particular, there exists an onto homomorphism $g \in \Hom_{\Lambda_\fp}(P_\fp,V_\fp)$ such that $g(P_1) = 0$ (recall that $P$ is embedded in $P_\fp$). Since $V$ is finitely generated, there exists $P_0'$, a finitely generated submodule of $P_\fp$, such that $g(P_0') = V_\fp$. Scaling by a suitable element of $R \setminus \fp$ (which is invertible over $\Lambda_\fp$), we observe that there exists a finitely generated $P_0 \subseteq P$ such that $P_0\Lambda_{\fp} = P_0'$.

Let $F$ be a finite subset of $\omega$ such that $\nu(P_0 + P_1)$ is contained in $\Lambda^{(F)}$. Further, let $\iota \in \Hom_{\Lambda}(\Lambda^{(F)}, \Lambda^{(\omega)})$ and $\varrho \in \Hom_{\Lambda}(\Lambda^{(\omega)},\Lambda^{(F)})$ be the canonical embedding and projection. Since $\Lambda^{(F)}$ is finitely generated over $\Lambda$ and torsion-free over $R$, by Lemma~\ref{isofg}, there exist $h \in \Hom_{\Lambda}(\Lambda^{(F)},V)$ and $q \in R_{\fp}$ such that $h \otimes q = g (\pi \iota)_\fp$. Consider $f:= h\varrho\nu$. Then $f \otimes q = h \otimes q \circ (\varrho \nu) \otimes 1 = g (\pi \iota\varrho\nu)_\fp$. Since $\Lambda^{(F)}$ contains $\nu(P_0 + P_1)$, $\pi \iota \varrho\nu$ is the identity on $P_0+P_1$,
and consequently, $f \otimes q$ is onto and $f \otimes q (P_1 \otimes 1) = 0$. This implies that $f_\fp$ is onto and $f(P_1) = 0$.
\end{Proof}

\begin{Remark}
We are going to apply Lemma~\ref{locm} in the following situation: $\Lambda$ is a locally semiperfect torsion-free algebra over an $h$-local domain $R$ such that $\Lambda_Q$ is a simple artinian ring. Assume that $P$ is a countably generated projective $\Lambda$-module and $0 \neq I$ is the least element of  $\mathcal{I} (P)=\{J \unlhd \Lambda \mid P/PJ\text{ is finitely generated}\}$. Take any $\fm \in \mSpec(R)$ and let us check that $P_{\fm}$ is a countable direct sum of finitely generated projective $\Lambda_{\fm}$-modules such that almost all these direct summands have a trace equal to $I_{\fm}$.

Indeed, since $P$ is not finitely generated, $P_{\fm}$ is not a finitely generated $\Lambda_{\fm}$-module by Proposition~\ref{minimallocallysemiperfect}. Since $\Lambda_{\fm}$ is semiperfect, $P_{\fm}$ is a direct sum of finitely generated modules $P_{\fm} = U_1^{m_1} \oplus \cdots \oplus U_{t}^{m_t}\oplus (U_{t+1} \oplus \cdots \oplus U_s)^{(\omega)}$, where $U_1,\dots,U_s$ is a representative set of indecomposable projective $\Lambda_{\fm}$-modules. Let $U' = U_{t+1} \oplus \cdots \oplus  U_{s}$. As explained in the proof of Proposition~\ref{minimallocallysemiperfect}(ii), $I_{\fm} = \Tr_{\Lambda_{\fm}}(U')$. Therefore, $P_{\fm} = U_1^{m_1} \oplus \cdots \oplus U_{t}^{m_t}\oplus U'^{(\omega)}$ is a decomposition of $P_{\fm}$ that has the required properties.    
\end{Remark}

Lemma \ref{locm} can be stated in simpler form when $\fp = (0)$ and $I = \Lambda$.

\begin{Lemma} \label{loc0}
Let $R$ be a commutative domain with field of fractions $Q$, and let $\Lambda$ be a torsion-free $R$-algebra. Let $P$ be a projective $\Lambda$-module such that $P_Q$ is countably generated but not finitely generated free, and let $V$ be a finitely generated $\Lambda$-module which is torsion-free over $R$. Then, for any  finitely generated submodule $P_1$ of $P$, there exists $f \in \Hom_{\Lambda}(P,V)$ such that $f(P_1) = 0$, and $f_Q$ is onto. 
\end{Lemma}

\begin{Prop} \label{Ibig2}
Let $R$ be an $h$-local domain with field of fractions $Q$, and let $\Lambda$ be a locally semiperfect torsion-free $R$-algebra such that $\Lambda_Q$ is simple artinian. Let $P$ be a countably generated projective $\Lambda$-module such that $I = \min \{J \unlhd \Lambda \mid P/PJ\text{ is finitely generated}\}$ exists and is the trace of a countably generated projective module.
Then $P$ is $I$-big.
\end{Prop}
\begin{Proof}
We have to prove that if $U$ is a countably generated projective right $\Lambda$-module and  $\Tr_\Lambda(U) \subseteq I$, then $U$ is a direct summand  of $P$. By Corollary~\ref{omega}, $U^{(\omega)}$ is a direct sum of finitely generated modules. Thus, we may assume that $U$ is a direct sum of finitely generated modules. 

We prove the following claim: If $P_0$ is a finitely generated submodule of $P$ and $V$ is a finitely generated projective $\Lambda$-module with trace contained in $I$, then there exists an epimorphism $f \colon P \to V$ such that $f(P_0) = 0$.

Assume we have proved the claim and let $U = \bigoplus_{i \in \N} U_i$, with each $U_i$ being finitely generated. Let $p_1,p_2,\dots \in P$ be a countable set of generators of $P$. Consider the following construction: Let $f_1\colon P \to U_1$ be an epimorphism, and let $P_1$ be a finitely generated submodule of $P$ such that $f_1(P_1) = U_1$ and $p_1 \in P_1$. Assume that $f_1,f_2,\dots,f_k$ and $P_1,P_2,\dots,P_k$ are defined such that $f_i \colon P \to U_i$ is an epimorphism, $P_1\subseteq \cdots \subseteq P_k$ are finitely generated submodules of $P$ such that $f_{i}(P_{i-1}) = 0$, $p_i \in P_i$, and $f_i(P_i) = U_i$ (we can set $P_0 = 0$). By the claim, there exists an epimorphism $f_{k+1} \colon P \to U_{k+1}$ such that $f_{k+1}(P_k) = 0$. Finally, let $P_{k+1}$ be a finitely generated submodule of $P$ satisfying $P_k \subseteq P_{k+1}$, $p_{k+1} \in P_{k+1}$, and $f_{k+1}(P_{k+1}) = U_{k+1}$. 

Consider the canonical homomorphism $f = \prod_{k \in \N} f_k \colon P \to \prod_{k} U_k$. Observe that $f_{k}(p_l) = 0$ if $k > l$, so ${\rm Im}\ f \subseteq \bigoplus_{i \in \N} U_i$. An easy inductive argument also gives $f(P_i) = \bigoplus_{j = 1}^i U_j$ for any $i \in \N$. Thus $f$ is an epimorphism.

Thus, we are left to prove the claim. First, note that we may assume $I \neq 0$ since any countably generated projective module is $0$-big.
Then $P$ is not finitely generated. Therefore, by Proposition~\ref{minimallocallysemiperfect}, $P_Q$ cannot be finitely generated. Hence $P_Q \simeq \Lambda_Q^{(\omega)}$.

By Lemma~\ref{loc0}, there exists $g \in \Hom_{\Lambda}(P,V)$ such that $g_Q$ is onto and $g(P_0) = 0$. Let $T_0$ be a finitely generated submodule of $P$ containing $P_0$ such that $g_Q(T_0\Lambda_Q) = (g(T_0))_Q = V_Q$. Since $V$ is finitely generated, $V_{\fm} = g(T_0)_{\fm}$ for almost all maximal ideals $\fm$ of $R$. Let $\fm_1,\fm_2,\dots,\fm_k$ be the maximal ideals of $R$ where the equality does not hold. If there are no such ideals, then we are done, since $g$ is onto. So we may assume that $k\ge 1$.

It remains to show that there are homomorphisms $f_1,f_2,\dots,f_k \in \Hom_{\Lambda}(P,V)$ and finitely generated submodules $T_1,T_2,\dots,T_k$ of $P$ such that 
\begin{enumerate}
    \item[(i)] $T_0 \subseteq T_1 \subseteq T_2 \subseteq \cdots \subseteq T_k$
    \item[(ii)] $f_i(T_{i-1}) = 0$ for $i = 1, \dots,k$
    \item[(iii)] $(g+\sum_{j = 1}^{i} f_j)_{\fm_i}((T_{i})_{\fm_i}) = V_{\fm_i}$ for $i = 1,\dots,k$
\end{enumerate}
Then $f = g + \sum_{i = 1}^k f_i \in \Hom_{\Lambda}(P,V)$ satisfies $f(T_0) = g(T_0)$. Hence, $f_{\fm}$ is onto if $\fm$ is not among $\fm_1,\dots,\fm_k$. Moreover, $f(T_i) = (g+\sum_{j = 1}^i f_j)(T_i)$ shows that $f_{\fm_i}$ is onto for every $i = 1,\dots,k$. This implies that $f$ is onto.

Assume $f_1,\dots,f_{\ell-1}$ and $T_0,T_1,\dots,T_{\ell-1}$ are found. Consider $g_{\ell} = g + \sum_{i = 1}^{\ell-1} f_i$. Recall that $V_{\fm_{\ell}}$ is a finitely generated projective module over a semiperfect ring $\Lambda_{\fm_{\ell}}$. Then $V_{\fm_\ell}/V_{\fm_\ell}J(\Lambda_{\fm_\ell})$ is a finitely generated (semisimple module). Therefore, if $\pi\colon V_{\fm_\ell} \to V_{\fm_\ell}/V_{\fm_\ell}J(\Lambda_{\fm_\ell})$ is the canonical projection, $\pi((g_{\ell}(P))_{\fm_{\ell}})$ is a finitely generated submodule of $V_{\fm_{\ell}}/V_{\fm_{\ell}}J(\Lambda_{\fm_{\ell}})$.  Let $T_{\ell}' \subseteq P$ be a finitely generated submodule of $P$ containing $T_{\ell - 1}$ such that $(g_{\ell}(T_{\ell}'))_{\fm_{\ell}} + V_{\fm_{\ell}}J(\Lambda_{\fm_{\ell}}) = (g_{\ell}(P))_{\fm_{\ell}} + V_{\fm_{\ell}}J(\Lambda_{\fm_{\ell}})$.

By Lemma \ref{locm}, there exists $f_{\ell} \in {\rm Hom}_{\Lambda}(P,V)$ such that $(f_{\ell})_{\fm_{\ell}}$ is onto and $f_{\ell}(T_{\ell}') = 0$. Then, the image of $g_{\ell} + f_{\ell}$ contains $g_{\ell}(T_{\ell}')$. Hence, $\im\ (g_{\ell} + f_{\ell})_{\fm_{\ell}} + V_{\fm_{\ell}}J(\Lambda_{\fm_{\ell}}) $ also contains $\im\ (g_{\ell})_{\fm_{\ell}}$, and since $(f_\ell)_{\fm_{\ell}}$ is onto, $\im\ (g_{\ell} + f_{\ell})_{\fm_{\ell}} + V_{\fm_{\ell}}J(\Lambda_{\fm_{\ell}}) = V_{\fm_{\ell}}$. Nakayama's Lemma gives that $(g_{\ell} + f_{\ell})_{\fm_{\ell}}$ is onto. Therefore, there exists a finitely generated submodule $T_{\ell}$ of $P$ containing $T_{\ell}'$ such that $((g_{\ell} + f_{\ell}))_{\fm_{\ell}} ((T_{\ell})_{\fm_{\ell}})= V_{\fm_{\ell}}$, that is, $(iii)$ holds if $i = \ell$.
\end{Proof}

\begin{Th} \label{projectivesls} 
Let $R$ be an $h$-local domain with field of fractions $Q$, and let $\Lambda$ be a locally semiperfect torsion-free $R$-algebra such that $\Lambda_Q$ is a simple artinian ring. Then $V_r^*(\Lambda)=V_r(\Lambda)\sqcup B_r(\Lambda)\cong V_\ell(\Lambda)\sqcup B_\ell(\Lambda)=V_\ell^*(\Lambda)$.
\end{Th}
\begin{Proof} 
Propositions~\ref{minimallocallysemiperfect} and~\ref{Ibig2} show that $V_r^*(\Lambda)=V_r(\Lambda)\sqcup B_r(\Lambda)$. The symmetry of the hypothesis implies that $V_\ell(\Lambda)\sqcup B_\ell(\Lambda)=V_\ell^*(\Lambda)$. By Remark~\ref{descriptiontracesalmostsemiperfect}, we can deduce from Corollary~\ref{isomonoidtraces} that both monoids are isomorphic.
\end{Proof}

Hinohara~\cite{hinohara} proved that every projective module over an $h$-local (or, more generally, weakly noetherian) domain is either finitely generated or free. We can extend this result to our setting. 

\begin{Cor} \cite{hinohara} \label{hinohara}
Let $R$ be an $h$-local domain with field of fractions $Q$, and let $\Lambda$ be a locally semiperfect torsion-free $R$-algebra such that $\Lambda_Q$ is simple artinian. For each maximal ideal $\fm$ of $R$, assume that $\Lambda_{\fm}$ is Morita equivalent to a local ring. Then any infinitely generated projective $\Lambda$-module is free.
\end{Cor}
\begin{Proof}
By Remark~\ref{descriptiontracesalmostsemiperfect}, ${\mathcal T}(\Lambda) = \{0,\Lambda\}$, and hence $V(\Lambda) \sqcup B(\Lambda) = V(\Lambda) \sqcup \{[\Lambda^{(\omega)}]\}$.
The result then follows from Theorem~\ref{projectivesls}.
\end{Proof}

\section{Rank and Genus of countably generated projective modules}\label{sec:7}

We continue with the same setting as in \S \ref{sec:6}. Let $R$ be an $h$-local domain with field of fractions $Q$. Let $\Lambda$ be a  locally semiperfect torsion-free algebra over  $R$ such that $\Lambda_Q$ is simple artinian. That is  $\Lambda _Q\cong M_k (D)$, for $D$ a division ring.

There are two natural monoid morphisms. The first one is the one associated with tensoring with $Q$:
    \[\begin{tikzcd}
        \Psi \colon V^*(\Lambda) \rar & V^*(\Lambda _Q)
    \end{tikzcd}\]
defined by $\Psi ([P])=[P_Q] =[S^r]$, where $S$ is a simple right $\Lambda _Q$-module. Then we have a monoid morphism that, by abuse of notation, we also call $\Psi$
    \[\begin{tikzcd}
        \Psi \colon V^*(\Lambda) \rar & \No^*
    \end{tikzcd}\]
defined by $\Psi ([P])=r$ if $P$ is finitely generated and $P_Q\cong S^r$, and $\Psi ([P])=\infty$ if $P$ is not finitely generated. We call $\Psi ([P])$ the \emph{rank} of $P$, and we refer to $\Psi$ as the \emph{rank map}.

Associated with the localization at the different maximal ideals, we have the monoid morphism
    \[\begin{tikzcd}
        \Phi \colon V^*(\Lambda) \rar & \prod_{\fm \in \mSpec (R)} V^*(\Lambda _\fm) 
    \end{tikzcd}\]
defined by $\Phi ([P])=([P_\fm])_{\fm \in \mSpec (R)}$. We shall refer to $\Phi$ as the \emph{genus morphism}, and to $\Phi ([P])$ as the genus of $[P]$.

We say that two countably generated projective right $\Lambda$-modules $P$ and $P'$ \emph{are in the same genus} if   $P_\fm\cong P'_\fm$ for every maximal ideal $\fm$ of $R$, that is, if and only if $\Phi ([P])=\Phi ([P'])$. 

In this section, we want to describe the images of $\Psi$ and $\Phi$.

At this point, we have plenty of information about $V^*(\Lambda)$. By Theorem~\ref{projectivesls}, $V^{*}(\Lambda) = V(\Lambda) \sqcup B(\Lambda)$, and by Lemma~\ref{semiperf.fgtrace}, any ideal of $\mathcal T(\Lambda)$ is a  trace ideal of a finitely generated projective module, and by Corollary~\ref{minimal}, $\mathcal T(\Lambda)$ is finite.

In general, there may be plenty of finitely generated projective right modules over $\Lambda$ in the same genus that are not isomorphic. For example, let $R=\Lambda$ be a Dedekind domain. Then any localization of a projective module at a maximal ideal is free. So two projective modules are in the same genus if and only if they have the same rank. Certainly, two finitely generated projective modules over a Dedekind domain with the same rank are not, in general, isomorphic. 

The situation changes for projective $\Lambda$-modules whose isomorphism class is in $B(\Lambda)$, or for finitely generated projective $\Lambda$-modules when $\Lambda$ is a suitable  algebra over a semilocal domain $R$.

\begin{Prop} \label{genusisomorphic}
Let $R$ be an $h$-local domain with field of fractions $Q$. Let $\Lambda$ be a locally semiperfect torsion-free $R$-algebra such that $\Lambda _Q$ is simple artinian. Then $\Phi (B(\Lambda))\subseteq \prod _{\fm \in \mSpec (R)}  B(\Lambda _\fm)$. Moreover, if $[P]$ and $[P']$ are two elements in $B(\Lambda)$, then $\Phi ([P])=\Phi ([P'])$ if and only if $[P]=[P']$. 
\end{Prop}
\begin{Proof} 
The first part of the proof follows from Corollary~\ref{relbig_loc} and Proposition~\ref{minimallocallysemiperfect}.

For the second part of the proof, we are going to use the notation of  Corollary~\ref{iso.monoid} that is, an element $[P]\in B(\Lambda)$ can be thought of as the isomorphism class of a pair $(P/PI, I)$ where  $I\in \mathcal{T} (\Lambda)\setminus 
\{0\}$, and  $P$ is relatively $I$-big. Hence, $\Phi ([(P/PI, I)])= ([(P_\fm/P_\fm I_\fm, I_\fm)])_{\fm \in \mSpec(R)}.$

By Lemma~\ref{almosttrace}, there exist $\fm _1,\dots ,\fm _\ell$ maximal ideals of $R$, such that $\Lambda /I\cong \prod _{i=1}^\ell \Lambda _{\fm _i}/I_{\fm _i}$. Hence, $P/PI\cong \prod _{i=1}^\ell P _{\fm _i}/P _{\fm _i}I_{\fm _i}$, and $P _{\fn}=P _{\fn}I_{\fn}$ is free and $I_{\fn}=\Lambda _\fn$, for  any $\fn \in \mSpec (R)\setminus \{\fm _1,\dots ,\fm _\ell\}$.  

By Theorem~\ref{dealtraces}, $I$ is uniquely determined by $I_\fm$, $\fm \in \mSpec (R)$.
Therefore, by Lemma~\ref{Jbigdetermined}, it follows that $\Phi$ is injective when restricted to $B(\Lambda)$.
\end{Proof}

\begin{Prop} \label{genusisomorphic_semilocal}
Let $\Lambda$ be a locally semiperfect, torsion-free algebra over a semilocal domain $R$. Then:
\begin{itemize}
    \item[(i)] $\Lambda$ is semilocal;
    \item[(ii)] Assume $\Lambda$ is module-finite over $R$, and let $P_1$ and $P_2$ be two finitely generated projective right $\Lambda$-modules.  Then $P_1\cong P_2$ if and only if $\Phi ([P_1])=\Phi ([P_2])$. 
\end{itemize}
\end{Prop}

\begin{Proof}
$(i)$  Let $\mSpec (R)=\{\fm_1,\dots, \fm _\ell\}$. Consider the embedding $$\alpha \colon \Lambda \to \prod _{i \in \{1,\dots ,\ell\}}\Lambda _{\fm _i}$$ induced by taking the   localization map at each component. If $\lambda \in \Lambda$ is such that $\alpha (\lambda)$ is invertible, then for each $i\in \{1,\dots ,\ell\}$, there exists $s_i\in R\setminus \fm _i$ and $\mu _i \in \Lambda$ such that $\lambda \mu _i=s_i$. Then $I=\sum _{i=1}^\ell s_iR$ is not contained in any maximal ideal of $R$, so $I=R$. Since $R=I\subseteq \lambda \Lambda$, we deduce that $\lambda $ is a right invertible element of $\Lambda$. By symmetry, $\lambda$ is an invertible element of $\Lambda$. This shows that $\alpha$ is a local ring homomorphism (the image of an element  is invertible if and only if it is an invertible element of $\Lambda$). Since $\prod _{i \in \{1,\dots ,\ell\}}\Lambda _{\fm _i}$ is a finite product of semiperfect rings, we deduce that $\Lambda$ is semilocal (cf. \cite[Theorem~4.2]{libro}).

$(ii)$ This statement is proved in \cite[Corollary~4.5]{AHP}.
\end{Proof}

The following straightened version of the Package Deal Theorem to our setting will be the key ingredient to describe the image of the genus map $\Phi$.

\begin{Prop} \label{fgprojlocsem} 
Let $R$ be an $h$-local domain with field of fractions $Q$. Let $\Lambda$ be a locally semiperfect torsion-free $R$-algebra such that $\Lambda_Q$ is simple artinian. Let $S$ be a simple right $\Lambda _Q$-module.

Fix $r\ge 1$. For each maximal ideal $\fm$ of $R$, let $P(\fm)$ be a finitely generated projective right $\Lambda _\fm$-module such that $P(\fm)_Q\cong S^r$. Then there exists a finitely generated projective right $\Lambda$-module $P$ such that $P_\fm \cong P(\fm)$ for every maximal ideal $\fm$ of $R$.

Notice that, by definition, $\Psi ([P])=r$.
\end{Prop}
\begin{Proof}
Let $k$ be such that $\Lambda _Q\cong S^k$. By Lemma~\ref{finitely_many}, there is a finite set of maximal ideals of $R$, $\fm_1,\dots,\fm_\ell$, such that $\Lambda_{\fm}\cong M_{k}(T_{\fm})$ for some local ring $T_{\fm}$
whenever $\fm $ is a maximal ideal not in $\mathcal{S}=\{\fm_1,\dots,\fm_\ell\}$. Note that if $\fm \not\in \mathcal{S}$, then $\Lambda _\fm$ is isomorphic to a ring of $k\times k$ matrices over a local ring, so it has a unique indecomposable projective module $X(\fm)$ such that $X(\fm)_Q\cong S$ and any projective $\Lambda (\fm)$-module is isomorphic to a suitable number of copies of $X(\fm)$. Therefore, for any $\fm \not \in \mathcal{S}$, $P(\fm)\cong X(\fm)^r$.

Let $q$ be a multiple of $k$ that is greater than $r$, and let $F = \Lambda^{q/k}$. Further, let $U$ be a finitely generated submodule of $F$ such that $U_Q \cong S^r$.
Since the embedding $U \hookrightarrow F$ splits over $\Lambda _Q$, it splits over almost all maximal ideals of $R$  \cite[Lemma~4.2]{AHP}.
Then $U$ is a finitely generated $\Lambda$-module which is torsion-free as an $R$-module, such that $U_{\fm}$ is projective for almost all maximal ideals $\fm$ of $R$ and $U_Q \cong S^r$. So let $B$ be a finite set of maximal ideals of $R$ such that $\mathcal{S}\subseteq B$, and satisfying that $U_{\fm}$ is projective for any $\fm \not \in B$. In particular, for any $\fm \not \in B$, $U_{\fm}\cong P(\fm)\cong X(\fm)^r$.

By Lemma~\ref{d}, for any $\fm \in B$, we can choose   $U(\fm)$   a submodule of $U_{\fm}$ such that 
$U(\fm) \cong P(\fm )$. For $\fm \not \in B$, we set $U(\fm)=U_\fm$. 

Note that for any maximal ideal $\fm$ of $R$, $U(\fm)_Q = U_Q$ since all these modules have length $r$ over $\Lambda _Q$.
Therefore, we can apply  Theorem~\ref{dealsubmodules} and conclude that  there exists a finitely generated 
$\Lambda$-module $P$ of $U$ such that $P(\fm)\cong U(\fm) = P_{\fm}$ for every $\fm \in \mSpec(R)$.
By Corollary \ref{loc.proj}, $P$ is a projective module over $\Lambda$.
\end{Proof}

We introduce the notation necessary to more thoroughly detail the image of the genus map $\Phi$ and the rank map $\Psi$.

\begin{notation} \label{not}
Let $R$ be an $h$-local domain with field of fractions $Q$. Let $\Lambda$ be a locally semiperfect torsion-free $R$-algebra such that $\Lambda_Q$ is simple artinian.
Let $S$ be a simple $\Lambda_Q$-module, and let $k$ be the length of $\Lambda_Q$, i.e., $\Lambda_Q \cong S^k$ or, equivalently, $\Psi (\Lambda)=k$.
By Lemma~\ref{finitely_many}, there is a finite set of maximal ideals $\fm_1,\dots,\fm_\ell$ such that $\Lambda_{\fm}\cong M_{k}(T_{\fm})$ for some local ring $T_{\fm}$
whenever $\fm $ is a maximal ideal not in $\mathcal{S}=\{\fm_1,\dots,\fm_\ell\}$.

If $\ell = 0$, then the assumptions of Corollary~\ref{hinohara} are satisfied. So in this case, the image of the rank map $\Psi$ is $\No ^*$ and the genus map composed with the natural isomorphism   $V^*(\Lambda _\fm) \cong \No ^*$
    \[\begin{tikzcd}
        \Phi \colon V^*(\Lambda) \rar & \prod_{\fm \in \mSpec (R)} V^*(\Lambda _\fm) \cong \prod_{\fm \in \mSpec (R)} \No ^*
    \end{tikzcd}\]
just sends $[P]$ to the constant sequence with value the rank of $P$.

So from now on we will assume that $\ell \ge 1$.

For $i = 1,\dots,\ell$, let $U_{i,1},\dots,U_{i,t_i}$ be representatives of indecomposable projective right $\Lambda_{\fm_i}$-modules, and let $r_{i,1},\dots,r_{i,t_i}$ be the lengths of their 
localizations over $\Lambda_Q$, i.e., $(U_{i,j})_{Q} \cong S^{r_{i,j}}$.

Note that, for any $i=1,\dots ,\ell$, if  the right $\Lambda _{\fm_i}$-module $\Lambda_{\fm_i}$ is isomorphic to $U_{i,1}^{m_{i,1}}\oplus \dots \oplus U_{i,t_i}^{m_{i,t_i}}$, then we  have the relation
\[\sum _{j=1}^{t_i}r_{i,j}m_{i,j}=(r_{i,1},\dots ,r_{i,t_i})\begin{pmatrix}
    m_{i,1}\\ \vdots \\ m_{i,t_i}
\end{pmatrix}=k.\]
\end{notation}

We proceed to construct examples of finitely generated locally semiperfect torsion-free $\Z$-algebras. In particular, they are algebras over an $h$-local domain with an infinite number of maximal ideals.  These examples show that there are no other restrictions on the parameters $k, \ell, t_i, r_{ij}, m_{ij}$ than the $\ell$ equations above.

\begin{Exs} \label{exs:locallysemiperfect}
(1) Let $k$ be a positive integer and let  $p$ be a prime. Assume we are given positive integers $r_1,\dots,r_t$, $m_1,\dots,m_t$ such that $\sum_{i = 1}^t m_ir_i = k$. Our aim is to construct a semiperfect $\mathbb{Z}_{(p)}$-order $\Lambda$ in ${\rm M}_k(\mathbb{Q})$ such that $\Lambda \cong \bigoplus_{i = 1}^t U_i^{m_i}$ as right modules, where $U_1,\dots,U_t$ is a representative set of indecomposable projective $\Lambda$-modules and $r_i$ is the length of $(U_i)_{Q}$.

Write $\{1,\dots,k\}$ as a disjoint union of sets $\bigcup_{i = 1}^t (A_{i,1} \cup A_{i,2} \cup \cdots \cup A_{i,m_i})$, where $|A_{i,j}| = r_i$. For any $i$ and any $1 \leq a,b \leq m_i$
we fix a bijection $f_{i,a,b} \colon A_{i,a} \to A_{i,b}$ such that $f_{i,a,a}$ is 
the identity of $A_{i,a}$ and $f_{i,b,c}f_{i,a,b} = f_{i,a,c}$ if $1 \leq a,b,c \leq m_i$.

For elements $i,a,b$, $1 \leq i \leq t$, $1 \leq a ,b \leq m_i$ define $e_{i,a,b} \in {\rm M}_k(\mathbb{Q})$ by $e_{i,a,b} = \sum_{j \in A_{i,a}} E_{f_{i,a,b}(j),j}$, where $E_{x,y}$
is a matrix unit in ${\rm M}_k(\mathbb{Q})$. 
Note that $e_{i,a,a}^2 = e_{i,a,a}$ and $e_{i,b,c}e_{i,a,b} = e_{i,a,c}$.
Further note that $e_{i,a,b}e_{j,c,d} = 0$ if $1 \leq i \neq j \leq k$, 
$1 \leq a,b \leq m_i$ and $1 \leq c,d  \leq m_j$ or if $i = j$ and $1 \leq d \neq a \leq m_i$.

Consider $J = {\rm M}_{k}(p\mathbb{Z}_{(p)})$ and $\Lambda = J+ \sum_{i = 1}^t \sum_{1 \leq a,b \leq m_i}
e_{i,a,b} \mathbb{Z}_{(p)}$. The following facts are easily verified:
\begin{enumerate}
    \item[(i)] $\Lambda$ is a subring of ${\rm M}_{k}(\mathbb{Q})$, is a finitely generated $\Z _{(p)}$-algebra and $\Lambda_Q = {\rm M}_{k}(\mathbb{Q})$.
    \item[(ii)] $J$ is an ideal of $\Lambda$, $J \subseteq J(\Lambda)$.
    \item[(iii)] For any $i,j$ $(e_{i,j,j}\Lambda)/(e_{i,j,j}J)$ is a simple $\Lambda$-module. 
    \item[(iv)] $e_{i,a,a} \Lambda \cong e_{j,b,b}\Lambda$ if and only if $i = j$
    \item[(v)] $e_{i,a,a} \Lambda e_{j,b,b} = e_{i,a,a} J e_{j,b,b}$ whenever $1 \leq i \neq j \leq t$, $1 \leq a \leq m_i,1 \leq b \leq m_j$
\end{enumerate}

From $(iii)$ it follows that $\Lambda/J$ is a direct sum of simple modules. Hence $\Lambda/J$
is semisimple artinian and therefore $J = J(\Lambda)$. Also note that the canonical 
projection $e_{i,j,j}\Lambda \to e_{i,j,j}\Lambda/e_{i,j,j}J$ is a projective cover of 
$e_{i,j,j}\Lambda/e_{i,j,j}J$. So $\Lambda$ is semiperfect and $e_{i,j,j}\Lambda$ are indecomposable
projectives. Also $e_{i,a,a}\Lambda \cong e_{j,b,b} \Lambda$ if and only if $i = j$. So we have 
that $\Lambda = \bigoplus_{i = 1}^t \bigoplus_{a = 1}^{m_i} e_{i,a,a} \Lambda \cong \bigoplus_{i = 1}^t (e_{i,1,1}\Lambda)^{m_i}$. Also, it is easy to check that $e_{i,a,a}{\rm M}_{k}(\mathbb{Q})$
has length $r_{i}$ since $|A_{i,a}| = r_{i}$.

(2) Let $k,\ell$ be positive integers $p_1,\dots,p_\ell$ pairwise different primes.
Assume that for each $i$ positive integers $t_i$, $r_{i,j},m_{i,j}$, $j \in \{1,\dots,t_i\}$
are given such that $\sum_{j = 1}^{t_i} r_{i,j}m_{i,j}=k$ for every $i = 1,\dots,\ell$.
By part (a) there exists a semiperfect $\mathbb{Z}_{(p_i)}$-order $\Lambda(p_i)$ in ${\rm M}_k(\mathbb{Q})$
such that $\Lambda(p_i) \cong \bigoplus_{j = 1}^{t_i} U_{i,j}^{m_{i,j}}$ where $U_{i,1},\dots U_{i,t_i}$ are pairwise 
non-isomorphic projective $\Lambda(p_i)$-modules and $(U_{i,j})_Q$ has rank $r_{i,j}$.

Now the Package Deal Theorem (here it is enough to use the noetherian version from \cite{levyodenthal2}) finds a finitely generated $\mathbb{Z}$-submodule $\Lambda \subseteq {\rm M}_{k}(\mathbb{Z})$ such that $\Lambda_{(p_i)} = \Lambda(p_i)$ for 
$i = 1, \dots, \ell$ and $\Lambda_{(q)} = {\rm M}_k(\mathbb{Z}_{(q)})$ if $q$ is a prime different from $p_1,\dots,p_{\ell}$. It is easy to see 
that $\Lambda$ has to be a subring of ${\rm M}_k(\mathbb{Z})$  and that $\Lambda_{(p_i)} \simeq \Lambda(p_i)$ as rings. 
\end{Exs}

\begin{Rem} 
We keep the  Notation~\ref{not}. First, we observe that 
    \[\begin{tikzcd}
        \Phi \colon V^*(\Lambda) \rar & \prod_{\fm \in \mSpec (R)} V^*(\Lambda _\fm)=\prod _{i=1}^\ell  V^*(\Lambda _{\fm _i})\times \prod _{\fm \not \in \mathcal {S}} V^*(\Lambda _\fm)
    \end{tikzcd}\]
So, for any $\fm \not \in  \mathcal {S}$, $V^*(\Lambda _\fm)\cong \No ^*$, and for any $[P]\in V^*(\Lambda)$,  this isomorphism sends $[P_\fm]$ to its rank, which is $\Psi ([P])$. As we are assuming that $\ell \ge 1$,  the genus map is completely determined by its image in  $\prod _{i=1}^\ell  V^*(\Lambda _{\fm _i})$.

Then, there is a monoid morphism
    $$\begin{tikzcd}[row sep=0ex]
        \Phi'\colon V(\Lambda)\rar & \prod_{i=1}^\ell V(\Lambda _{\fm _i})\cong \prod _{i=1}^\ell \left (\No ^{t_i}\right)
    \end{tikzcd}$$
    $$\Phi ' ([P])=([P_{\fm _1}], \dots , [P_{\fm _\ell}])=\Big(\textstyle\sum _{j=1}^{t_1} [U_{1,j}]^{x_{1,j}}, \dots ,  \sum _{j=1}^{t_\ell} [U_{\ell,j}]^{x_{\ell,j}}\Big).$$
Where the  isomorphism $\prod_{i=1}^\ell V(\Lambda _{\fm _i})\cong \prod _{i=1}^\ell \left (\No ^{t_i}\right)$ assigns to $\Phi ' ([P])$ the element 
\[((x_{1,1},\dots ,x_{1,t_1}),\dots , (x_{\ell,1},\dots ,x_{\ell,t_\ell }))\in \prod _{i=1}^\ell \left (\No ^{t_i}\right).\]
If $P_Q\cong S^r$, then for any $i=1,\dots ,\ell$, we have the identity
    \begin{equation}\label{eq:systemfg}
        \sum _{j=1}^{t_i}r_{i,j}x_{i,j}=(r_{i,1},\dots ,r_{i,t_i})\begin{pmatrix}
            x_{i,1}\\ \vdots \\ x_{i,t_i}
        \end{pmatrix}=r
    \end{equation}
Conversely,  by Proposition~\ref{fgprojlocsem}, an element $\mathbf{X}\in \prod _{i=1}^\ell \left (\No ^{t_i}\right)$ corresponds to a finitely generated projective right $\Lambda$-module if and only if there exists $r$ such that the components of $\mathbf{X}$ satisfy the system of equations \eqref{eq:systemfg}. This describes all the possible genus of finitely generated projective $\Lambda$-modules. 
 \end{Rem}
 
In summary, we have proved the following.

\begin{Prop} \label{systemfg}
Let $\mathcal{A}=\Phi '(V(\Lambda))$ which is a submonoid of $\prod _{i=1}^\ell \left (\No ^{t_i}\right)$. Then $\mathbf{X}\in 
  \prod _{i=1}^\ell \left (\No ^{t_i}\right)$  is an element of $\mathcal{A}$ if and only if there exists $r\in \No$ such that the components of $\mathbf{X}$ satisfy
 \[\sum _{j=1}^{t_i}r_{i,j}x_{i,j}=(r_{i,1},\dots ,r_{i,t_i})\begin{pmatrix}
    x_{i,1}\\ \vdots \\ x_{i,t_i}
\end{pmatrix}=r\]
for any $i=1,\dots ,\ell$.  

If $P$ is a finitely generated projective right $\Lambda$-module such that   $\Phi '([P])=\mathbf{X}$, then $P_Q\cong S^r$. So 
    $$\Psi (V(\Lambda))=\{r\in \No \mid \mbox{there exists $\mathbf{X}\in 
  \prod _{i=1}^\ell \left (\No ^{t_i}\right)$ satisfying  \eqref{eq:systemfg}}\}.$$
\end{Prop}

\begin{Remark} \label{lisone} 
If $\ell =1$ and $U_1, \dots , U_t$ are the indecomposable projective modules of $\Lambda _\fm$, up to isomorphism, and $\Psi (U_i)=r_i$ for $i=1,\dots , t$. Then $\Psi (V(\Lambda))$ is the submonoid of $\No$ generated by $r_1, \dots ,r_\ell$. By Proposition~\ref{fgprojlocsem}, $\Phi ' (V(\Lambda))\cong \No ^t$, because for any $i=1,\dots ,\ell$, there exists an (indecomposable) finitely generated projective right $\Lambda$-module $P_i$ such that $(P_i)_{\fm _i}\cong U_i$.

Clearly, $\Psi (V^*(\Lambda))= \Psi (V(\Lambda))\cup \{\infty\}$  and $\Phi ' (V^*(\Lambda))\cong (\No ^*) ^t$. Notice that, in this case, by Proposition~\ref{genusisomorphic}, any projective $\Lambda$-module is isomorphic to a direct sum of finitely generated projective $\Lambda$-modules. 
\end{Remark}

Now we proceed to extend this description of the genus to countably generated projective modules. In view of Remark~\ref{lisone}, we may assume $\ell\ge 2$.

\begin{Prop} \label{systemcg}
Extend the monoid morphism $\Phi '$  to $V^*(\Lambda)$ and enlarge $\No$ to $\No^*= \No\cup \{\infty\}$. Then $\mathbf{Y}\in \prod _{i=1}^\ell \left ((\No^*) ^{t_i}\right)$ is an element in $\Phi ' (B(\Lambda))=\Phi '(V^*(\Lambda)\setminus V(\Lambda))$ if and only if its components are solutions of 
    \begin{equation}\label{eq:systemcg}
        \sum _{j=1}^{t_i}r_{i,j}y_{i,j}=(r_{i,1},\dots ,r_{i,t_i})\begin{pmatrix}
            y_{i,1}\\ \vdots \\ y_{i,t_i}
        \end{pmatrix}=\infty
    \end{equation}
for any $i=1,\dots ,\ell$. 
\end{Prop}
\begin{Proof} 
By Theorem~\ref{projectivesls}, $V^{*}(\Lambda) = V(\Lambda) \sqcup B(\Lambda)$, so $\Phi ' (B(\Lambda))=\Phi '(V^*(\Lambda)\setminus V(\Lambda))$. 

Let  $\mathbf{Y}\in \prod _{i=1}^\ell \left ((\No^*) ^{t_i}\right)$. For $i=1,\dots ,\ell$, let $S_i (\mathbf{Y}) =\{j\mid y_{i,j}=\infty\}$. Observe that $\mathbf{Y}$ is a solution of the equations \eqref{eq:systemcg} if and only if $S_i(\mathbf{Y})\neq \emptyset$ for any $i=1,\dots ,\ell$.

Let $P$ be a countably generated right $\Lambda$-module such that it is not finitely generated. Then $\Psi (P)=\infty$, which implies that, for any $i=1,\dots ,\ell$, there exists $j$ such that the indecomposable module $U_{i,j}$ appears infinitely many times in the decomposition of $P_\fm$ as a direct sum of indecomposable modules. Therefore, if $\mathbf{Y}=\Phi '([P])$, then $S_i(\mathbf{Y})\neq \emptyset$ for any $i=1,\dots ,\ell$. So $\mathbf{Y}$ is a solution to the system \eqref{eq:systemcg}. 

To prove the converse, let $\mathbf{Y}\in \prod _{i=1}^\ell \left ((\No^*) ^{t_i}\right)$ such that $S_i(\mathbf{Y})\neq \emptyset$ for any $i=1,\dots ,\ell$. Then we show that there is a countably generated projective module $P$ such that $\mathbf{Y}=\Phi '([P])$.

For $i=1,\dots ,\ell$, let $P'(\fm _i)=\bigoplus _{j\in S_i(\mathbf{Y})}U_{i,j}$ and let $r_i$ be the rank of $P'(\fm _i)$. Let $r$ be the least common multiple of $r_1,\dots ,r_\ell ,k$. For $i=1,\dots ,\ell$, set $P(\fm _i)=P'(\fm _i)^{\frac m{r_i}}$ and set $P(\fm)=\Lambda _\fm ^{\frac rk}$ if $\fm$ is a maximal ideal of $R$ different from $\fm_1,\dots,\fm_\ell$. By Package deal Theorem~\ref{dealprojective}, there exists a finitely generated projective $\Lambda$-module $P'$ such that $P'_\fm \cong P(\fm)$ for any maximal ideal $\fm$ of $R$. 

Let $I$ be the trace ideal of $P'$. Then $I_{\fm _i} =\sum_{j\in S_i(\mathbf{Y})} \mathrm{Tr} (U_{i,j})$ for $i=1,\dots ,\ell$, and $I_{\fm}=\Lambda _\fm$ for the other maximal ideals of $R$. As $R$ is $h$-local,  
$$\Lambda /I \cong \prod _{i=1}^\ell \Lambda _{\fm _i}/I_{\fm _i}.$$
For any $i=1,\dots ,\ell$, set $A_i =\{1,\dots ,t_i\}\setminus  S_i(\mathbf{Y})$. Then the indecomposable projective modules over $\Lambda _{\fm _i}/I_{\fm _i}$ are $\{U_{i,j}/U_{i,j}I_{\fm _i}\} _{j\in A_i}$.

Now consider the finitely generated projective $\Lambda /I$-module 
    $$\overline {P} =\bigoplus_{i=1}^\ell  \Big( \bigoplus _{j\in A_i} (U_{i,j}/U_{i,j}I_{\fm _i})^{y_{i,j}}\Big) .$$
By Lemma~\ref{traces.quo}, there exists a countably generated projective $\Lambda$-module $P$ such that $P/PI\cong \overline{P}$. Then $P\oplus (P')^{(\omega)}$ is an $I$-big countably generated projective $\Lambda$-module such that $\Psi' (P\oplus (P')^{(\omega)})= \mathbf{Y}.$
\end{Proof}

\begin{Remark} In Proposition~\ref{systemcg}, the condition on $\mathbf{Y}\in \prod _{i=1}^\ell \left ((\No^*) ^{t_i}\right)$ is equivalent to saying that, for each $i\in \{1,\dots ,\ell\}$, at least one of the components of $(y_{i,1}, \dots, y_{i,t_i})$ is $\infty$.

We have kept the notation of a system of linear equations because it is consistent with the characterization of the monoids that can be realized as $V (\Lambda)\sqcup B (\Lambda)$ for $\Lambda$ a semilocal ring given in \cite{crelle}. Namely, the description of $V (\Lambda)\sqcup B (\Lambda)$ for $\Lambda$ a semilocal ring can always be given in terms of the solutions in $\No ^*$ of a system of linear equations and linear congruences with coefficients in $\No$. However, it is not at all clear which invariants of the semilocal ring  determine the coefficients of the system defining the monoid of relatively big projectives. 

In our situation, by Proposition~\ref{genusisomorphic_semilocal}, $\Lambda _\Sigma$ where $\Sigma =R\setminus \bigcup _{i=1}^\ell \fm _i$ is a locally semiperfect, semilocal algebra over the semilocal domain $R_\Sigma$.  We are proving that, at least for the genus, the system of equations that determines the monoid can be deduced from the structural properties of $\Lambda$ and its localizations at the maximal ideals of $R$. Note that if $\Lambda$ is finitely generated over $R$, then, by Proposition~\ref{genusisomorphic_semilocal}, $V^*(\Lambda _\Sigma)=V(\Lambda)\sqcup B(\Lambda)\cong \Phi ' (V^*(\Lambda))$.
\end{Remark}

\section{Infinite sums of finitely generated projective modules} \label{sec:8}

In this section, we continue working with an algebra $\Lambda$ over an $h$-local domain $R$, with field of fractions $Q$, and such that $\Lambda _Q$ is simple artinian. 

We keep using the notation from the previous section. In particular, the one introduced in  Notation~\ref{not}. As in Proposition~\ref{systemfg}, let $\mathcal{A}=\Phi ' (V(\Lambda))$. If $\ell$ is either $0$ or $1$, then every projective module is a direct sum of finitely generated ones (cf. Remark~\ref{lisone}). So we assume $\ell \ge 2$.

Let $\mathcal{B}$ be the image via $\Phi '$ of $B(\Lambda)=V^*(\Lambda)\setminus V(\Lambda)$ in $\prod _{i=1}^\ell \left ((\No^*) ^{t_i}\right)$.

By Proposition~\ref{genusisomorphic}, the genus classifies the modules of $B(\Lambda)$ up to isomorphism. Therefore, the modules in  $B(\Lambda)$ that are  isomorphic to a direct sum of finitely generated modules are those such that the image of their isomorphism class via $\Phi '$ lies in 
$$ \mathcal{A}+\infty \cdot (\mathcal{A}\setminus \{0\})=\{a+\infty\cdot a'\mid a,a'\in \mathcal{A} \mbox{ and } a'\neq 0\}.$$

In this section we will examine the question of whether any projective $\Lambda$-module is a direct sum of finitely generated modules. By Proposition~\ref{monoidbtrivial} and Lemma~\ref{semiperf.fgtrace},  we only need to characterize when finitely generated projective modules modulo an ideal of $\mathcal T(\Lambda)$ lift to finitely generated projective modules. Our next Lemma shows that we only need to ensure this property for the minimal elements in $\mathcal T(\Lambda)\setminus \{0\}$.

\begin{Lemma} \label{criteria}
Let $R$ be an $h$-local domain with field of fractions $Q$, and let $\Lambda$ be a locally semiperfect torsion-free $R$-algebra such that $\Lambda_Q$ is simple artinian. The following are equivalent
\begin{enumerate}
  \item[(i)] For any $I \in \mathcal T(\Lambda)$ and any finitely generated projective $\Lambda/I$-module $P'$, there exists a finitely generated projective $\Lambda$-module $P$
  such that $P' \cong P/PI$.
  \item[(ii)] For any $I$ which is a minimal element of $\mathcal T(\Lambda) \setminus \{ 0 \}$ and any finitely generated projective $\Lambda/I$-module $P'$, there exists a finitely generated projective $\Lambda$-module $P$ such that $P' \cong P/PI$.
\end{enumerate}
\end{Lemma}
\begin{Proof}
Of course, only $(ii) \Rightarrow (i)$ needs a proof. To prove $(ii) \Rightarrow (i)$, let $I$ be a non-zero trace of $\Lambda$ containing a minimal trace $X \in T(\Lambda)\setminus \{0\}$.
Such an $X$ exists by Corollary~\ref{minimal}.

By Lemma~\ref{almosttrace}, there exist $\fm_1,\dots ,\fm _\ell$ maximal ideals of $R$ such that $\Lambda /X=\prod _{i=1}^\ell \Lambda _{\fm _i}/X_{\fm _i}$. Thus, it is a finite product of semiperfect rings.  Therefore, $\Lambda/X$ is semiperfect.

Let $P'$ be a finitely generated projective $\Lambda/I$-module. By Lemma~\ref{lifting}, there exists a finitely generated projective $\Lambda/X$-module 
$P_1$ such that $P_1/P_1I \cong P'$. Now we can use $(ii)$ to find a finitely generated projective $\Lambda$-module $P$ such that $P/PX = P_1$. Then $P/PI \cong (P/PX)/(PI/PX) \cong (P_1/P_1I) \cong P'$.
\end{Proof}

The next lemma will be useful for checking when the criteria of Lemma~\ref{criteria} is satisfied.

\begin{Lemma} \label{equations}
Let $2\le \ell \in \N$. Let $r_1,\dots ,r_{\ell }, s \in \N$. Let $m=\mathrm{lcm} (r_1,\dots ,r_{\ell-1})>0$. Then the system 
    \begin{equation}\label{eq:equations1}
        \begin{pmatrix}
            r_1&0&\dots&0 &-s\\0&r_2&\dots&0 &-s\\\vdots & &\ddots&&\vdots\\
            0&0& \dots &r_{\ell -1}&-s
        \end{pmatrix} \begin{pmatrix}
            x_1\\x_2\\ \vdots \\ x_\ell
        \end{pmatrix}=\begin{pmatrix}
            r_\ell\\ \vdots \\ r_\ell
        \end{pmatrix}
    \end{equation}
has a solution in $\N$ if and only if   the equation $my-sx_\ell=r_\ell$ has a solution in $\N$. If $(y,s)$ is such a solution, then $$\left(x_1=\frac m{r_1}y, \dots ,  x_{\ell-1}=\frac m{r_{\ell-1}}y, x_\ell\right)$$ is a solution of \eqref{eq:equations1}. Therefore, \eqref{eq:equations1} has a solution in $\N$ if and only if $\gcd (m, s)$ divides $r _\ell$.
\end{Lemma}
\begin{Proof}
The claim follows because 
    \[\begin{pmatrix}
        1&0&\dots&0 \\-1&1&\dots &0\\ \vdots & &\ddots&\vdots\\
        -1&0& \dots &1
    \end{pmatrix}\begin{pmatrix}
        r_1&0&\dots&0 &-s\\0&r_2&\dots&0 &-s\\\vdots & &\ddots&&\vdots\\
        0&0& \dots &r_{\ell -1}&-s
    \end{pmatrix}=\begin{pmatrix}
        r_1&0&\dots&0 &-s\\-r_1&r_2&\dots&0 &0\\\vdots & &\ddots&&\vdots\\
        -r_1&0& \dots &r_{\ell -1}&0
    \end{pmatrix}\]
so the system \eqref{eq:equations1} is equivalent to the system
    \begin{equation*}\label{eq:equations2}
        \begin{pmatrix}
            r_1&0&\dots&0 &-s\\-r_1&r_2&\dots&0 &0\\\vdots & &\ddots&&\vdots\\
            -r_1&0& \dots &r_{\ell -1}&0
        \end{pmatrix}\begin{pmatrix}
            x_1\\ x_2\\\vdots \\ x_\ell
        \end{pmatrix}=\begin{pmatrix}
            r_\ell\\0\\ \vdots \\ 0
        \end{pmatrix}
    \end{equation*}
which has a solution in $\N$ if and only if   the equation $my-sx_\ell=r_\ell$ has a solution in $\N$. Also, it is clear that the solutions are related, as claimed in the statement.
\end{Proof}

\begin{Prop} \label{char}
Under the Notation~\ref{not}, the following are equivalent 
\begin{enumerate}
    \item[(i)] Every projective $\Lambda$-module is a direct sum of finitely generated modules.
    \item[(ii)] For every family of natural numbers $a_1,a_2,\dots,a_\ell$, $1 \leq a_i \leq t_i$, and for every $j\in \{1,\dots ,\ell\}$ and  $1 \leq b \neq a_j \leq t_j$, there are non-negative integers $x_1,x_2,\dots,x_\ell$ such that $r_{j,b} + x_{j}r_{j,a_j} = x_ir_{i,a_i}$ for every 
    $i \neq j$.
    \item[(iii)] For every family of natural numbers $a_1,a_2,\dots,a_\ell$, $1 \leq a_i \leq t_i$, and for every $j\in \{1,\dots ,\ell\}$ and  $1 \leq b \neq a_j \leq t_j$, $\gcd(\mathrm{lcm} (\{r_{1,a_1},\dots ,r_{\ell,a_\ell}\}\setminus \{r_{j,a_j}\}), r_{j,a_j})$ divides $r_{j,b}$.
\end{enumerate}
\end{Prop}
\begin{Proof}
Assume $(i)$ is true, so finitely generated projective modules of $\Lambda$ modulo a (minimal) element of $\mathcal T(\Lambda) \setminus \{0\}$ can be lifted to finitely generated projective $\Lambda$-modules.

By Remark~\ref{descriptiontracesalmostsemiperfect}, $\Lambda$ has an ideal  $I\in \mathcal{T} (\Lambda)$ such that $I_{\fm_i} = \Tr_{\Lambda_{\fm_i}}(U_{i,a_i})$ for each $i = 1,\dots,\ell$. As $\Lambda _Q$ is simple artinian and $R$ is $h$-local,  $\Lambda/I \cong \bigoplus_{i = 1}^{\ell} \Lambda_{\fm_i}/I_{\fm_i}$.  Therefore, every (finitely generated) projective $\Lambda/I$-module is of the form $\bigoplus_{i = 1}^{\ell} M_i$, where $M_i$ is a (finitely generated) projective $(\Lambda/I)_{\fm_i}$-module. Note that $(\Lambda/I)_{\fm_i} \cong \Lambda_{\fm_i}/I_{\fm_i}$, which is a semiperfect ring, so projective modules over it are described in Proposition~\ref{tracessemiperfect} (iv).

Consider a projective $\Lambda/I$-module $P'$ such that $P_{\fm_j}' \cong U_{j,b}/U_{j,b}I_{\fm_j}$ and $P_{\fm_{i}}' = 0$ if $i \neq j$. By assumption, there exists $P$  a finitely generated projective $\Lambda$-module such that $P/PI \cong P'$.  Then $P_{\fm_i} \cong U_{i,a_i}^{x_i}$ if $i \neq j$ and $P_{\fm_j} \cong U_{j,b} \oplus U_{j,a_j}^{x_j}$. Since $(P_{\fm_i})_Q \cong P_Q \cong (P_{\fm_j})_Q$, we obtain $x_ir_{i,a_i} = r_{j,b}+x_{j}r_{j,a_{j}}$. So $(ii)$ holds.

Conversely, let $I$ be a minimal non-zero trace of $\Lambda$. Then $I_{\fm_i} = \Tr_{\Lambda_{\fm_i}}(U_{i,a_i})$ for some $a_i$, $i \in \{1,\dots,\ell_i\}$. Recall that $\Lambda/I \cong \bigoplus_{i = 1}^{\ell} \Lambda_{\fm_i}/I_{\fm_i}$ is a semiperfect ring. Hence, it is sufficient to lift an indecomposable projective $\Lambda/I$-module to a finitely generated projective $\Lambda$-module. So let $P'$ be a finitely generated indecomposable projective $\Lambda/I$-module. Regarding the decomposition of $\Lambda/I$, $P'_{\fm_i} = 0$ for all $i$'s with exactly one exception. Assume $P_{\fm_j}' \cong U_{j,b}/U_{j,b}I_{\fm_{j}}$, where $b \neq a_j$. Take $x_1,\dots,x_{\ell} \in \N_0$ such that the equations in  $(ii)$ are satisfied. Let $r = r_{j,b} + x_jr_{j,a_j}$. 

For $i = 1,\dots,\ell$ set $U(\fm_i) = U_{i,a_i}^{x_i}$ if $i \neq j$ and $U(\fm_j) = U_{j,b} \oplus U_{j,a_j}^{x_j}$. For $\fn \in \mSpec(R) \setminus \{\fm_1,\dots,\fm_{\ell}\}$, let $U(\fn)$ be a projective $\Lambda_{\fn}$-module such that $U(\fn)_Q$ has length $r$ (it exists by Lemma~\ref{finitely_many}). Because of the choice of $x_1,\dots,x_{\ell}$, we can use Proposition~\ref{fgprojlocsem} to show the existence of a finitely generated projective $\Lambda$-module $P$ such that $P_{\fm} \cong U(\fm)$ for any maximal ideal $\fm$ of $R$. Then $P_{\fm_i}I_{\fm_i} = P_{\fm_i}$ if $i \neq j$ and $P_{\fm_j}/P_{\fm_j}I_{\fm_j} \cong U_{j,b}/U_{j,b}I_{\fm_j}$. Therefore $P/PI \cong \bigoplus_{i = 1}^{\ell} P_{\fm_i}/P_{\fm_i}I_{\fm_i}\cong P'$. This proves $(ii) \Rightarrow (i)$.

The equivalence of $(ii)$ and $(iii)$ follows from Lemma~\ref{equations}.
\end{Proof}

\begin{Ex} \label{exr1}
Under the assumptions in Notation~\ref{not}, let $\ell=2$, and assume that $U_{1,1},U_{1,2}$ and $U_{2,1},U_{2,2}$ are representatives of the finitely generated indecomposable projective right modules over $\Lambda_{\fm_1}$ and $\Lambda_{\fm_2}$, respectively. Suppose that 
    \[r_{1,1}=2,\quad r_{1,2}=4,\quad r_{2,1}=3,\quad r_{2,2}=9.\]

Since  $$12=r_{1,1} 2+r_{1,2} 2=r_{2,1} + r_{2,2},$$
 $\Lambda$ can be taken to be an algebra over an $h$-local domain $R$ with infinitely many maximal ideals, such that 
 \begin{itemize}
     \item[(i)] $\Lambda _Q\cong M_{12} (T_Q)$;
     \item[(ii)] $\Lambda _{\fm _1}$ is a semiperfect ring such that modulo its Jacobson radical is isomorphic to $M_2(T_{1,1})\times M_2(T_{1,2})$, and $\Lambda _{\fm _2}$ is a semiperfect ring that is a product of two division rings modulo the Jacobson radical;
     \item[(iii)]for any other maximal ideal $\fm$ of $R$ the localization  $\Lambda _\fm\cong M_{12} (T_\fm)$ where $T_\fm$ is a local ring.
\end{itemize}
By Examples~\ref{exs:locallysemiperfect}, such  a locally semiperfect  algebra $\Lambda$ exists and can be even taken to be a finitely generated $\Z$-algebra. 

By Proposition~\ref{systemfg}, to compute $\mathcal A=\Phi'(V(\Lambda))$ in $\N_0^2\times\N_0^2$  we need to find the $((x_{1,1},x_{1,2}),(x_{2,1},x_{2,2}))\in\N_0^2\times \N_0^2$ which are   solutions of the system of equations:
    \begin{align}\label{eq:example}
    \begin{split}
        2x_{1,1}+4x_{1,2}&=r \\
        3x_{2,1}+9x_{2,2}&=r,
    \end{split}
    \end{align}
for some fixed $r\in\mathbb N_0$. 

Since $\mathrm{lcm}(\gcd(2,4),\gcd(3,9))=6\mid r$, the image of the rank function $\Psi$ lies in $6\N_0$ and, since there are finitely generated projective modules of rank $6$, in fact $\Psi (V(\Lambda))= 6\N_0$. 

As each $r_{i,j}$ is coprime to $r_{i',j'}$ for every $i'\neq i$ and  $j'\neq j$,  statement (iii) in Proposition~\ref{char} is satisfied. So every projective $\Lambda$-module is a direct sum of finitely generated modules.

Denote an element in $\mathcal A$ by $((x_{1,1},x_{1,2}),(x_{2,1},x_{2,2}))_r$, where $x_{1,1},x_{1,2},x_{2,1},x_{2,2},r\in\N_0$ satisfy \eqref{eq:example}. Hence, a minimal set of generators of $\mathcal A$ is
    \begin{gather*}
        ((3,0),(2,0))_6,\quad ((1,1),(2,0))_6,\quad ((0,3),(4,0))_{12},\\ 
        ((6,0),(1,1))_{12},\quad ((4,1),(1,1))_{12},\quad ((2,2),(1,1))_{12},\quad ((0,3),(1,1))_{12},\\ 
        ((9,0),(0,2))_{18},\quad ((7,1),(0,2))_{18},\quad ((5,2),(0,2))_{18},\\
        ((3,3),(0,2))_{18},\quad ((1,4),(0,2))_{18},\quad ((0,9),(0,4))_{36}.
    \end{gather*}
So the isomorphism class of  finitely generated projective $\Lambda$-module has as an image in $\mathcal A$ a non-negative integral linear combination of elements in this generating set.

Recall that if $\Lambda$ is taken to be finitely generated over a noetherian domain $R$ then, by Proposition~\ref{genusisomorphic_semilocal}, $\mathcal{A}\cong V(\Lambda)$.

If we extend the monoid morphism $\Phi'$ to $V^*(\Lambda)$ and enlarge $\N_0$ to $\N_0^*=\N_0\cup\{\infty\}$, then we can also compute $\mathcal B=\Phi'(B(\Lambda))$ in $(\N_0^*)^2\times(\N_0^*)^2$. A minimal set of generators of $\mathcal B$ is
    \begin{gather*}
        ((\infty,0),(\infty,0)),\quad ((\infty,0),(0,\infty)),\quad ((0,\infty),(\infty,0)),\quad ((0,\infty),(0,\infty)),\\ 
        ((\infty,1),(\infty,0)),\quad ((\infty,1),(0,\infty)),\quad ((1,\infty),(\infty,0)),\quad ((1,\infty),(0,\infty)),\\ 
        ((\infty,0),(\infty,1)),\quad ((\infty,0),(1,\infty)),\quad ((0,\infty),(\infty,1)),\quad ((0,\infty),(1,\infty)).
    \end{gather*}
and, as explained before, all of them correspond to the isomorphism class of a countably generated projective module that is a direct sum of finitely generated modules.  For example,
$$((\infty,1),(\infty,0))= ((1,1),(2,0))_6+\infty \cdot ((3,0),(2,0))_6.$$
\end{Ex}

\begin{Ex}
Under the assumptions in Notation~\ref{not}, there are examples of algebras $\Lambda$ where there is a projective $\Lambda$-module which is not a direct sum of finitely generated modules. Set now
    \[r_{1,1}=2,\quad r_{1,2}=8,\quad r_{2,1}=2,\quad r_{2,2}=4.\]
By Examples~\ref{exs:locallysemiperfect}, we know there exist  locally semiperfect algebras $\Lambda$ with such invariants.

To determine $\mathcal A=\Phi'(V(\Lambda))$ we have to solve the following system of linear equations:
    \begin{align}\label{eq:example2}
    \begin{split}
        2x_{1,1}+8x_{1,2}&=r \\
        2x_{2,1}+4x_{2,2}&=r.
    \end{split}
    \end{align}
Since $\gcd(\{r_{1,2},r_{2,2}\}\setminus\{r_{1,2}\},r_{1,2})=\gcd(4,8)=4\nmid 2=r_{1,1}$, Proposition~\ref{char} (iii) fails. Therefore, in this case, there are projective $\Lambda$-modules which are not a direct sum of finitely generated modules. We are going to compute an explicit example of such a module.

As in Example~\ref{exr1}, if we extend the monoid morphism $\Phi'$ to $V^*(\Lambda)$ and enlarge $\N_0$ to $\N_0^*=\N_0\cup\{\infty\}$, then we can compute $\mathcal B=\Phi'(B(\Lambda))$ in $(\N_0^*)^2\times(\N_0^*)^2$. The generating set is the same as in Example~\ref{exr1}, because it is independent of the $r_{i,j}$'s chosen, it only depends on how many of them there are. 

Let $[P]\in B(\Lambda)$ be such that $\Phi'([P])=((0,\infty),(1,\infty))$. Then there is no $\mathbf{X}\in\N_0^2\times\N_0^2$ of the form $\mathbf{X}=((0,x_{1,2}),(1,x_{2,2}))$ being a solution of \eqref{eq:example2} since $8x_{1,2}\equiv0\mod4$ and $2+4x_{2,2}\equiv2\mod4$. Therefore, there is no finitely generated projective $\Lambda$-module $P'$ such that $\Phi'([P'])=((0,x_{1,2}),(1,x_{2,2}))$, so $P$ is a countably generated projective $\Lambda$-module which is not a direct sum of finitely generated modules.

Let $P_1$ be a finitely generated $\Lambda$-module such that $\Phi'([P_1])=((0,1),(0,2))_8$, and let $P_2$ be such that $\Phi'([P_2])=((0,1),(2,1))_8$. By Proposition~\ref{genusisomorphic}, $P\cong P\oplus P_1^{(\omega)}$ and $P^2\cong P_2\oplus P_1^{(\omega)}$. Thus $P^2$ is a direct sum of finitely generated projective modules.
\end{Ex}

\section{Torsion-free modules over \texorpdfstring{$h$}{h}-local domains} \label{sec:9}

In this section we will use extensively the notation  introduced in \S ~\ref{subsec:indnoetherian} and Proposition~\ref{equivalencia}. 
We apply the previous results to $\Add(M)$ for $M$ a finitely generated torsion-free module $M$ over an $h$-local domain. 

We need to keep in mind the following lemma.

\begin{Lemma} \label{endfg}
\cite[Lemma 2.18]{AHP}
Let $R$ be a commutative domain with field of fractions $Q$. Let $M$ be a non-zero finitely generated torsion-free module with endomorphism ring $\Lambda$. Then
\begin{itemize}
    \item[(i)] $\Lambda$ is a torsion-free $R$-module containing $R$.
    \item[(ii)] $\Lambda_\fm\cong \End_{R_\fm}(M_\fm)$ is a semilocal ring for every maximal ideal $\fm$ of $R$.
    \item[(iii)] $\Lambda_Q\cong M_k(Q)$ is a simple artinian ring, where $k$ is the rank of $M$.
    \item[(iv)] For any non-zero two-sided ideal $I$ of $\Lambda$, $I\cap R\neq \{0\}$.
\end{itemize}
If, in addition, $R$ is $h$-local, then 
\begin{itemize}
    \item[(v)] $M_\fm$ is free for almost all maximal ideals of $R$.
    
    \item[(vi)] For any non-zero two-sided ideal $I$ of $\Lambda$, $I\neq \Lambda$, there is a canonical isomorphism $\Lambda/I\longrightarrow (\Lambda/I)_{\fm_1}\times\dots\times (\Lambda/I)_{\fm_t}$,
    where $\{\fm_1,\dotsc,\fm_t\}$ are the unique maximal ideals of $R$ containing $I\cap R$.
\end{itemize}
\end{Lemma}

Let $M_R$ be a finitely generated torsion-free module of rank $k$ over an $h$-local domain $R$. Following Lemma~\ref{endfg}, let $\fm_1,\dots,\fm_{\ell}$ be the list of those maximal ideals of $R$ over which $M_{\fm}$ is not free. Assume further that $\Lambda$ is locally semiperfect. That is, each $M_{\fm}$ is a direct sum of modules with local endomorphism rings, cf. Example~\ref{ex:locsemi}. For $i = 1,\dots,\ell$, let $M_{\fm_i} \cong \bigoplus_{j = 1}^{t_i} M_{i,j}^{m_{i,j}}$, where $M_{i,1},\dots,M_{i,t_i}$ are pairwise non-isomorphic indecomposable modules and $m_{i,j} \in \N$. By Proposition~\ref{equivalencia}, $\Lambda_{\fm_{i}} \cong \bigoplus_{j = 1}^{t_i} U_{i,j}^{m_{i,j}}$, where $U_{i,j}  = \Hom_{R_{\fm_i}}(M_{\fm_i},M_{i,j})$ is an indecomposable decomposition of $\Lambda_{\fm_i}$. Let $r_{i,j}$ be the length of $(U_{i,j})_Q$. As a $Q$-module, $(U_{i,j})_Q  = Q^{kr_{i,j}}$. On the other hand, by Lemma~\ref{isofg},  $(U_{i,j})_Q \cong \Hom_Q(M_Q,(M_{i,j})_Q) \cong (M_{i,j})_Q^k = Q^{k\rank (M_{i,j})}$. Therefore, $r_{i,j} = \rank (M_{i,j})$. 

\begin{Lemma} \label{generator}
Let $R$ be an $h$-local domain of Krull dimension $1$, and let $M$ be a finitely generated torsion-free $R$-module which is a generator of $\Mod R$, i.e., $\Tr_R(M) = R$. Assume that every module of $\Add(M)$ is a direct sum of finitely generated modules. Then
  \begin{enumerate}
    \item[(i)] $\End_R(M)$ is a locally semiperfect ring.
    \item[(ii)] for every $\fm \in \mSpec (R)$, $M_{\fm}$ contains a direct summand isomorphic to $R_{\fm}$.
  \end{enumerate}
\end{Lemma}
\begin{Proof}
The assumption $\Tr_R(M) = R$ is equivalent to $R \in \add(M)$. In particular, $\Add(M)$ contains $R^{(\omega)} \oplus M$, so $(i)$ follows from \cite[Proposition~6.1]{AHP}. By Lemma~\ref{loc.traces}, $\Tr_{R_{\fm}}(M_{\fm}) = R_{\fm}$ for any maximal ideal $\fm \in \mSpec(R)$, hence $R_{\fm} \in \add(M_{\fm})$. Since $R_{\fm}$ is local, $R_{\fm}$ has to be isomorphic to a direct summand of $M_{\fm}$ (otherwise, $\Tr_{R_{\fm}}(M_{\fm}) \subseteq \fm R_{\fm}$).   
\end{Proof}

Now we are ready to prove our main application to direct sum decompositions of modules, which complements the results in \cite{AHP}. 

\begin{Th} \label{th:Add}
Let $R$ be an $h$-local domain,  and let $M$ be a non-zero  finitely generated torsion-free 
$R$-module with endomorphism ring $\Lambda = \End_R(M)$. Consider the following statements,
\begin{enumerate}
    \item[(i)] $\Lambda$ is locally semiperfect and if $X$ and $Y$ are indecomposable direct summands of $M_{\fm}$ and $M_{\fn}$, respectively, for some $\fm \neq \fn \in \mSpec (R)$. Then $\rank (X)$ and $\rank (Y)$ are coprime.
    \item[(ii)] Every module in $\Add(M)$ is a direct sum of finitely generated modules. 
\end{enumerate}
Then $(i)$ always implies $(ii)$ and the converse is true provided $M$ is a generator and $R$ has Krull dimension $1$.
\end{Th}
\begin{Proof} 
Notice that by Lemma~\ref{endfg}, $\Lambda$ is torsion-free as an $R$-module, and $\Lambda _Q$ is simple artinian. 

Let us show that  $(i) \Rightarrow(ii)$.  Assume that $(i)$ is true. We will use the notation introduced after Lemma~\ref{endfg}.

Note that  $r_{i,a}$ and $r_{j,b}$ are coprime whenever $i \neq j$, $a \in  \{1,\dots,t_i\}$, $b \in \{1,\dots,t_j\}$. We want to check that $(ii)$ of Proposition~\ref{char} is true. 

Consider $a_1,\dots,a_{\ell}$, $a_i \in \{1,\dots,t_{i}\}$, and $A = \prod_{i = 1}^{\ell} r_{i,a_i}$. Fix $1 \leq j \leq \ell$. Since $r_{1,a_1},\dots,r_{\ell,a_{\ell}}$ are pairwise coprime, there exists $t \in \{0,\dots,A-1\}$ such that $t \equiv 0\ ({\rm mod}\ r_{i,a_i})$ if $i \neq j$ and $t \equiv r_{j,b}\ ({\rm mod}\ r_{j,a_j})$. If $t< r_{j,b}$ we can increase $t$ by a suitable multiple of $A$, so we may assume $t \geq r_{j,b}$. Set $x_i:= t/r_{i,a_i}$ if $i \neq j$ and $x_j:=(t - r_{j,b})/r_{j,a_j}$. Then $x_1,\dots,x_{\ell}$ are non-negative integers and $t = x_i r_{i,a_i}$ if $ i\neq j$ and $t = x_j r_{j,a_j}+r_{j,b}$. In other words, the condition $(ii)$ of Proposition~\ref{char} is true. Then any projective $\Lambda$-module is a direct sum of finitely generated modules. $(ii)$ then follows from Proposition~\ref{equivalencia}.

Now, assume that $(ii)$ is true, and suppose that $\Tr_R(M) = R$, and $R$ has Krull dimension $1$. Then, by Lemma~\ref{generator}, $\Lambda$ is locally semiperfect. Keeping the notation of Proposition~\ref{char}, we have to prove that $r_{i,a}$ and $r_{j,b}$ are coprime whenever $i \neq j$, $a \in \{1,\dots,t_i\}$, and $b \in \{1,\dots,t_j\}$. Note that it is trivial if either $r_{i,a} = 1$ or $r_{j,b} = 1$. By Lemma~\ref{generator}, for any $\fm \in \mSpec (R)$, $R_\fm$ is a direct summand of $M_\fm$ (so $M_\fm$ has a direct summand of rank $1$, cf. the remarks before Lemma~\ref{generator}). Therefore, if $r_{i,a} > 1$, then there exists $c \in \{1,\dots,t_i\}$ such that $r_{i,c} = 1$. Then $r_{i,c} + x r_{i,a} = y r_{j,b}$ is solvable only if $r_{i,a},r_{j,b}$ are coprime. Now $(i)$ follows from Proposition~\ref{char}.
\end{Proof}

\bibliographystyle{amsplain}
\bibliography{references}

\begin{bibdiv}
\begin{biblist}

\bib{AHP}{article}{
      author={\'Alvarez, R.},
      author={Herbera, D.},
      author={Příhoda, P.},
       title={Torsion-free modules over commutative domains of {K}rull
  dimension one},
        date={2025},
     journal={Revista Matemática Iberoamericana},
       pages={On\ndash line first},
}

\bib{andersonfuller}{book}{
      author={Anderson, F.~W.},
      author={Fuller, K.~R.},
       title={Rings and categories of modules},
      series={Graduate Texts in Mathematics {13}},
   publisher={Springer New York},
        date={1992},
}

\bib{AABPV}{article}{
      author={Antoine, Ramon},
      author={Ara, Pere},
      author={Bosa, Joan},
      author={Perera, Francesc},
      author={Vilalta, Eduard},
       title={The {C}untz semigroup of a ring},
        date={2025},
        ISSN={1022-1824,1420-9020},
     journal={Selecta Math. (N.S.)},
      volume={31},
      number={1},
       pages={Paper No. 6, 46},
         url={https://doi.org/10.1007/s00029-024-01002-9},
}

\bib{bass}{article}{
      author={Bass, H.},
       title={Big projective modules are free},
        date={1963},
     journal={Illinois Journal of Mathematics},
      volume={7},
      number={1},
       pages={24–31},
}

\bib{libro}{book}{
      author={Facchini, A.},
       title={Module theory: Endomorphism rings and direct sum decompositions
  in some classes of modules},
      series={Progress in Mathematics {167}},
   publisher={Birkh\"{a}user Basel},
        date={1998},
}

\bib{FHS1}{article}{
      author={Facchini, A.},
      author={Herbera, D.},
      author={Sakhajev, I.},
       title={Finitely generated flat modules and a characterization of
  semiperfect rings},
        date={2003},
     journal={Communications in Algebra},
      volume={31},
      number={9},
       pages={4195–4214},
}

\bib{FHL}{book}{
      author={Fontana, M.},
      author={Houston, E.},
      author={Lucas, T.},
       title={Factoring ideals in integral domains},
      series={Lecture Notes of the Unione Matematica Italiana {14}},
   publisher={Springer Berlin Heidelberg},
        date={2013},
}

\bib{fuchssalce}{book}{
      author={Fuchs, L.},
      author={Salce, L.},
       title={Modules over non-{N}oetherian domains},
      series={Mathematical Surveys and Monographs {84}},
   publisher={American Mathematical Society},
        date={2000},
}

\bib{GW}{book}{
      author={Goodearl, K.~R.},
      author={Warfield, R.~B., Jr.},
       title={An introduction to noncommutative {N}oetherian rings},
     edition={Second Edition},
      series={London Mathematical Society Student Texts},
   publisher={Cambridge University Press, Cambridge},
        date={2004},
      volume={61},
        ISBN={0-521-83687-5; 0-521-54537-4},
         url={https://doi.org/10.1017/CBO9780511841699},
}

\bib{crelle}{article}{
      author={Herbera, D.},
      author={Příhoda, P.},
       title={Big projective modules over {N}oetherian semilocal rings},
        date={2010},
     journal={Journal f\"{u}r die reine und angewandte Mathematik (Crelle's
  Journal)},
      volume={648},
       pages={111–148},
}

\bib{TAMS}{article}{
      author={Herbera, D.},
      author={Příhoda, P.},
       title={Infinitely generated projective modules over pullbacks of rings},
        date={2014},
     journal={Transactions of the American Mathematical Society},
      volume={366},
      number={3},
       pages={1433–1454},
}

\bib{traces}{article}{
      author={Herbera, D.},
      author={Příhoda, P.},
       title={Reconstructing projective modules from its trace ideal},
        date={2014},
     journal={Journal of Algebra},
      volume={416},
       pages={25–57},
}

\bib{wiegand}{article}{
      author={Herbera, D.},
      author={Příhoda, P.},
      author={Wiegand, R.},
       title={Big pure projective modules over commutative {N}oetherian rings:
  Comparison with the completion},
        date={2025},
     journal={Forum Mathematicum},
      volume={37},
      number={4},
       pages={1103\ndash 1146},
}

\bib{hinohara}{article}{
      author={Hinohara, Y.},
       title={Supplement to ``{P}rojective modules over weakly {N}oetherian
  rings''},
        date={1963},
     journal={Journal of the Mathematical Society of Japan},
      volume={15},
      number={1},
       pages={474\ndash 475},
}

\bib{jaffard}{article}{
      author={Jaffard, P.},
       title={Théorie arithmétique des anneaux du type de {D}edekind},
        date={1952},
     journal={Bulletin de la Société Mathématique de France},
      volume={80},
       pages={61–100},
}

\bib{kaplansky}{article}{
      author={Kaplansky, I.},
       title={Projective modules},
        date={1958},
     journal={Annals of Mathematics},
      volume={68},
      number={2},
       pages={372–377},
}

\bib{lam2}{book}{
      author={Lam, T.~Y.},
       title={Lectures on modules and rings},
      series={Graduate Texts in Mathematics {189}},
   publisher={Springer New York},
        date={1999},
}

\bib{lazard}{article}{
      author={Lazard, D.},
       title={Liberté des gros modules projectifs},
        date={1974},
     journal={Journal of Algebra},
      volume={31},
      number={3},
       pages={437–451},
}

\bib{levyodenthal2}{article}{
      author={Levy, L.~S.},
      author={Odenthal, C.~J.},
       title={Package deal theorems and splitting orders in dimension 1},
        date={1996},
     journal={Transactions of the American Mathematical Society},
      volume={348},
      number={9},
       pages={3457–3503},
}

\bib{matlis3}{book}{
      author={Matlis, E.},
       title={Torsion-free modules},
      series={Lectures in Mathematics},
   publisher={University of Chicago Press},
        date={1973},
}

\bib{MR}{book}{
      author={McConnell, J.},
      author={Robson, J.},
       title={Noncommutative {N}oetherian rings},
      series={Pure and Applied Mathematics (New York)},
   publisher={John Wiley \& Sons, Ltd., Chichester},
        date={1987},
}

\bib{nazemian_smertnig}{article}{
      author={Nazemian, Z.},
      author={Smertnig, D.},
       title={A monoid-theoretical approach to infinite direct-sum
  decompositions of modules},
        date={2024},
}

\bib{olberding}{inproceedings}{
      author={Olberding, Bruce},
       title={Characterizations and constructions of {$h$}-local domains},
        date={2008},
   booktitle={Models, modules and abelian groups},
   publisher={Walter de Gruyter, Berlin},
       pages={385\ndash 406},
         url={https://doi.org/10.1515/9783110203035.385},
}

\bib{PS}{article}{
      author={Procesi, C.},
      author={Small, L.~W.},
       title={Endomorphism rings of modules over {PI}-algebras},
        date={1968},
     journal={Mathematische Zeitschrift},
      volume={106},
      number={3},
       pages={178–180},
}

\bib{P2}{article}{
      author={Příhoda, P.},
       title={Fair-sized projective modules},
        date={2010},
     journal={Rendiconti del Seminario Matematico della Universit\`a di
  Padova},
      volume={123},
       pages={141–168},
}

\bib{P3}{article}{
      author={Příhoda, P.},
       title={Generalized lattices over one-dimensional {N}oetherian domains},
        date={2022},
     journal={Journal of Commutative Algebra},
      volume={14},
      number={3},
       pages={443–453},
}

\bib{PP}{article}{
      author={Příhoda, P.},
      author={Puninski, G.},
       title={Iterated power intersections of ideals in rings of iterated
  differential polynomials},
        date={2013},
     journal={Journal of Algebra and Its Applications},
      volume={12},
      number={07},
       pages={1350020, 10},
}

\bib{SW}{article}{
      author={Sexauer, N.~E.},
      author={Warnock, J.~E.},
       title={The radical of the row-finite matrices over an arbitrary ring},
        date={1969},
     journal={Transactions of the American Mathematical Society},
      volume={139},
       pages={287–295},
}

\bib{small}{article}{
      author={Small, L.~W.},
       title={Localization in {PI}-rings},
        date={1971},
     journal={Journal of Algebra},
      volume={18},
      number={2},
       pages={269–270},
}

\bib{RS}{article}{
      author={Small, L.~W.},
      author={Robson, J.~C.},
       title={Idempotent ideals in {P.I.} rings},
        date={1976},
     journal={Journal of the London Mathematical Society},
      volume={14},
      number={1},
       pages={120–122},
}

\bib{Vasconcelos2}{article}{
      author={Vasconcelos, W.~V.},
       title={On projective modules of finite rank},
        date={1969},
     journal={Proceedings of the American Mathematical Society},
      volume={22},
      number={2},
       pages={430\ndash 433},
}

\bib{whitehead}{article}{
      author={Whitehead, James~M.},
       title={Projective modules and their trace ideals},
        date={1980},
     journal={Communications in Algebra},
      volume={8},
      number={19},
       pages={1873–1901},
}

\end{biblist}
\end{bibdiv}

\end{document}